\documentclass[11pt, reqno]{amsart}
\pdfoutput=1
\usepackage{amssymb,amsmath,amsthm,amsfonts, upgreek,latexsym}
\usepackage{epsfig,esint}
\usepackage{mathrsfs}
\usepackage{mathtools,textcomp}
\usepackage{graphicx}
\usepackage{epic,eepic,wrapfig,color, ifthen} 
\usepackage{url,cite}
\usepackage[colorlinks=true, citecolor=blue, filecolor=cyan, linkcolor=red, urlcolor=magenta]{hyperref}
\usepackage{pmboxdraw}
\usepackage{verbatim}
\usepackage[normalem]{ulem}
\usepackage{enumerate}
\usepackage[bottom]{footmisc}

\usepackage{imakeidx}

\newtheorem{thm}{Theorem}[section]
\newtheorem{theorem}[thm]{Theorem}
\newtheorem{lemma}[thm]{Lemma}

\newtheorem{prop}[thm]{Proposition}
\newtheorem{corollary}[thm]{Corollary}
\newtheorem{defn}[thm]{Definition}
\newtheorem{example}[thm]{Example}
\def\As{{\bf (A)}}
\def\diam{{\text{\rm diam}}}

\numberwithin{equation}{section}
\allowdisplaybreaks

\def\rX {\overline X}

\DeclareFontFamily{U}{mathx}{}
\DeclareFontShape{U}{mathx}{m}{n}{<-> mathx10}{}
\DeclareSymbolFont{mathx}{U}{mathx}{m}{n}
\DeclareMathAccent{\wc}{0}{mathx}{"71}

\def\1{{\bf 1}}
\def\ind{{\bf 1}}
\def\nn{\nonumber}
\def\bee{\begin{equation}}
\def\eee{\end{equation}}
\def\sA {{\mathcal A}} \def\sB {{\mathcal B}} \def\sC {{\mathcal C}}
 \def\sE {{\mathcal E}} \def\sF {{\mathcal F}}

 \def\sN {{\mathcal N}} \def\sO {{\mathcal O}}

 \def\bH {{\mathbb H}}

\def\R {{\mathbb R}}  \def\bH {{\mathbb H}}
 \def\Z {{\mathbb Z}}

\def\esssup{{\mathop{\rm ess \; sup \, }}}

\newcommand{\overbar}[1]{\mkern 0.5mu\overline{\mkern-0.5mu#1\mkern-0.5mu}\mkern 0.5mu}

\def\eps{\varepsilon}
\def\e{\mathbf e}

\def\FF{{\mathcal F}}

\def\R{{\mathbb R}}

\def\E{{\mathbb E}}
\def\F{{\bf F}}
\def\P{{\mathbb P}}
\def\N{{\mathbb N}}
\def\eps{\varepsilon}
\def\wh{\widehat}
\def\wt{\widetilde}

\def\ub{{\overbar \beta}}

\begin{document}
	\title[]
{ \bf   
Markov processes with jump kernels blowing up  at the boundary of $C^{1, \eta}$ open sets
}

\author{Soobin Cho, \quad Panki Kim, \quad Renming Song \quad and \quad Zoran Vondra\v{c}ek}

\address[Cho]{School of Mathematics and Statistics, Beijing Institute of Technology, Beijing 100081, China}
\curraddr{}
\email{soobincho.math@gmail.com}

\address[Kim]{Department of Mathematical Sciences and Research Institute of Mathematics,
	Seoul National University,	Seoul 08826, Republic of Korea}
\curraddr{}
\thanks{Panki Kim is supported by the National Research Foundation of Korea(NRF) grant funded by the Korea government(MSIP) (No. RS-2026-25518141).}
\email{pkim@snu.ac.kr}

\address[Song]{
	Department of Mathematics, University of Illinois Urbana-Champaign, Urbana, IL 61801,
	USA}
\curraddr{}
\thanks{Research of Renming Song is supported in part by a grant from
	the Simons Foundation \#960480.}
\email{rsong@illinois.edu}

\address[Vondra\v{c}ek]
{
	Dr.~Franjo Tu{\dj}man Defense and Security 
	University
	 and Department of Mathematics, Faculty of Science, University of Zagreb, Zagreb, Croatia
}
\curraddr{}
\thanks{Zoran Vondra\v{c}ek is supported in part by the Croatian Science Foundation under the project IP-2025-02-8793.} 
\email{vondra@math.hr}

\date{}

\begin{abstract}
We study the boundary behavior of nonnegative harmonic functions for a class of 
purely discontinuous symmetric Markov processes in $C^{1,\eta}$ open sets $D\subset\mathbb R^d$. The jump kernels are allowed to blow up at the boundary of $D$ and 
the processes may have 
a critical killing potential. 
We identify the exact boundary decay exponent $p$, characterized by the coefficient of the critical killing potential. 
We establish a boundary Harnack principle which yields the comparability of $f(x)/\delta_D(x)^p$ for nonnegative harmonic functions that vanish continuously on a portion of the boundary. 
The key ingredient of the proof combines estimates of the exact decay rate for the survival probability of the killed process with a modified box method and uses a regularized distance function to handle the limited boundary regularity. 
\end{abstract}

\bigskip

\maketitle



\bigskip
\tableofcontents

%%%%%%%%%%%%%%%%%%%%%%%%%%%%%%%%%%%%%%%%%%%%%%%%%%%%%%%%%%%%%%%%%%%%%%%%%%%%%%%%%%%%%%%%%%%%%%%%%%%%%%%%%%%%%%%%%%%%%%%%%%%%%%%%%%%%%
%%%%%%%%%%%%%%%%%                      Introduction and main results                    %%%%%%%%%%%%%%%%%%%%%%%%%%%%%%%%%%%%%%%%%%%%%
%%%%%%%%%%%%%%%%%%%%%%%%%%%%%%%%%%%%%%%%%%%%%%%%%%%%%%%%%%%%%%%%%%%%%%%%%%%%%%%%%%%%%%%%%%%%%%%%%%%%%%%%%%%%%%%%%%%%%%%%%%%%%%%%%%%%%
\section{Introduction and main results}\label{intro}

Over the last few decades, substantial progress has been made in the potential theory of stable-like processes in open subsets $D$ of $\mathbb{R}^d$. Roughly speaking, stable-like processes are 
purely discontinuous 
Markov processes whose jump intensities are governed by kernels comparable to 
the jump intensity $|x-y|^{-d-\alpha}$ of isotropic $\alpha$-stable processes for some $\alpha\in(0,2)$. 
Among the fundamental results established in this setting are the boundary Harnack principle \cite{B97, BBC03, BKK08, CaSi07, SW99}, sharp two-sided estimates of Green functions 
\cite{CK02, CS98, K97}, and sharp two-sided estimates of Dirichlet heat kernels 
\cite{BGR10, CKS10, CKS10a, CT11}. 
 Related developments in the regularity theory of linear and nonlinear integro-differential elliptic equations can be found in \cite{FR} and the references therein.

A common feature of all these works is that the jump kernel of the process in $D$ is either the restriction to $D\times D$ of the jump kernel of an underlying process in $\mathbb{R}^d$ or is comparable to such a restriction. 
Consequently, the associated integro-differential operators are non-degenerate and 
satisfy a uniform ellipticity condition: their jump kernels are bounded above and below by constant multiples of the corresponding  
stable kernel $|x-y|^{-d-\alpha}$.

An important source of examples of degenerate 
operators comes from subordination. In particular, the generators of subordinate killed Brownian motions and subordinate killed stable processes in an open set $D$ are 
operators obtained by subordinating the corresponding killed processes with stable subordinators. The generator of a subordinate killed Brownian motion is the spectral fractional Laplacian, which has been extensively studied in the PDE community; see, for instance, \cite{BFV18, BV16, CaSt16, Grubb, StTo10}. A characteristic feature of these processes is that their jump kernels admit the representation
$$
J(x,y)=\sB(x,y)|x-y|^{-d-\alpha},
$$
where $\sB(x,y)$ 
is bounded between two positive constants in any compact subset of $D\times D$, 
but tends to zero as
$x$ or $y$ approaches the boundary of $D$; see \cite{KSV19, KSV20, SV03}. Thus, the corresponding jump intensities 
are comparable to $|x-y|^{-d-\alpha}$ away from the boundary and degenerate near the boundary.

Motivated in part by these examples, we developed a general framework for the potential theory of purely discontinuous Markov processes with jump kernels of the above form, where $\sB$ decays to zero at the boundary. This theory was first developed in the half-space
in \cite{KSV23a}, and subsequently extended to general smooth open sets in \cite{CKSV23}. 
The potential theory of these degenerate processes has been investigated in a series of papers \cite{CKSV22, CKSV23, CKSV25b, KSV23a, KSV24, KSV24b}. More specifically, \cite{CKSV22, KSV23a, KSV24, KSV24b} address elliptic aspects of the theory, including the boundary Harnack principle and Green function estimates, whereas \cite{CKSV23, CKSV25b} establish sharp two-sided heat kernel estimates.

Boundary degeneracy, however, can also arise in the opposite regime, when the factor $\sB(x,y)$ blows up rather than vanishes near the boundary. Important examples include the jump kernels of the trace of 
a stable process on a smooth open set (see, for instance \cite{BGPKR}) 
and of the process associated with a nonlocal Neumann problem introduced in \cite{DRV17, Von21}.
 The latter has also been investigated from an analytic perspective, with various extensions of the nonlocal Neumann problem developed in \cite{Ab20, AFR20, MPL, DRVW26}. 
 In both cases, the jump kernel can be expressed in the form
$$
J(x,y)=\sB(x,y)|x-y|^{-d-\alpha},
$$
with $\sB$ blowing up near the boundary. 
Study of the potential theory of this class of processes was initiated in \cite{KSV23}, 
where fundamental interior potential-theoretic properties were established. Subsequently, \cite{KSV25} treated the half-space and established the boundary Harnack principle and sharp Green function estimates.

A significant advance in this direction was made in \cite{CK25+}. In a fairly general setting, this work established the stability of Hölder regularity for caloric and harmonic functions, introduced a new formulation of a tail condition for jump measures, and developed associated weighted functional inequalities. Building on these results, \cite{CKSV25} established sharp two-sided heat kernel estimates for a large class of purely discontinuous conservative symmetric Markov processes on closed subsets $F$ of $\R^d$. The jump kernels considered there have the form
$$
J(x,y)=\sB(x,y)|x-y|^{-d-\alpha},
$$
where $\sB$ blows up on a Borel subset $\Sigma$ of $F$, while $F\setminus\Sigma$ is a dense $\kappa$-fat subset of $F$. A fundamental difficulty in this setting is that the tails of the jump measures are not uniformly bounded. In particular, the parabolic Harnack principle fails, and the standard methods of heat kernel analysis do not provide the necessary a priori off-diagonal estimates.

The process $\overline X$ constructed in \cite{CKSV25} is conservative, that is, it has no killing. In the present paper, we investigate the potential-theoretic and analytic properties of two subprocesses of $\overline X$ in a $C^{1,\eta}$ open set $D\subset\R^d$: the process obtained by killing $\overline X$ upon approaching the boundary $\partial D$, and the process obtained by introducing a critical killing potential.

The analysis of these processes presents several new challenges. The boundary blow-up of the jump intensity leads to highly nonuniform jump measures, while the introduction of killing further complicates the behavior of harmonic functions near the boundary. These difficulties are particularly pronounced in domains of limited regularity, where the construction of suitable barrier functions and the derivation of precise boundary decay estimates require new ideas.

A key ingredient in the proof is the use of the regularized distance introduced in \cite{Li85} and more recently employed in \cite{CS25}, rather than the ordinary distance to the boundary, in the construction of barrier functions. This choice is particularly important for obtaining the precise boundary decay rate under the relatively weak regularity assumption on $D$. For comparison, in the nondegenerate setting, where $\sB(x,y)$ is bounded above and below by positive constants, the precise boundary decay rate of nonnegative harmonic functions vanishing on the boundary of a $C^{1,\eta}$ open set was recently investigated in \cite{KW24}.

\medskip

Let us now describe the setting and introduce the assumptions. 
\emph{Throughout the paper we assume that 
$d\ge 1$ and 
$\alpha\in (0,2)$.}

We say that $\Phi:[0,\infty) \to [1,\infty)$ is a \emph{blow-up weight function} if it is continuous, increasing, 
satisfies  $\Phi(0)=1$ and has the following upper scaling property: there exist constants $\beta\ge 0$ and $C>1$ such that
\begin{align}\label{e:upper-scaling}
	\Phi(r)/\Phi(s) \le C(r/s)^{ \beta} \quad \text{for all $0<s\le r$}.
\end{align}
The upper  Matuszewska index $\overline \beta$ of $\Phi$ is defined as  the infimum of all $ \beta$ 
for which \eqref{e:upper-scaling} holds.

Suppose that $D\subset \R^d$ is a Lipschitz open set with localization radius $R_0$, and that the Assouad dimension ${\rm dim}_{\rm A}(\partial D)$ of $\partial D$ is $d-1$.  
We refer to \cite[Chapter 2]{Fr20} for definitions and  equivalent formulations of the Assouad dimension. 
Note that
 if $D$ is a bounded Lipschitz open set, or a half-space-like 
Lipschitz open set, or an exterior Lipschitz open set, then ${\rm dim}_{\rm A}(\partial D)=d-1$. 
However, there exist unbounded smooth open sets $D$ such that ${\rm dim}_{\rm A}(\partial D)>d-1$. Examples include $D=\R^d\setminus \bigcup_{k\in \Z^d}\overline{B(k,1/4)}$ and  $D=\{x\in \R^d: x_d\in \cup^{\infty}_{n=0}(2n, 2n+1)\}$, for which ${\rm dim}_{\rm A}(\partial D)=d$.

Let $\sB:D \times D \to (0,\infty)$ be a \emph{symmetric} Borel function satisfying the following assumption:

\medskip

\noindent  \hypertarget{A}{{\bf (A)}} 
 There exist a blow-up weight function $\Phi$  with upper Matuszewska index $\overline \beta$ satisfying
\begin{align}\label{e:beta-range}
\overline \beta<1 \wedge \alpha
\end{align} 
and constants $A_0 \in (0,\infty]$ and  $C_0>1$  such that for all $x,y \in D$,
\begin{align}\label{e:J-Phi}
	C_0^{-1}\Phi\left( \frac{(|x-y|\wedge A_0)^2}{(\delta_D(x)\wedge A_0) (\delta_D(y) \wedge A_0)}\right)\le 	\sB(x,y) 
	\le C_0\Phi\left( \frac{(|x-y|\wedge A_0)^2}{(\delta_D(x)\wedge A_0) (\delta_D(y) \wedge A_0)}\right).
\end{align}
Here $\delta_D(x)=\mathrm{dist}(x, D^c)$. 
\bigskip

Set $J(x,y):=\sB(x,y)|x-y|^{-d-\alpha}$, $x,y\in D$, and consider the bilinear form
\begin{align*}
	\sE^0(u,v):= \frac12 \int_{D\times D} (u(x)-u(y))(v(x)-v(y)) J(x,y) dxdy,
\end{align*}

Let ${\rm Lip}_c(\overline D)$ be the family of all Lipschitz functions on $\overline D$ with compact supports. Note that the $d$-dimensional Lebesgue measure of $\partial D$ is zero. Hence we can identify $L^2(\overline D)$
with $L^2( D)$. Under \As, $\sE^0(u,u)<\infty$ for all $u\in {\rm Lip}_c(\overline D)$, see \cite{CKSV25}.

Let  $\overline \sF$ be the closure of $\mathrm{Lip}_c(\overline D)$ in 
$L^2(D)$
under  the norm $(\sE^0_1)^{1/2}$ where  $\sE^0_1:=\sE^0 + \lVert \cdot \rVert_{L^2(D)}^2$.  
Then  $(\sE^0, \overline \sF)$ is a regular Dirichlet form on 
$L^2(D)$.

It follows from \cite[Theorem 1.2]{CKSV25} that $(\sE^0, \overline \sF)$ admits a jointly continuous heat kernel 
$\overline{p}(t,x,y)$ on $(0,\infty)\times \overline{D}\times \overline{D}$ satisfying the following sharp two-sided estimates: For every $T \in (0, \infty)$, there exists a constant $C \ge 1$ such that, for all $(t,x,y) \in (0, T \vee R_0^\alpha) \times \overline{D}\times \overline{D}$,
	\begin{align}\label{e:main-HKE1}
	&C^{-1} \left[t^{-d/\alpha} \wedge  \left( \frac{ t }{|x-y|^{d+\alpha}}\Phi\left( \frac{(|x-y| \wedge A_0)^2}{((\delta_D(x) \vee t^{1/\alpha}) \wedge A_0)((\delta_D(y) \vee t^{1/\alpha}) \wedge A_0) }\right) \right)\right] \nn\\
	&\le \overline{p}(t,x,y)  
	\le	C \left[t^{-d/\alpha} \wedge  \left( \frac{ t }{|x-y|^{d+\alpha}}\Phi\left( \frac{(|x-y| \wedge A_0)^2}{((\delta_ D(x) \vee t^{1/\alpha}) \wedge A_0)((\delta_D(y) \vee t^{1/\alpha}) \wedge A_0) }\right) \right) 
	\right].
	\end{align}
As a consequence of joint continuity of the heat kernel, we can conclude that the Hunt process $\overline{X}$ 
associated with the form $(\sE^0, \overline \sF)$ can be refined to be strongly Feller and can start from every point in 
$\overline{D}$.

In this paper we are interested in studying 
potential theoretic properties 
of two subprocesses of $\overline{X}$ in $C^{1,\eta}$ 
open sets:
the process obtained by killing $\overline{X}$ upon approaching the boundary $\partial D$, and the process obtained by killing $\overline{X}$ with a critical killing potential.

To be more precise, let $\sF^0$ be the closure of  $\mathrm{Lip}_c(D)$ in $L^2(D)$ under 
$(\sE^0_1)^{1/2}$. Then $(\sE^0,\sF^0)$ is a regular Dirichlet form. Let $X^0=(X^0_t,t \ge 0; \P_x,x \in D\setminus \sN_0)$ be the Hunt process associated with $(\sE^0,\sF^0)$, where  
$\sN_0$ is an exceptional set for $X^0$. By the same argument as that in \cite[Subsection 4.2]{CKSV23} we see that $X^0$ can be regarded as a part of $\overline{X}$, hence can be refined to be a Hunt process which can start from every point in $D$. 
In the case $\alpha\in (1,2)$, we can follow the proofs of \cite[Lemmas 4.12 and 4.13]{CKSV23} (only lower boundedness of $\sB$ is used) to conclude that $\sF^0$ is strictly contained in $\sF$, hence $X^0$ is a proper subprocess of $\overline{X}$. 
We use $\zeta^0$ to denote the lifetime of $X^0$. 

The other 
process is obtained by killing $X^0$ via  a killing potential $\kappa$. We assume that $\kappa$ is a non-negative 
locally integrable  Borel function on $D$, 
and consider a symmetric form $(\sE^\kappa, \sF^\kappa)$ defined by
	\begin{align*}
		\sE^\kappa(u,v)&=\sE^0(u,v) + 
		\int_{D} u(x)v(x)\kappa(x)dx, \\
		\sF^\kappa&= 	\wt \sF^0 \cap L^2(D, 
			\kappa(x) dx),\nn
	\end{align*}
where $\wt \sF^0$ is the family of all $\sE^0_1$-quasi-continuous functions in $\sF^0$. Then $(\sE^\kappa,\sF^\kappa)$ is a regular Dirichlet form on $L^2(D)$ with Lip$_c(D)$ as a special standard core,  see \cite[Theorems 6.1.1 and 6.1.2]{FOT}.
Let  $X^\kappa=(X^\kappa_t,t \ge 0; \P_x,x \in D \setminus  \sN_\kappa)$ be the Hunt process associated with 
$(\sE^\kappa, \sF^\kappa)$ where $\sN_\kappa$ is an exceptional set for $X^\kappa$. We denote by $\zeta^\kappa$ the lifetime of $X^\kappa$, and define $X^\kappa_t=\partial$ for $t \ge \zeta^\kappa$, where $\partial$ is a cemetery point added to the state space $D$. Again, $X^{\kappa}$ can be regarded as a part of $X^0$, hence can be refined to be a Hunt process that can start from every point in $D$. 

From now on we will use the notation $X$ to denote either $X^0$ (when there is no killing, 
i.e., $\kappa\equiv 0$) 
or $X^{\kappa}$,  and
use $\zeta$ to denote the lifetime of $X$. 
We use the usual convention that a function $u$ on $D$ is automatically extended to $D\cup \{\partial\}$ by setting $u(\partial)=0$.
In order to analyze $X$ and study its potential theory we need additional assumptions
on the smoothness of the state space $D$, as well as on the jump kernel $J(x,y)$, or, more precisely, the function $\sB(x,y)$.

Let us start with some geometric definitions.

\begin{defn}\label{df-c11}
\rm 	For $\eta\in (0,1]$,	we say that $D$ is a $C^{1,\eta}$ open set with characteristics  $(\wh R, \Lambda)$, if 
	for each $Q \in \partial D$, there exist a $C^{1,\eta}$ function $\Psi=\Psi^Q:\R^{d-1} \to \R$ with 
	\begin{align}\label{e:local-map}
		\Psi(\wt 0)= |\nabla \Psi(\wt 0)|=0 \quad \text{and} \quad |\nabla \Psi(\wt y)-\nabla\Psi(\wt z)| \le \Lambda |\wt y-\wt z|^\eta \; \text{ for all } \,  \wt y, \wt z \in \R^{d-1},
	\end{align}
	and an orthonormal coordinate system CS$_{Q}$ with origin at $Q$ such that
	\begin{align}\label{e:local-coordinate-system}
		B_D(Q, \wh R) = \left\{ y= (\wt y, y_d) \in B(0, \wh R) \text{ in CS$_{Q}$} : y_d>\Psi(\wt y) \right\}.
	\end{align}
\end{defn}

\smallskip

\emph{Throughout this paper, we assume that $D$ is 
a $C^{1,\eta}$ open set with characteristics $(\wh R,\Lambda)$ such that 
${\rm dim}_{\rm A}(\partial D)=d-1$. Without loss of generality, we also assume $\wh R \le 1 \wedge (1/(4\Lambda))^{1/\eta}$. }

\bigskip

For $Q\in \partial D$, $\nu\in (0,1]$, and 
$R\in (0, \wh R/4]$,
let
\begin{align}\label{e:cone-like}
E^Q_{\nu}(R):=\left\{y=(\wt{y},y_d) \text{ in CS$_Q$}:\, |\wt{y}|<R/4, 4R^{-\nu}|\wt{y}|^{1+\nu}<y_d<R/2\right\}
\end{align}
be the paraboloid-like subset of $D$ near the boundary point $Q$. 
In case $d=1$ we interpret $E^Q_{\nu}(R)$ as the interval $(0,R/2)$.
From the assumption $\wh R^\eta \le1 \wedge (1/(4\Lambda))$ and
\eqref{e:Psi-maximum} below, 
it is easy to see that  
$E^Q_{\nu}(R)\subset E^Q_{\nu'}(R)\subset D$ for 
$0<\nu\le \nu'\le  \eta $.

\medskip
We now state additional 
assumptions on $\sB$: 
\medskip

\noindent\hypertarget{B1}{{\bf (B1)}} If $\alpha \ge 1$, then there exist constants $\theta_0>\alpha-1$ and $C>0$ such that
	\begin{align}\label{e:ass-B1}
		|\sB(x,x)-\sB(x,y)| \le C \bigg(\frac{|x-y|}{\delta_D(x) \wedge \delta_D(y) \wedge  \wh R} \bigg)^{\theta_0} 
			\quad \text{for all }  x,y \in D.
	\end{align}

This condition is needed to assure that $C^2$ functions are in the domain of the nonlocal operator associated with the form $\sE^0$. It appeared already in \cite{CKSV23, CKSV25b, KSV23, KSV25}.

Condition \hyperlink{B2}{{\bf (B2)}} is more delicate -- 
roughly 
it serves as a substitute of the flattening of the boundary method. This condition was introduced in \cite{CKSV23} where a more detailed explanation and  
justification are given  in the case when the jump kernel decays to zero at the boundary. The justification for the 
current
setting is given in the appendix where we show that this condition is satisfied in case of a trace process and the process related to the non-local Neumann problem, as well as in the case when in \eqref{e:J-Phi} we have exact equality.

For $a \in \R$, let $\bH_{a}:=\{x=(\wt{x}, x_d)\in \R^d: \, x_d>a\}$ and denote $\bH_0$ by $\bH$. 

\medskip
\noindent  \hypertarget{B2}{{\bf (B2)}} There  exist  constants $C>0$,  
$
\nu_0 
\in (0,\eta)$, 
$\theta_1>0$, $\theta_2>0$, and a non-negative Borel function $\F_0$ on 
$\bH_{-1}$   
with $\F_0(0)=1$
such that  for any $Q \in \partial D$ and    $x,y \in E^Q_{\nu_0}(\wh R/8)$ with $x=(\wt x,x_d)$ in CS$_Q$,
	\begin{align}\label{e:ass-B5-2}
			& \big|\sB(x,y)- \sB(x,x)\F_0((y-x)/x_d) \big| + \big|\sB(x,y)- \sB(y,y)\F_0((y-x)/x_d) \big|\nn \\	
			&\le C\left(\frac{|x-y|\vee \delta_D(x)\vee \delta_D(y)}{|x-y|\wedge \delta_D(x)\wedge \delta_D(y)}\right)^{\theta_1} \, 
			(|x-y|\vee\delta_D(x)\vee \delta_D(y))^{\theta_2}.
\end{align}
\medskip

Define a  function $\F$ on $\bH_{-1}$ by 
	\begin{align}\label{e:F-symmetrization}
		\F(y) =\frac{\F_0(y) + \F_0(-y/(1+y_d))}{2}, \quad \;\;y  = (\wt y,y_d) \in \bH_{-1}.
	\end{align} 
	By  construction, we have $\F(y)=\F(-y/(1+y_d))$ for all $y \in \bH_{-1}$. This identity is related to the symmetry property of $\sB$; see the discussion at the beginning of Section \ref{s:F}. 
With the function $\F$ in \eqref{e:F-symmetrization} and 
$q \in [(\alpha-1)_+, \alpha- \overline{\beta})$, we associate a  constant $C(\alpha,q, \F)$  defined by
\begin{align}\label{e:def-killing-constant-2}
			C(\alpha,q,\F) =\begin{cases}	 
	\displaystyle \int_{\R^{d-1}}\frac{1}{(|\wt u|^2+1)^{(d+\alpha)/2}}\int_0^1  
		\frac{(s^q -1)(1-s^{\alpha-1-q})}{(1-s)^{1+\alpha}} \F\big(((s-1)\wt u, s-1)\big) ds \, d \wt u, 
	 &  d\ge 2 \\[8pt]
	\displaystyle \int_0^1 \frac{(s^q-1)(1-s^{\alpha-q-1})}{(1-s)^{1+\alpha}} 
	\F(s-1)ds,  & d=1.
	\end{cases}
\end{align}
We will show in Lemma \ref{l:constant} that	$q \mapsto C(\alpha,q, \F)$ is a well-defined strictly increasing continuous function 
on $[(\alpha-1)_+, \alpha- \overline{\beta})$, $C(\alpha,(\alpha-1)_+,\F)=0$. 

The last condition concerns  the killing potential $\kappa$.

\medskip
\noindent  \hypertarget{K}{{\bf (K)}} 
There exist a non-negative constant $\kappa_0 <\lim_{q \uparrow \alpha- \overline{\beta}} C(\alpha,q,\F)$ and positive constants   $\theta_3, C>0$	such that for all $x \in D$,
	\begin{align}\label{e:def-K}
	\begin{cases}			
	|\kappa(x) - \kappa_0 \sB(x,x) \delta_D(x)^{-\alpha}| \le C\delta_D(x)^{-\alpha+\theta_3} &\text{ if }  \delta_D(x) < 1,\\[3pt]			\kappa(x) \le C &\text{ if } \delta_D(x) \ge 1.		
	\end{cases}	
	\end{align} 		
When $\alpha \le 1$, we further assume that $\kappa_0>0$.

\medskip
Note that condition \hyperlink{K}{{\bf (K)}} states that, in the case when $\kappa_0>0$, the killing potential is \emph{critical} in the sense that near the boundary the main contribution comes from the term $\kappa_0\sB(x,x)\delta_D(x)^{-\alpha}$. 
Assuming $\kappa_0>0$ when $\alpha \le 1$ ensures that, in this case, the killing potential $\kappa$ is non-trivial.

The following observation is fundamental as it connects the constant $\kappa_0$ with a parameter 
 $p\in [(\alpha-1)_+, \alpha-\ub)$ which governs the decay rate of harmonic functions. Since we assume 
\begin{align}\label{e:K-upper-bound}	
	\kappa_0< \lim_{q \uparrow \alpha- \overline{\beta}} C(\alpha,q,\F),
\end{align}
there exists a unique constant 
$p\in [(\alpha-1)_+, \alpha - \overline{\beta})$
such that 
\begin{align}\label{e:C(alpha,p,F)}
	\kappa_0=C(\alpha,p, \F).
\end{align}
We note that in the case $\Phi(r)=(1\vee r^{\overline{\beta}}) (\log(e \vee r))^a$ with $a\ge 0$ an arbitrary constant, \eqref{e:K-upper-bound} is automatically satisfied.
Indeed, 
under   \hyperlink{A}{{\bf (A)}}  and \hyperlink{B2}{{\bf (B2)}}, the function $\F$ satisfies the following estimates in terms of the blow-up weight function $\Phi$ (see Lemma \ref{l:F-basic}(ii) below):
	\begin{align*}
	C^{-1} \Phi\left(\frac{|y|^2}{y_d+1}\right)\le  \F(y)  \le C \Phi\left(\frac{|y|^2}{y_d+1}\right) \quad \text{for all $y = (\wt y, y_d) \in \bH_{-1}$}.
\end{align*} 
Using this, one can verify that $\lim_{q\uparrow \alpha-\overline\beta}C(\alpha, q, \F)=\infty$, so that \eqref{e:K-upper-bound} 
is automatically satisfied. 
Blow-up weight functions of this form arise, in particular, for the trace of a stable process and the process associated with the nonlocal Neumann problem.

The main result of this paper is the following boundary Harnack 
principle.

\begin{thm}\label{t:BHPnew} \!{\rm(Boundary Harnack principle)} 
Suppose  that $D\subset \R^d$
is a $C^{1,\eta}$ open set
with  $\eta\in (0,1]$, that ${\rm dim}_{\rm A}(\partial D)=d-1$, 
  and that conditions \hyperlink{A}{{\bf (A)}},  \hyperlink{B1}{{\bf (B1)}}, 
\hyperlink{B2}{{\bf (B2)}} and  \hyperlink{K}{{\bf (K)}}  hold.
Then there exists a  constant $C \ge 1$
such that, for  any $Q \in \partial D$, $0<r \le \wh{R}$, and any non-negative Borel function $f$ in $D$ which is harmonic in $B_D(Q,r)$ with respect to $X$ and vanishes continuously on $\partial D\cap B(Q,r)$, 	we have 
	\begin{equation}\label{e:TAMSe1.8new}
		\frac{f(x)}{\delta_D(x)^p}\le C	\frac{f(y)}{\delta_D(y)^p} \quad \text{for} \;\, x,y\in B_D(Q,r/2).
	\end{equation}
\end{thm}

We emphasize that the decay rate $p$ is related to $\kappa_0$ through equation \eqref{e:C(alpha,p,F)}.
Heat kernel estimates of the killed processes will be studied in the companion paper \cite{CKSV26b}.

\medskip

In the remaining part of the introduction, we describe the organization of the paper and explain the methods.

In Section \ref{s:est-geometric}, we establish auxiliary results on the geometry of $C^{1,\eta}$ open sets, introduce the 
regularized distance function and provide basic estimates. We end the section with several integral estimates involving the weight function $\Phi$.

Section \ref{s:F} gives results on the function $\F$ and the related constant $C(\alpha, q, \F)$ from \eqref{e:def-killing-constant-2}. 

In Section \ref{s:dynkin-X}, we provide analytic and probabilistic properties of the processes $\overline{X}$ and $X$. These include preliminary estimates involving the heat kernel of $\overline{X}$ and its part processes, as well as estimates of the $\lambda$-resolvent. In Proposition \ref{p:censored}, we show that when $\alpha \le 1$, $X^0$ and $\overline{X}$ have the same law when started in $D$, analogous to the well-known case of the censored process (cf.~\cite{BBC03}).

Section \ref{s:barriers} contains the most significant analytic results of the paper. Here we introduce the operator
$$
L^0_{\sB}u(x) 
:= \operatorname{p.v.} \int_D (u(y) - u(x)) J(x,y) \, dy,
$$
and analyze in Proposition \ref{p:barrier} its action on functions that behave near the boundary as powers of the regularized distance. This result is an analog of the fundamental \cite[Proposition 6.9]{CKSV23}. Its proof is one of the main novelties of our approach, in which the regularized distance function plays a crucial role. Using Proposition \ref{p:barrier}, we further construct barrier functions and establish their properties in Lemma \ref{l:Dynkin1}. Based on these barrier functions, Proposition \ref{p:explicit-1} provides preliminary boundary estimates involving certain exit probabilities.

In Section \ref{s:box}, we prove one of the main results of the paper, Theorem \ref{t:explicit-equiv}, which gives the exact decay rate of a certain survival probability. The key ingredient in the proof is Proposition \ref{p:explicit-2}, showing that this survival probability is dominated by the survival probability with vertical escape. The proof relies on a modification (based on \cite{CC24}) of the time-honored box method.

In Section \ref{s:green-pot}, based on the explicit decay rate from Theorem \ref{t:explicit-equiv}, we establish Green potential estimates in Proposition \ref{p:bound-for-integral-new}. 

In Section \ref{s:bhp} we state and prove the Carleson inequality and the boundary Harnack principle (Theorem \ref{t:BHPnew}). 
Although the arguments are similar to those in \cite{KSV25} once the results from Sections \ref{s:box} and \ref{s:green-pot} are established, we spell out the proofs for completeness.

Finally, in Section \ref{s:appl}, we provide two classes of examples satisfying conditions \hyperlink{A}{{\bf (A)}}, \hyperlink{B1}{{\bf (B1)}}, and \hyperlink{B2}{{\bf (B2)}}. These are processes related to nonlocal Neumann boundary conditions and traces of the $\alpha$-stable process. The proof that condition \hyperlink{A}{{\bf (A)}} holds was established in \cite{CKSV25}, whereas the rather long proofs for conditions \hyperlink{B1}{{\bf (B1)}} and \hyperlink{B2}{{\bf (B2)}} are deferred to the Appendix.

We use $m_d$ to denote the $d$-dimensional Lebesgue measure. 
 For a Borel set $E\subset \R^d$, $x\in \R^d$ and $r>0$, we use the shorthand notation $B_{E}(x, r):=B(x,r)\cap E$.
Throughout this paper, $A\asymp B$ means that 
 the quotient of $A$ and $B$
 is bounded between two positive constants on the specified region.
 Lower case letters  $c_i$, $i=1,2,  \dots$, are used to denote
positive constants in the proofs and the labeling of these constants starts anew in each proof.

We end the introduction with the list of notation used in the paper. 

\smallskip

\begin{center}
{\bf List of Notations}
\end{center}
$J(x,y):=|x-y|^{-d-\alpha}\sB(x,y)$,\\
$\sE^0(u,v):= \frac{1}{2}\iint_{D\times D}  (u(x)-u(y))(v(x)-v(y)) J(x, y)dydx$,\\
$\sE^{\kappa}(u,v):=\sE^0(u,v)+\int_{D}u(x)v(x)\kappa(x)\, dx$,\\
$\overline{\FF}$:   the closure of  $C_c^\infty(\overline{D})$ in $L^2(D,dx)$ under the norm  
$(\sE^0_1)^{1/2}$,\\
$\sF^0$:   the closure of  $C_c^\infty(D)$ in $L^2(D,dx)$ under the norm  
$(\sE^0_1)^{1/2}$,\\
 $\sF^\kappa:= 	\wt \sF^0 \cap L^2(D, \kappa(x)dx)$,\\
$\overline X$:  Hunt process associated with  $(\sE^0, \overline \FF)$; \ \ $\overline p(t,x,y)$:    the heat kernel of $\rX$,\\
$X$:   Hunt process associated with $(\sE^\kappa, \sF^\kappa)$ (including $\kappa\equiv 0$); \ \ $p(t,x,y)$:    the heat kernel of $X$,\\
$X^0$: Hunt process associated with  $(\sE^0, \sF^0)$,\\
$L^0_{\sB}u(x):= \textrm{p.v.} \int_{D}(u(y)-u(x))J(x,y)\, dy$; \ \  
$L^{\kappa}_{\sB}u(x):=L^0_{\sB}u(x)-\kappa(x)u(x)$,\\
$\Phi$: blow-up weight function, \eqref{e:Phi-scaling}\\
$\Psi$: $C^{1,\eta}$ function, Definition \ref{df-c11}\\
$\F_0$ and $\F$: \hyperlink{B2}{{\bf (B2)}}, \eqref{e:ass-B5-2}, and \eqref{e:F-symmetrization}\\
$\bH_{a}:=\{x=(\wt{x}, x_d)\in \R^d: \, x_d>a\}$;  $\bH=\bH_0$; $\bH^Q:=\left\{ y= (\wt y, y_d) \in \R^d \text{ in CS$_{Q}$} : y_d>0 \right\}$
\\[6pt]
\textbf{Constants:}\\
$\overline{\beta}$: The upper Matuszewska index of $\Phi$, 
 \eqref{e:Phi-scaling}\\
$C(\alpha,q,\F)$: \eqref{e:def-killing-constant-2}\\
$\theta_0$: constant in \hyperlink{B1}{{\bf (B1)}}, 
\eqref{e:ass-B1}\\
$\theta_1$ and $\theta_2$: constants in \hyperlink{B2}{{\bf (B2)}}, \eqref{e:ass-B5-2}\\
$\kappa_0$ and $\theta_3$: constants
in \hyperlink{K}{{\bf (K)}}, \eqref{e:def-K} \\
$\theta_4(p)$: constant in Proposition \ref{p:barrier} \\
$C_0$: constant in \eqref{e:J-Phi}\\
$C_1$: constant in 
 \eqref{e:regular-distance-1} \\
$C_2$: constant in 
 \eqref{e:regular-distance-2} \\
$\eps_1$: constant in Lemma \ref{l:Dynkin1}\\

%%%%%%%%%%%%%%%%%%%%%%%%%%%%%%%%%%%%%%%%%%%%%%%%%%%%%%%%%%%%%%%%%%%%%%%%%%%%%%%%%%%%%%%%%%%%%%%%%%%%%%%%%%%%%%%%%%%%%%%%%%%%%%%%%%%%
%%%%%%%%%%%%%%%%%                                  Estimates on $\Phi$                    %%%%%%%%%%%%%%%%%%%%%%%%%%%%%%%%%%%%%%%%%%%
%%%%%%%%%%%%%%%%%%%%%%%%%%%%%%%%%%%%%%%%%%%%%%%%%%%%%%%%%%%%%%%%%%%%%%%%%%%%%%%%%%%%%%%%%%%%%%%%%%%%%%%%%%%%%%%%%%%%%%%%%%%%%%%%%%%%%

\section{Estimates related to  $C^{1,\eta}$ open sets}\label{s:est-geometric}
 Throughout this paper, we assume that \textit{$D$ is a $C^{1,\eta}$ open set  with characteristics $(\wh R,\Lambda)$ such that  
  $\wh R \le 1 \wedge (1/(4\Lambda))^{1/\eta}$, ${\rm dim}_{\rm A}(\partial D)=d-1$, 
 $\Phi$ is a blow-up  weight function  with upper Matuszewska index $\overline \beta$, and that \hyperlink{A}{{\bf (A)}} holds.} 
 Observe that,  since $\Phi$ has upper Matuszewska index $\overline \beta$, for every  $\beta>\overline \beta$, there exists a constant 
 $C=C(\beta)\ge 1$ such that
 \begin{align}\label{e:Phi-scaling}
 	\frac{\Phi(r)}{\Phi(s)} \le C\bigg(\frac{r}{s}\bigg)^{ \beta} \quad \text{for all $0<s\le r$}.
 \end{align}
Since $\Phi$ is increasing, it follows immediately that
\begin{align}\label{e:Phi-scaling-2}
\frac{\Phi(r)}{\Phi(s)} \le
C\bigg(1+\bigg(\frac{r}{s}\bigg)^{\beta}\bigg) \quad \text{for all $r,s>0$}.
\end{align}

In this section, we establish several estimates related to the   $C^{1,\eta}$ open set $D$. In Subsection \ref{ss:C1eta}, we recall two lemmas from \cite{CKSV25} and establish a relation between $\delta_D(x)$ and $x_d$ in the local coordinate systems for $x\in E_\nu^Q(\wh R/8)$.  
 In Subsection \ref{ss:regular-distance},  we introduce a regularized distance function for $D$ and derive   some estimates about it. In Subsection \ref{ss:integral-estimates}, we derive several integral estimates over subsets of $D$ involving the weight function $\Phi$.  
The results of this section will play important roles in the subsequent analysis.

\subsection{Some results for $C^{1,\eta}$ open sets}\label{ss:C1eta}
We begin by recalling the following two lemmas, which were established in \cite[Lemmas 2.6 and 2.8]{CKSV25}.

\begin{lemma}\label{l:Aikawa-integrals}
(i) For any $q \in [0,1)$,	there exists $C=C(q)>0$ such that for all  $x\in  {\overline D}$, $A\in (0,\infty]$ and  $r>0$,
\begin{align*}
	\int_{B_D(x,r)} \frac{dy}{(\delta_D(y) \wedge A)^{q}}   \le \frac{Cr^{d}}{((\delta_D(x) \vee r) \wedge A)^{q}}. 
\end{align*}
\noindent (ii) There exists $C>0$ such that  for all $x\in  {\overline D}$, $A\in (0,\infty]$, $B\in (0,\infty)$ and  $r>0$,  
	\begin{align*}
		\int_{B_D(x,r)}\Phi\bigg( \frac{B}{\delta_D(y) \wedge A}\bigg)  dy \le Cr^{d} \Phi\bigg( \frac{B}{(\delta_D(x) \vee r) \wedge A}\bigg). 
	\end{align*}
\end{lemma}

\begin{lemma}\label{l:p-estimate0}
 (i) For any  $q\in [0,1)$ and  $\theta >q$,  there exists $C=C(q,\theta)>0$  such that  for all 
 $x \in  \overline D$, $A\in (0,\infty]$, and $r>0$,
\begin{align*}
	\int_{D\setminus B(x,r)}	 \frac{ 1 }{|x-y|^{d+\theta}}\left( \frac{(|x-y| \wedge A)^2}{\delta_D(y) \wedge A }\right)^q dy \le \frac{C(r\wedge A)^q}{r^{\theta}} .
\end{align*}

\noindent (ii) For every 
	$\theta >\ub$,  there exists $C=C(\theta)>0$  such that 
	 for all $x \in  \overline D$, $A\in (0,\infty]$, $B\in (0,\infty)$  
and $r>0$,
	\begin{align*}
		\int_{D\setminus B(x,r)}	 \frac{ 1 }{|x-y|^{d+\theta}}\Phi\left( \frac{(|x-y| \wedge A)^2}{B(\delta_D(y) \wedge A) }\right)dy \le \frac{C}{r^{\theta}} \Phi\left( \frac{r \wedge A}{B}\right).
	\end{align*}
\end{lemma}

 We note that Lemmas \ref{l:Aikawa-integrals} and \ref{l:p-estimate0} hold under weaker geometric assumptions on $D$, for instance, when $D$ is a Lipschitz open set. 

Observe that  
\begin{align}\label{e:without-A0}
		\frac{(r \vee B)\wedge A}{(s \vee B) \wedge A} \le \frac{r}{s} \quad \text{for all $A\in (0,\infty]$,  $B \in [0,\infty)$ and $0<s\le r$}.\end{align}
As a consequence of Lemma \ref{l:p-estimate0} and \eqref{e:J-Phi}, we obtain the following tail estimates for the jump kernel over exterior balls: there exists $C>0$ such that for all $x\in \overline D$ and $r>0$,
\begin{align}\label{e:jump-tail-estimate}
		&	\int_{D\setminus B(x,r)} J(x,y)dy \nn\\
		&\le 
		C_0	\int_{D\setminus B(x,r)} \frac{\sB(x,y)}{|x-y|^{d+\alpha}}dy \le \frac{C}{r^{\alpha}}  \Phi\left( \frac{r \wedge A_0}{\delta_D(x) \wedge A_0}\right) \le \frac{C}{r^{\alpha}}  \Phi\left( \frac{r}{\delta_D(x)}\right).
\end{align}
The last inequality follows from \eqref{e:without-A0} and the monotonicity of $\Phi$.

 	Let $Q\in \partial D$ and let $\Psi=\Psi^Q$ be the function in Definition \ref{df-c11}.   
 By \eqref{e:local-map}, for all $\wt z \in \R^{d-1}$,	\begin{align}\label{e:Psi-maximum}		|\Psi(\wt z)| = |\Psi(\wt z) - \Psi(\wt 0)| \le \sup_{|\xi|\le |\wt z|} |\nabla \Psi(\wt \xi)||\wt z| =  \sup_{|\xi|\le |\wt z|} |\nabla \Psi(\wt \xi)- \nabla \Psi(\wt 0)||\wt z| \le \Lambda |\wt z|^{1+\eta}.	\end{align}

 Recall that $E^Q_\nu(R)$ is defined in \eqref{e:cone-like}.
 
 \begin{lemma}\label{l:distance-in-Enu}
	(i) For any $Q\in \partial D$ and  $x=(\wt x, x_d) \in E_\eta^{Q}(\wh R/8)$ in the coordinate system \text{\rm CS}$_{Q}$,  
 	\begin{align*}
 		\frac{4}{5}\delta_D(x) \le  x_d \le 2(1+2\Lambda)\delta_D(x).
 	\end{align*} 
 	\noindent (ii) For any $\nu \in (0,\eta)$, there exists $C=C(\nu)>1$ such that for any $Q\in \partial D$ and $x=(\wt x, x_d) \in E_\nu^{Q}(\wh R/8)$ in the coordinate system \text{\rm CS}$_{Q}$,
 	\begin{align*}
 		|x_d -\delta_D(x)| \le C \delta_D(x)^{(1+\eta)/(1+\nu)}.
 	\end{align*}
 \end{lemma}
 \begin{proof}
 	Let $Q\in \partial D$.	We use the coordinate system CS$_{Q}$ and let $\Psi=\Psi^Q$ be the function in Definition \ref{df-c11}.  
 	Let $x=(\wt x, x_d) \in  E^Q_\eta(\wh R/8)$ and let $Q_x=(\wt Q_x, \Psi(\wt Q_x)) \in \partial D$
 	satisfy $\delta_D(x)=|x-Q_x|$.  
 	Since $ (\wt x, \Psi(\wt x)) \in \partial D$,  we get from \eqref{e:Psi-maximum}  that
 	\begin{align}\label{e:distance-in-Enu-2}
 		\delta_D(x) \le |x - (\wt x, \Psi(\wt x))| \le x_d + |\Psi(\wt x)| \le x_d + \Lambda |\wt x|^{1+\eta}.
 	\end{align}
 	Using  \eqref{e:Psi-maximum}, we also get that
 	\begin{align*}
 		x_d \le | (\wt x,x_d) - (\wt Q_x, 0) |  \le |x-Q_x|  + |\Psi (\wt Q_x)|  \le \delta_D(x) + \Lambda |\wt Q_x|^{1+\eta}.
 	\end{align*}
 	Note that $|\wt Q_x|^{1+\eta} \le  (|\wt x - \wt Q_x| + |\wt x|)^{1+\eta} \le  (\delta_D(x) + |\wt x|)^{1+\eta}\le 2^{1+\eta} (\delta_D(x)^{1+\eta} + |\wt x|^{1+\eta}).$ Thus, 
 	\begin{align}\label{e:distance-in-Enu-3}
 		x_d  \le \delta_D(x) + 2^{1+\eta} \Lambda (\delta_D(x)^{1+\eta} + |\wt x|^{1+\eta}).
 	\end{align} 
 	
 	\noindent (i) 		Since $\Lambda \wh R^\eta \le 1$, we see from the definition of $E^Q_\eta(\wh R/8)$ that
 	$\Lambda	|\wt x|^{1+\eta} \le 2^{-2}\Lambda ((\wh R/8)^{\eta}x_d) \le 2^{-2-3\eta}x_d.$ 
 	Hence, since $\delta_D(x) \le |x|<1/8$, we get from \eqref{e:distance-in-Enu-2} and \eqref{e:distance-in-Enu-3} that $\delta_D(x) \le (5/4)x_d$ and that
 	\begin{align*}
 		x_d\le (1+ 2^{1-2\eta} \Lambda )	\delta_D(x)  +  2^{-1-2\eta} x_d < (1+2\Lambda)\delta_D(x) + 2^{-1}x_d.
 	\end{align*}
 	These yield the desired result.
 	
 	\noindent (ii) Let $\nu \in (0,\eta)$ and suppose that $x\in E_\nu^Q(\wh R/8)$. Since $x\in E_\nu^Q(\wh R/8) \subset E_\eta^Q(\wh R/8)$, using (i), we obtain
 	$\Lambda	|\wt x|^{1+\eta} \le \Lambda ((\wh R/8)^\nu x_d/4)^{(1+\eta)/(1+\nu)} \le c_1 \delta_D(x)^{(1+\eta)/(1+\nu)}.$ 
 	Hence, by \eqref{e:distance-in-Enu-2} and \eqref{e:distance-in-Enu-3}, 
 	\begin{align*}
 		|x_d-\delta_D(x)| \le 2^{1+\eta}\Lambda \delta_D(x)^{1+\eta} + c_2  \delta_D(x)^{(1+\eta)/(1+\nu)}.
 	\end{align*} 
 	Since $\delta_D(x)<1/8$ and $1+\eta>(1+\eta)/(1+\nu)$, the desired result follows.
 \end{proof}

\subsection{Regularized distance function and its basic estimates}\label{ss:regular-distance}
Recall that $\delta_D(x)$ denotes the distance from $x$ to $D^c$.   For each $x\in D$,  we fix a $Q_x\in \partial D$  such that $\delta_D(x) = |x-Q_x|$.  For each $Q\in \partial D$, let $n_Q$ denote the inward unit normal to $\partial D$ at $Q$.

For $r>0$, let
\begin{align*}
	D(r):=\{x\in D: \delta_D(x)<r\}.
\end{align*} The signed distance to $\partial D$ is defined by
\begin{align}\label{e:signed-distance}
	\wt\delta_D(x) := \begin{cases}
		\delta_D(x) &\mbox{ if $x\in \overline D$},\\
		-\delta_{\overline D^c}(x) &\mbox{ if $x\in \overline D^c$}.
	\end{cases}
\end{align}

Recall that we have assumed $\wh R \le 1 \wedge (1/(4\Lambda))^{1/\eta}$, so that  $\Lambda \wh R^\eta \le 1/4$. 
By the construction at the beginning of \cite[Section 5]{CS25}, which goes back to \cite{Li85}, there exists a local \textit{regularized distance function}  $\rho \in C^1(\R^d) \cap C^2(\R^d \setminus \partial D)$  for $D$  satisfying the following properties: there are  constants $C_1>1$, $C_2>0$ and $C_3>0$ 
depending only on  
 $d,\wh R$ and $\Lambda$   such that  
\begin{gather}
		\text{$C_1^{-1} \le \frac{\rho(x)}{(\wt \delta_D(x) \vee (-\wh R)) \wedge \wh R}\le C_1$ } \quad \text{for all $x\in \R^d$,}\label{e:regular-distance-1}\\
	\text{$|\nabla\rho(x) - \nabla \rho(y)|\le C_2 |x-y|^\eta$} \quad \text{for all $x,y\in \R^d$,} \label{e:regular-distance-2}\\
	\text{$|D^2 \rho(x)| \le C_3(1+  |\wt \delta_D(x)|^{\eta-1})$} \quad \text{for all $x\in \R^d \setminus \partial D$.}\label{e:regular-distance-3}
\end{gather} 
Moreover,  it holds that 
\begin{align}\label{e:regular-disntance-inward-normal}
	\frac{\nabla \rho(Q)}{|\nabla \rho(Q)|}  = n_Q \quad \text{for all $Q\in \partial D$},
\end{align}
and there exists
  $R_0=R_0(d,\wh R)\in (0, \wh R/8]$ 
such that
\begin{align}\label{e:regular-distance-4}
	\frac14 \le |\nabla \rho(x)| \le 2  \quad \text{for all $x\in \overline{D(R_0)}$}.
\end{align} 

In the remainder of this paper,  $R_0$ always stands for the constant in \eqref{e:regular-distance-4}. 
Note that $R_0\le 1/8$. 
As an immediate consequence of  \eqref{e:regular-distance-1}, we have 
\begin{align}\label{e:regular-distance_C1}
C_1^{-1}\delta_D(x)\le \rho(x)\le C_1\delta_D(x) \quad \text{for all $x\in \overline{D(R_0)}$}.
\end{align} 

 For $Q\in \partial D$, let
\begin{align*}
\bH^Q:=\left\{ y= (\wt y, y_d) \in \R^d \text{ in CS$_{Q}$} : y_d>0 \right\}.
	\end{align*}

	Let $\e_d:=(\wt 0,1) \in \R^d$.

\begin{lemma}\label{l:geometry}
 The following  hold for every $x\in D(R_0)$.

\noindent (i) It holds that
\begin{equation}\label{e:geometry-result-0}
	\left|\nabla\rho(x)-\frac{\rho(x)}{\delta_D(x)}n_{Q_x}\right|
	\le \frac{(2+\eta)C_2}{1+\eta}\delta_D(x)^\eta .
\end{equation}

\noindent (ii) There exists $C>0$ independent of $x$ such that for all $y =(\wt y, y_d)\in B_D(Q_x,R_0) \cap \bH^{Q_x}$ in {\rm CS}$_{Q_x}$, 
\begin{align}
		\rho(y) + \frac{\rho(x)}{\delta_D(x)}y_d & \le C( |x-y| \vee \delta_D(x)), \label{e:geometry-result-1}\\
	\left|\rho(y)-\frac{\rho(x)}{\delta_D(x)}y_d\right|
&	\le C \left(  |x-y|\vee \delta_D(x)\right)^\eta |x-y|.\label{e:geometry-result-2}
\end{align}
\end{lemma} 
\begin{proof}
	Let $x\in D(R_0)$ and $r:=\delta_D(x)$. We work in the coordinate system CS$_{Q_x}$, so that $Q_x=0$, $\bH^{Q_x}=\bH$ and $n_{Q_x}=\e_d$.

	\noindent (i)  
	Since $x=r\e_d$, $\rho(0)=0$ and 
	$\nabla\rho(0)=|\nabla\rho(0)|\e_d$, we have
	\begin{align*}
		\frac{\rho(x)}{r} - |\nabla \rho(0)| =\frac1r \left( \rho(r\e_d)  -  r \nabla \rho(0)\cdot \e_d \right)  =\frac1r\int_0^r
		\bigl(\nabla\rho(s\e_d)-\nabla\rho(0)\bigr)\cdot \e_d\,ds.
	\end{align*}
	Hence, by \eqref{e:regular-distance-2}, we get
	\begin{align}\label{e:geometry-1}
		\left|\frac{\rho(x)}{r}\e_d-\nabla\rho(0)\right|=	\left|\frac{\rho(x)}{r}-|\nabla\rho(0)|\right|\le \frac1r\int_0^r
			\bigl|\nabla\rho(s\e_d)-\nabla\rho(0)\bigr|ds \le \frac{C_2}{r}\int_0^r s^\eta\,ds
		=\frac{C_2}{1+\eta}r^\eta.
	\end{align}
	Using  \eqref{e:regular-distance-2} and \eqref{e:geometry-1}, we arrive at
	\begin{align*}
			\left|\nabla\rho(x)-\frac{\rho(x)}{r}\e_d\right|
		\le |\nabla\rho(r\e_d)-\nabla\rho(0)|+\left|\nabla\rho(0)-\frac{\rho(x)}{r}\e_d\right|
		\le C_2r^\eta+\frac{C_2}{1+\eta}r^\eta.
	\end{align*}

	\noindent (ii) Note that for any $y=(\wt y, y_d) \in \bH$, 
	\begin{equation}\label{eq:affine}
		\frac{\rho(x)}{r} y_d=	\frac{\rho(x)}{r} \left( r + (y-x) \cdot \e_d\right)  = \rho(x)+	\frac{\rho(x)}{r} (y-x) \cdot \e_d.
	\end{equation}
	Thus, by  \eqref{e:regular-distance-1}, we get that for all $y = (\wt y, y_d) \in B_D(0, R_0) \cap \bH$, 
	\begin{align}\label{e:geometry-result-1-1}
			\frac{\rho(x)}{r} y_d\le C_1(r+|x-y|) .
	\end{align}
	On the other hand, by  \eqref{e:regular-distance-1},  we have 
	 $\rho(y)  \le C_1 \delta_D(y) \le  C_1 (r+|x-y|)$ for all $y \in B_D(0, R_0) \cap \bH$.
	  Combining this with \eqref{e:geometry-result-1-1} yields  \eqref{e:geometry-result-1}.
	
	Next, we prove \eqref{e:geometry-result-2}.
By \eqref{eq:affine},  for all $y = (\wt y, y_d) \in B_D(0, R_0) \cap \bH$, 
\begin{align*}
	&	\left|\rho(y)-\frac{\rho(x)}{r}y_d\right| = 
	\left|\rho(y) - \rho(x) - 	\frac{\rho(x)}{r} (y-x) \cdot \e_d\right|\\
	& \le \left| \rho(y) - \rho(x)  - \nabla \rho(x) \cdot (y-x)\right| + \left| \left( \nabla \rho(x) -\frac{\rho(x)}{r} \e_d \right) \cdot (y-x)\right|=:I_1+I_2.
\end{align*}
By \eqref{e:geometry-result-0}, we have $I_2 \le c_1 r^\eta |x-y|$. For $I_1$, using 	 \eqref{e:regular-distance-2}, we obtain
\begin{align*}
	I_1 & = \left| \int_0^1
	\bigl(\nabla\rho(x+s(y-x))-\nabla\rho(x)\bigr)\cdot(y-x)\,ds\right| \\
	&\le  |x-y| \int_0^1
	\bigl|\nabla\rho(x+s(y-x))-\nabla\rho(x)\bigr|\, ds \le C_2|x-y|^{1+\eta} \int_0^1 s^\eta ds = \frac{C_2|x-y|^{1+\eta}}{1+\eta}.
\end{align*}
This completes the proof of \eqref{e:geometry-result-2}.
\end{proof}

\subsection{Integral estimates with the weight function}\label{ss:integral-estimates}
 
In this subsection, we establish several integral estimates involving the weight function $\Phi$. These estimates will be used later in the construction of barrier functions.

We begin with the following simple consequence of Lemma \ref{l:Aikawa-integrals}.

\begin{lemma}\label{l:tail-F}
Let $\nu \in (0,\eta)$ and 
	$\theta >\ub$.  There exists $C>0$  such that  for all $R\in (0,R_0]$, $x\in D(R/2)$ and $r\in (0,R)$,
	\begin{align*}
		\int_{E_\nu^{Q_x}(R)\setminus B(x,r)} \frac{1}{|x-y|^{d+\theta}}\Phi\left(\frac{|x-y|^2}{\delta_D(x) y_d}\right)dy   \le \frac{C}{r^\theta} \Phi \left( \frac{r}{\delta_D(x)}\right),
	\end{align*}
where, on the left-hand side,  we work in the coordinate system \text{\rm CS}$_{Q_x}$.
\end{lemma}
\begin{proof}
	Let $R\in (0,R_0]$,  $x\in D(R/2)$ and $r\in (0,R)$. We work in the coordinate system CS$_{Q_x}$. By \eqref{e:Phi-scaling} and Lemma \ref{l:distance-in-Enu}(i), 
	\begin{align*}
		&		\int_{E_\nu^{Q_x}(R)\setminus B(x,r)} \frac{\Phi(|x-y|^2/(\delta_D(x) y_d))}{|x-y|^{d+\theta}}dy
		\le c_1	\int_{E_\nu^{Q_x}(R)\setminus B(x,r)} \frac{\Phi(|x-y|^2/(\delta_D(x)\delta_D(y)))}{|x-y|^{d+\theta}}dy\\
		& \le c_1	\int_{D\setminus B(x,r)} \frac{\Phi(|x-y|^2/(\delta_D(x)\delta_D(y)))}{|x-y|^{d+\theta}}dy.
	\end{align*}
	The result now follows from Lemma \ref{l:p-estimate0}(ii) with $A=\infty$ and $B=\delta_D(x)$. 
\end{proof}

 The following lemma  will be used later in the construction of barrier functions to show that the contribution of the nonlocal part  
 $L^{0}_{\sB}$ of the operator $L^{\kappa}_{\sB}$ 
 from the complement of $E_\nu^Q(R)$ is small. This allows us to reduce, up to a small error, the computation to integration over $E_\nu^Q(R)$, where the additional assumption \hyperlink{B2}{{\bf (B2)}} applies. 
Recall that for $x\in D(R_0)$, 
$Q_x\in \partial D$ denotes the point 
such that $\delta_D(x)=|x-Q_x|$. 

\begin{lemma}\label{l:Enu-difference}
 Let $\nu \in (0,\eta)$ and $\theta>\ub$. 
For any $0<b<(\theta - \ub) \wedge ((1-\ub)\nu)$, there exists $C=C(\nu, \theta, b)>0$ such that the following estimates hold for all $R\in (0,R_0]$ and $x\in D(R/8)$:
\begin{gather}
 \int_{D\setminus E^{Q_x}_\nu(R)}
 \frac{\Phi(|x-y|^2/(\delta_D(x)\delta_D(y)))}{|x-y|^{d+\theta}}dy
 \le \frac{C}{\delta_D(x)^{\theta-b} R^b},\label{e:Enu-difference-result-1}\\
\int_{\{ y = (\wt y, y_d) \text{ in {\rm CS}$_{Q_x}$}:\, y \in \bH^{Q_x}\setminus E^{Q_x}_\nu(R)\}}
 \frac{\Phi(|x-y|^2/(\delta_D(x) y_d))}{|x-y|^{d+\theta}}dy
 \le \frac{C}{\delta_D(x)^{\theta-b} R^b}.\label{e:Enu-difference-result-2}
 \end{gather}
\end{lemma}
\begin{proof}
Let $b\in (0,(\theta - \ub) \wedge ((1-\ub)\nu))$ and let $\beta \in (\ub,  \theta\wedge 1)$ be such that $b<	(\theta - \beta) \wedge ((1-\beta)\nu)$. It suffices to prove \eqref{e:Enu-difference-result-1}, since \eqref{e:Enu-difference-result-2} follows by applying  \eqref{e:Enu-difference-result-1} to the case $D=\bH$.

Let $R\in (0,R_0]$, $x\in D(R/8)$ and $r:=\delta_D(x)$.  We work in the coordinate system CS$_{Q_x}$ and write $E_\nu:=E_\nu^{Q_x}(R)$  for simplicity. We present the proof only for $d\ge 2$, as the case $d=1$ follows by a similar but simpler argument.  
  By definition of $E_\nu$,  we have
\begin{align*}
	&\int_{D \setminus E_\nu}    \frac{1}{|x-y|^{d+\theta}} \Phi\bigg(\frac{|x-y|^2}{r\delta_D(y)}\bigg) dy\\
	&\le \int_{D \setminus B(0,R/4)}    \frac{1}{|x-y|^{d+\theta}} \Phi\bigg(\frac{|x-y|^2}{r\delta_D(y)}\bigg) dy \\
	&\quad  +   \int_{\{y=(\wt y,y_d) \in B_D(0,R/4): 
	y_d \le 4R^{-\nu}|\wt y|^{1+\nu}\}}  \frac{1}{|x-y|^{d+\theta}} \Phi\bigg(\frac{|x-y|^2}{r\delta_D(y)}\bigg) dy \\
	&=:I_1+I_2.
\end{align*}	
Since $|x|= r<R/8$ and $\theta>\ub$, by
Lemma \ref{l:p-estimate0}(ii) and \eqref{e:Phi-scaling}, we obtain that 
\begin{align}\label{e:Enu-difference-I1}
	I_1 &\le \int_{D\setminus B(x,R/8)}  \frac{1}{|x-y|^{d+\theta}} \Phi\bigg(\frac{|x-y|^2}{r\delta_D(y)}\bigg) dy \le  \frac{c_1}{R^\theta}  \Phi \bigg( \frac{R}{r}\bigg) \le \frac{c_2}{r^{\beta}R^{\theta-\beta}} \le \frac{c_2}{r^{\theta-b}R^{b}},
\end{align} 
where we used the facts that $b<\theta-\beta$ and $r<R/8$ in the last inequality.

Let $\Psi=\Psi^{Q_x}$ be the function in Definition \ref{df-c11}. By \eqref{e:local-coordinate-system}, we have 
\begin{align*}
	I_2=\sum_{n=0}^\infty \int_{U_n}   \frac{1}{|x-y|^{d+\theta}} \Phi\bigg(\frac{|x-y|^2}{r\delta_D(y)}\bigg) dy,
\end{align*}
where
\begin{align*}
	U_n:=\left\{(\wt  y ,y_d) \in B_D(0,R/4) : 2^{-n-3}R\le |\wt y|<2^{-n-2}R \text{ and }  
	\Psi(\wt y)
	< y_d\le 4R^{-\nu}|\wt y|^{1+\nu}  \right\} .
\end{align*}	
By \eqref{e:local-map}, since $\nu<\eta$, we get $|\Psi(\wt y)| \le \Lambda |\wt y|^{1+\eta} \le \Lambda R^{-\nu} |\wt y|^{1+\nu}$ for all $\wt y \in \R^{d-1}$ with $|\wt y|\le 1$. Thus, 
\begin{align*}
	U_n\subset \left\{(\wt  y ,y_d) \in B_D(0,R/4): 2^{-n-3}R\le |\wt y|<2^{-n-2}R \text{ and } |y_d|\le (\Lambda \vee 4)2^{-(1+\nu)(n+2)}R \right\}.
\end{align*}
Consequently,
$U_n$ can be covered by open balls 
$\{B(z_{n,i}, 2^{-(1+\nu)(n+2)}R) : 1\le i\le \ell_n\}$ with centers $z_{n,i}\in \overline{U_n}$, where $\ell_n$ satisfies $\ell_n\le c_3 (2^{\nu (n+2)})^{d-1}$ and $c_3>0$ is a constant depending only on $d$.

Let $N:=\min\{n: 2^{-n}R <r/2 \}$.
For all $0\le n\le N$ and $y\in U_n$, we have $2^{-n-3}R\le |\wt y|\le |x-y|\le |\wt y| + |r-y_d| \le 2^{-n-2}R + r + 4R^{-\nu}(2^{-n-2}R)^{1+\nu} \le c_42^{-n}R$.	 Hence,  using \eqref{e:Phi-scaling} and  Lemma \ref{l:Aikawa-integrals}(ii), we obtain
\begin{align*}
	&	\sum_{n=0}^N \int_{U_n}   \frac{1}{|x-y|^{d+\theta}} \Phi\bigg(\frac{|x-y|^2}{r\delta_D(y)}\bigg) dy \le  \frac{c_5}{R^{d+\theta}} \sum_{n=0}^N  2^{n(d+\theta)}\int_{U_n} \Phi\bigg(\frac{(2^{-n}R)^2}{r\delta_D(y)}\bigg) dy\nn\\
	&\le  \frac{c_5}{R^{d+\theta}} \sum_{n=0}^N  2^{n(d+\theta)} \sum_{i=1}^{
	\ell_n
	}\int_{B_D(z_{n,i}, 2^{-(1+\nu)(n+2)}R)} \Phi\bigg(\frac{(2^{-n}R)^2}{r\delta_D(y)}\bigg) dy\nn\\
	&\le \frac{c_6}{R^\theta}\sum_{n=0}^N  2^{n(\theta-d\nu)} 
	\ell_n
	\Phi\bigg(\frac{2^{-(1-\nu)n}R}{r }\bigg)  \le \frac{c_7}{R^{\theta-\beta}r^\beta }\sum_{n=0}^N  2^{n(\theta-\nu - \beta(1-\nu))} .
\end{align*}
Since $b<(1-\beta)\nu$ and $2^N\le 4R/r$, it follows  that
\begin{align}\label{e:diff-D-Enu-11}
	&	\sum_{n=0}^N \int_{U_n}   \frac{1}{|x-y|^{d+\theta}} \Phi\bigg(\frac{|x-y|^2}{r\delta_D(y)}\bigg) dy \nn\\
	&\le  \frac{c_7 2^{N( \theta-\beta - b)} }{R^{\theta-\beta}r^\beta }\sum_{n=0}^N  2^{-n( (1-\beta)\nu  -b) }  \le \frac{c_8 2^{N(\theta- \beta- b )} }{R^{\theta-\beta}r^\beta }  \le \frac{c_9}{r^{\theta-b}R^b}.
\end{align}
On the other hand, for all $n>N$ and $y\in U_n$, we get $|\wt y|<r/8$, $y_d < 4R^{-\nu}(r/8)^{1+\nu} <r/2$ and $y_d > -\Lambda R^{-\nu}|\wt y|^{1+\nu} >- (\Lambda/8) r$. In particular, $|x-y| \ge |r-y_d|>r/2$ and $|x-y|\le |\wt y| + |r-y_d| \le c_{10} r$. Hence, using \eqref{e:Phi-scaling}, $\ell_n\le c_3 (2^{\nu (n+2)})^{d-1}$  and  Lemma \ref{l:Aikawa-integrals}(ii), we obtain
\begin{align*}
	&	\sum_{n=N+1}^\infty \int_{U_n}   \frac{1}{|x-y|^{d+\theta}} \Phi\bigg(\frac{|x-y|^2}{r\delta_D(y)}\bigg) dy \le  \frac{c_{11}}{r^{d+\theta}} \sum_{n=N+1}^\infty \int_{U_n} \Phi\bigg(\frac{r}{\delta_D(y)}\bigg) dy\nn\\
	&\le  \frac{c_{11}}{r^{d+\theta}}  \sum_{n=N+1}^\infty \sum_{i=1}^{
	\ell_n
	}\int_{B_D(z_{n,i}, 2^{-(1+\nu)(n+2)}R)} \Phi\bigg(\frac{r}{\delta_D(y)}\bigg) dy\nn\\
	&\le   \frac{c_{12}R^d}{r^{d+\theta}}   \sum_{n=N+1}^\infty 2^{-n(d+d\nu)} 
	\ell_n
	\Phi\bigg(\frac{r}{2^{-(1+\nu)(n+2)}R}\bigg)\le   \frac{c_{13}R^{d-\beta}}{r^{d+\theta-\beta}}   \sum_{n=N+1}^\infty 2^{-n(d+\nu  - (1+\nu)\beta )}  .
\end{align*}
Since $d+\nu  - (1+\nu)\beta \ge (1+\nu)(1-\beta)>0$ and $2^N \ge 2R/r$, it follows that
\begin{align}\label{e:diff-D-Enu-2}
	&\sum_{n=N+1}^\infty \int_{U_n}   \frac{1}{|x-y|^{d+\theta}} \Phi\bigg(\frac{|x-y|^2}{r\delta_D(y)}\bigg) dy \le \frac{c_{14}2^{-N(d+\nu  - (1+\nu)\beta )} R^{d-\beta}}{r^{d+\theta-\beta}} \nn\\
	&\le \frac{c_{15}}{r^{\theta-(1-\beta)\nu}R^{(1-\beta)\nu}} \le \frac{c_{15}}{r^{\theta-b}R^{b}},
\end{align}
where we used the facts that $b<(1-\beta)\nu$ and $r<R/8$ in the last inequality. 
Combining \eqref{e:Enu-difference-I1}, \eqref{e:diff-D-Enu-11} and  \eqref{e:diff-D-Enu-2}, we arrive at the desired result.
\end{proof}

%%%%%%%%%%%%%%%%%%%%%%%%%%%%%%%%%%%%%%%%%%%%%%%%%%%%%%%%%%%%%%%%%%%%%%%%%%%%%%%%%%%%%%%%%%%%%%%%%%%%%%%%%%%%%%%%%%%%%%%%%%%%%%%%%%%
%%%%%%%%%%%%%%%%                 Some results related to the function $\F$ %%%%%%%%%%%%%%%%%%%%%%%%%%%%%%%%%%%%%%%%%%%%%%%%%%%%%%%%
%%%%%%%%%%%%%%%%%%%%%%%%%%%%%%%%%%%%%%%%%%%%%%%%%%%%%%%%%%%%%%%%%%%%%%%%%%%%%%%%%%%%%%%%%%%%%%%%%%%%%%%%%%%%%%%%%%%%%%%%%%%%%%%%%%%
\section{Some results related to the function $\F$}\label{s:F}

Suppose that $K:\bH\times\bH\to[0,\infty)$ 
is   invariant under  dilations and horizontal translations. That is, for all $x,y\in\bH$, $a>0$, and $\wt z\in\R^{d-1}$,
\begin{align}\label{e:kernel-K}
	K(x,y)  = K(ax,ay) =  K(x+ (\wt z,0), y+(\wt z,0)) .
\end{align}
Recall that $\e_d=(\wt 0,1) \in \R^d$.  
 By \eqref{e:kernel-K}, for any $x=(\wt x,x_d) \in \bH$ and $y = (\wt y,y_d) \in \bH$,
\begin{align}\label{e:kernel-K-property}
	K(x,y) = K( x/x_d, y/x_d) = K( \e_d, (y-x)/x_d + \e_d ).
\end{align}
Hence, $K(x,y)$ depends on $(x,y)$ only through the quantity  $(y-x)/x_d$.   Define $\F_K:\bH_{-1}\to [0,\infty)$ by  
$$\F_K(y) = K( \e_d, y  + \e_d), \quad y \in \bH_{-1}. $$
Then \eqref{e:kernel-K-property} can be rewritten as 
\begin{align}\label{e:kernel-K-property-equiv}
	K(x,y)  = \F_K((y-x)/x_d) \quad \text{for all $x,y  \in \bH$}.
\end{align}
In particular, if $K(x,y) \asymp \Phi(|x-y|^2/(x_dy_d))$ for $x, y \in \bH$, then
\begin{equation*}
 \F_K(y)  \asymp  \Phi\left(\frac{|y|^2}{y_d+1}\right), \quad\;\; y \in \bH_{-1}.
\end{equation*}

Suppose, in addition, that $K$ is symmetric, that is, $K(x,y) = K(y,x)$ for all $x,y \in\bH$. Then, by  symmetry and  \eqref{e:kernel-K-property-equiv},
\begin{align*}
	\F_K(y) = K(\e_d, y+\e_d) = K(y+\e_d, \e_d) = \F_K( - y/(y_d+1)), \quad \;\; y=(\wt y,y_d) \in \bH_{-1}.
\end{align*}

Motivated by the preceding discussion, we introduce the following definition.
\begin{defn}
	\rm We say that a Borel function $\F:\bH_{-1} \to [0,\infty)$ is a \textit{symmetric kernel profile  with reference function $\Phi$} if the following conditions are satisfied:
	\begin{enumerate}[(i)]
		\item For all $y=(\wt y,y_d) \in \bH_{-1}$,
		\begin{align}\label{e:profile-1}
			\F(y) = \F(-y/(y_d+1)).
		\end{align} 
		
		\item There exists $C>1$ such that  for all $y=(\wt y,y_d) \in \bH_{-1}$,
			\begin{align}\label{e:profile-2}
			C^{-1} \Phi\left(\frac{|y|^2}{y_d+1}\right)\le  \F(y)  \le C \Phi\left(\frac{|y|^2}{y_d+1}\right).
		\end{align} 
	\end{enumerate}
	We use $\mathfrak P_\Phi$ to denote the family of all symmetric kernel profiles with reference function $\Phi$.
\end{defn}

Recall that $\sB(x,y)$ satisfies \eqref{e:J-Phi}.
If  \hyperlink{B2}{{\bf (B2)}} additionally holds, recall also that   $\F:\bH_{-1} \to [0,\infty)$ is defined by \eqref{e:F-symmetrization}.

\begin{lemma}\label{l:F-basic}
	Assume that \eqref{e:J-Phi} and \hyperlink{B2}{{\bf (B2)}} hold, and let $\F$ be defined by \eqref{e:F-symmetrization}. Then the following hold.

\noindent (i)	For any $Q \in \partial D$ and     $x,y \in E^{Q}_{\nu_0}(\wh R/8)$ with $x=(\wt x, x_d)$ in \text{\rm CS}$_Q$,	\begin{align}\label{e:F-basic-1}	
	\begin{split}
		& \big|\sB(x,y)-\sB(x,x) \F((y-x)/x_d) \big|\\		&\le C\left(\frac{|x-y|\vee \delta_D(x)\vee \delta_D(y)}{|x-y|\wedge \delta_D(x)\wedge \delta_D(y)}\right)^{\theta_1} \, 		(|x-y|\vee\delta_D(x)\vee \delta_D(y))^{\theta_2}.
	\end{split}		\end{align} 
	
\noindent (ii) $\F \in \mathfrak P_\Phi$.
\end{lemma}
\begin{proof} 
	(i)  The result follows directly from \hyperlink{B2}{{\bf (B2)}} by the same argument as in \cite[Lemma 6.1]{CKSV23}.
	
	\noindent (ii) Condition \eqref{e:profile-1} is immediate from the definition of $\F$. We show that \eqref{e:profile-2} holds. 
	
	 We fix  $Q\in \partial D$ and use the  coordinate system CS$_{Q}$. 
	 Let $y=(\wt y,y_d)\in \bH_{-1}$.  For each $\eps \in (0,\wh R/8]$, 
	define $y_\eps=\eps(y+\e_d) =(\wt y_\eps, (y_\eps)_d)$. Then we have
	\begin{align}\label{e:y_eps}
	\wt y_\eps = \eps \wt y \quad \text{and} \quad  (y_\eps)_d=\eps(y_d+1).
	\end{align}
Let $\eps_0\in (0,\wh R/16)$ be a constant satisfying $\eps_0 |y| <\wh R /32$, 
	$\eps_0(y_d+1)<\wh R /16$ and  $\eps_0(y_d+1)>2^{2+3\nu}\eps_0^{1+\nu} \wh R^{-\nu} |\wt y|^{1+\nu}$. 
	Then  for all $\eps \in (0,\eps_0]$,  $y_\eps \in E^Q_\nu(\wh R/8)$. Thus, by Lemma \ref{l:distance-in-Enu}(i) and \eqref{e:y_eps}, there exists $c_1=c_1(\nu)>1$ such that 
	\begin{align}\label{e:y_eps-2}
c_1^{-1}\eps (y_d+1) \le 	\delta_D(y_\eps) \le c_1\eps (y_d+1).
	\end{align}
Note that 
 $\delta_D(\eps \e_d) \asymp \eps$ 
 and $|\eps \e_d - y_\eps |=\eps|y|$. 	For all $\eps \in (0,\eps_0]$, using 	\eqref{e:F-basic-1}	 
	and \eqref{e:y_eps-2}, we get
	\begin{align}\label{e:F-B5-1}
		& \big|\sB(\eps \e_d,y_\eps)-\sB(\eps \e_d, \eps \e_d)\F(y) \big|
		= \big|\sB(\eps \e_d,y_\eps)- \sB(\eps \e_d, \eps \e_d) \F((y_\eps - \eps \e_d)/\eps) \big|\nn \\
		& \le c_1 \eps^{\theta_2}\bigg(\frac{ |y| \vee 1 \vee (y_d+1)}{|y| \wedge 1 \wedge (y_d+1)   }\bigg)^{\theta_1} 
		\big(|y| \vee 1 \vee  (y_d+1) \big)^{\theta_2}	=:c_2(y) \eps^{\theta_2}. 
	\end{align}
	
On the other hand, let $\eps_0' \in (0,\eps_0]$ be such that  $\eps_0' \left( |y| \vee (c_1(y_d+1)) \right)<A_0$. Then, by \eqref{e:y_eps-2},  
$ |\eps \e_d - y_\eps|  \vee \delta_D(\eps \e_d) \vee \delta_D(y_\eps) <A_0$ for all $\eps\in (0,\eps_0']$. 
 Using  \eqref{e:y_eps-2} again, we get
$$
\frac{|\eps \e_d - y_\eps|^2}{\delta_D(\eps \e_{d})\delta_D(y_\eps)}\asymp \frac{\eps^2 |y|^2}{\eps^2 (y_d+1)}=\frac{|y|^2}
 {y_d+1}. 
$$
Hence by \eqref{e:J-Phi} and \eqref{e:Phi-scaling} we have that  for all $\eps \in (0,\eps_0']$,
	\begin{align}\label{e:F-B5-2}
		\frac{1}{c_3} \Phi\left(\frac{|y|^2}{y_d+1}\right)\le	\sB(\eps \e_d, y_\eps) \le c_3 \Phi\left(\frac{|y|^2}{y_d+1}\right),
	\end{align}
	where $c_3>1$  is a constant independent of $y$.
	
	Now, by choosing $\eps\in (0,\eps_0']$ small enough so that $c_2(y)\eps^{\theta_2}<\Phi(|y|^2/(y_d+1))/(2c_3)$, we deduce from  
	\eqref{e:F-B5-1},  \eqref{e:F-B5-2} and  the triangle inequality 
	that 
	\begin{align*}
		\frac{1}{2c_3} \Phi\left(\frac{|y|^2}{y_d+1}\right)\le \sB(\eps \e_d, \eps \e_d)\F(y)  		\le 2c_3 \Phi\left(\frac{|y|^2}{y_d+1}\right).
	\end{align*}
	Since $C_0^{-1}\le \sB(\eps \e_d, \eps \e_d) \le C_0$ by \eqref{e:J-Phi},  this yields  \eqref{e:profile-2}. 
\end{proof}

In the remainder of this section, we establish some results related to
$\F\in\mathfrak P_\Phi$. Although we only need these results for the specific
function $\F$ defined in \eqref{e:F-symmetrization} under \hyperlink{B2}{{\bf (B2)}},
we work with the general class $\mathfrak P_\Phi$ to emphasize that our
arguments rely only on \eqref{e:profile-1} and \eqref{e:profile-2}.

\medskip

Let $\F\in\mathfrak P_\Phi$. Recall that for $q \in [(\alpha-1)_+, \alpha-\overline{\beta})$, 
\begin{align*}
	C(\alpha,q,\F) =\int_{\R^{d-1}}\frac{1}{(|\wt u|^2+1)^{(d+\alpha)/2}}\int_0^1  
	\frac{(s^q -1)(1-s^{\alpha-1-q})}{(1-s)^{1+\alpha}} \F\big(((s-1)\wt u, s-1)\big) ds \, d \wt u.
\end{align*}
 
One of our first objectives is to establish the well-definedness of $C(\alpha,q,\F)$.  To this end, we first collect some elementary estimates for the integrand appearing in the definition of $C(\alpha,q,\F)$.

Observe that,  for any $q \in [(\alpha-1)_+, \alpha-\overline{\beta})$,
\begin{align}\label{e:integrand-estimate-1}
	\frac{(s^q -1)(1-s^{\alpha-1-q})}{(1-s)^{1+\alpha}}  \le 
	2^{1+\alpha} s^{\alpha-1-q} \quad \text{for all $s\in (0,1/2)$}.
\end{align}
Moreover, by the mean value theorem, 
\begin{align}\label{e:integrand-estimate-2}
		\frac{(s^q -1)(1-s^{\alpha-1-q})}{(1-s)^{1+\alpha}}   \le q(q-\alpha+1) (2^{1-q} \vee 1) 	2^{ q+2-\alpha} (1-s)^{1-\alpha} \quad \text{for all $s\in [1/2,1)$}.
\end{align}
We also note that, for every $\beta \in [0, (\alpha+1)/2)$,
\begin{align}\label{e:integrand-estimate-3}
	\int_{\R^{d-1}}\frac{d\wt u}{(|\wt u|^2+1)^{(d+\alpha)/2-\beta}} = C\int_0^{\infty} \frac{r^{d-2}}{(r^2+1)^{(d+\alpha)/2-\beta}}\, dr\le c_3 \int_0^\infty ( r^{d-2} \wedge r^{2\beta - \alpha -2} )dr  <\infty.
\end{align}

\begin{lemma}\label{l:constant}	
 For any $\F\in\mathfrak P_\Phi$,   the map	$q\mapsto C(\alpha,q,\F)$ is  well-defined, strictly increasing, and  continuous  on 
	$[(\alpha-1)_+, \alpha-\overline{\beta})$. Moreover, $C(\alpha,(\alpha-1)_+,\F)=0$.
\end{lemma}
\begin{proof}
	Let $\beta\in (\overline{\beta}, 1\wedge \alpha)$ be such that  $q \in [(\alpha-1)_+, \alpha-\beta)$. 
	 Since $q \ge \alpha-1$ and  $\F$ is non-negative, the integrand in the definition of  $C(\alpha,q,\F)$ is non-negative.   Using \eqref{e:profile-2}, the monotonicity of $\Phi$ and \eqref{e:Phi-scaling},  we obtain that for all $\wt u \in \R^{d-1}$ and  $s\in (0,1)$,
	\begin{equation*}
			\F\big(((s-1)\wt u, s-1)\big) \le c_1\Phi\left(\frac{ (|\wt u|^2 + 1)}{s}\right) \le c_2 \left(\frac{ (|\wt u|^2 + 1)}{s}\right) ^\beta. \end{equation*}
	Thus,  
	\begin{align}\label{e:constant-1}
 C(\alpha,q,\F) &\le 		c_2\left( \int_{\R^{d-1}}\frac{d \wt u}{(|\wt u|^2+1)^{(d+\alpha)/2-\beta}}  \right) \left(  
\int_0^1		\frac{(s^q-1)(1-s^{\alpha-1-q})}{s^\beta (1-s)^{1+\alpha}}  ds \right)\nn\\
&\le c_3\int_0^1		\frac{(s^q-1)(1-s^{\alpha-1-q})}{s^\beta (1-s)^{1+\alpha}}  ds.
	\end{align}
Here, the second inequality follows from \eqref{e:integrand-estimate-3}, since  $\beta< 1\wedge \alpha \le (\alpha+1)/2$.

	Let $f(s):=s^{-\beta} (1-s)^{-1-\alpha} (s^q-1)(1-s^{\alpha-1-q}) $ for $s\in (0,1)$.  By \eqref{e:integrand-estimate-1} and \eqref{e:integrand-estimate-2},  
	\begin{align}\label{e:constant-3}
		f(s)  \le  c_4 \begin{cases}
			s^{\alpha-1-q-\beta} & \text{for $s\in (0,1/2)$},\\
			 (1-s)^{1-\alpha}& \text{for $s\in [1/2,1)$}.
		\end{cases}
	\end{align}
	Since $q<\alpha-\beta$ and $\alpha<2$, we have $\alpha-1-q-\beta>-1$ and $1-\alpha>-1$. Hence,  by \eqref{e:constant-3},  $\int_0^1 f(s)ds<\infty$. Therefore, by  \eqref{e:constant-1},  $C(\alpha,q, \F)$ is well-defined.

	Continuity can be proved using the dominated convergence theorem. 	For any fixed $s \in (0,1)$, the map  $g_s(q):= (s^q-1)(1-s^{\alpha-1-q})$ is strictly increasing on 
	$[(\alpha-1)/2,\infty)$ and satisfies  $g_s(0)=g_s(\alpha-1)=0$. Thus, $q\mapsto C(\alpha,q,\F)$ is strictly increasing on 
	$[(\alpha-1)_+,\alpha-\ub)$ and $C(\alpha,(\alpha-1)_+,\F)=0$. The proof is complete.
\end{proof}

In our later construction of barriers for the regional operator 
$L^{\kappa}_{\sB}$
 on
$D$, we will approximate  its nonlocal part 
$$
L^0_\sB u(x)
= \lim_{\eps\to 0}\int_{D,\, |y-x|>\eps}(u(y)-u(x))\frac{\sB(x,y)}{|x-y|^{d+\alpha}}\, dy
$$ 
by the following regional nonlocal operator on $\bH$:
$$
L^0_\F u(x)
:= \lim_{\eps\to 0}\int_{\bH,\, |y-x|>\eps}(u(y)-u(x))\frac{\sB_\F(x,y)}{|x-y|^{d+\alpha}}\, dy, 
$$
where $ \sB_\F(x,y):=\F((y-x)/x_d)$ is the symmetric dilation- and horizontal-translation-invariant kernel on $\bH \times \bH$ associated with the function $\F$ defined in  \eqref{e:F-symmetrization} under \hyperlink{B2}{{\bf (B2)}}. 

A key feature of  
$L^0_\F u$ 
is that its action on  the  monomials $y_d^q$, $q\in[0,\alpha-\ub)$, can be computed explicitly. The following proposition establishes this explicit formula, together with an error estimate for the convergence of the corresponding principal value integral.  Combined with the above approximation, these results will later be used to estimate  
$L^\kappa_\sB u$
for a suitable class of functions $u$.

\begin{prop}\label{p:I2}
	For any $\F\in\mathfrak P_\Phi$ and $q \in [(\alpha-1)_+, \alpha-\overline{\beta})$, 
	there exists $C>0$  such that for all $x=(\wt 0, x_d) \in \bH$ and $\delta \in (0, (x_d \wedge 1)/2)$,
	\begin{align}\label{e:I2-result}	
		\left|	\int_{\bH, \, |x-y|> \delta} (y_d^q-x_d^q)  
		\frac{\F((y-x)/x_d)}{|x-y|^{d+\alpha}}dy -	C(\alpha,q,\F)x_d^{q-\alpha} \right| 
		\le C(\delta/x_d)^{2-\alpha} x_d^{q-\alpha}.	
	\end{align}
\end{prop}

 Before giving the proof of   Proposition \ref{p:I2}, we  show that the integral in \eqref{e:I2-result} is well-defined. 
\begin{lemma}\label{l:int-F}
	If $\F \in \mathfrak P_\Phi$ and $q \in [0, \alpha-\overline{\beta})$, then
	for all  $x=(\wt 0,x_d) \in \bH$ and $\delta>0$, 
	\begin{align*}
		\int_{\bH, \, |x-y|>\delta} \frac{y_d^q\, \F((y-x)/x_d)}{|x-y|^{d+\alpha}} dy<\infty.
	\end{align*}
\end{lemma}
\begin{proof}
Let $x=(\wt 0,x_d)\in \bH$ 
and $\delta>0$. By \eqref{e:profile-2}, we have
\begin{align*}
	&	\int_{\bH, \, |x-y|>\delta} \frac{y_d^q\, \F((y-x)/x_d)}{|x-y|^{d+\alpha}} dy \le c_1 	\int_{\bH, \, |x-y|>\delta} \frac{y_d^q}{|x-y|^{d+\alpha}} \Phi \left( \frac{|x-y|^2}{x_dy_d}\right) dy \\
	&= c_1 \bigg(  	\int_{\bH, \, |x-y|>\delta, \, |x-y|^2\le x_dy_d}  +\int_{\bH, \, |x-y|>\delta, \, |x-y|^2> x_dy_d} \bigg) \frac{y_d^q}{|x-y|^{d+\alpha}} \Phi \left( \frac{|x-y|^2}{x_dy_d}\right) dy  \\
	&=:c_1(I_1+I_2).
\end{align*}
We first note that, with $\wh{x}:=(\wt{0}, (3/2)x_d)$, 
 \begin{align}\label{e:int-F-1}
  \{y\in \bH:  |x-y|^2<x_d y_d\} 
	= B( \wh x, (\sqrt 5/2)x_d).
\end{align}
	By the monotonicity of $\Phi$ and \eqref{e:int-F-1}, 
\begin{align*}
	I_1\le \Phi(1) \int_{ B( \wh x, (\sqrt 5/2)x_d), \, |x-y|>\delta} \frac{y_d^q}{|x-y|^{d+\alpha}} dy \le \frac{\Phi(1)((3+\sqrt 5)/2)^qx_d^q}{\delta^{d+\alpha}}\int_{B( \wh x, (\sqrt 5/2)x_d)} dy <\infty.
\end{align*}

To estimate $I_2$, let $\beta  \in  (\ub, 1 )$ be such that $q<\alpha-\beta$. 
We see from \eqref{e:int-F-1} that,  if  $|x-y|^2 > x_d y_d$, then $|x-y|\ge |\wh{x}-y|-|\wh{x}-x|\ge \sqrt{5}x_d/2-x_d/2\ge x_d/2$. 
Hence, using \eqref{e:Phi-scaling}, we obtain 
\begin{align}\label{e:int-F-2}
	I_2&\le \sum_{n=0}^\infty \int_{\bH \cap B(x, 2^nx_d) \setminus B(x,2^{n-1}x_d)}
	\frac{y_d^q}{|x-y|^{d+\alpha}} \Phi\left(\frac{|x-y|^2}{x_d y_d}\right)\, dy\nn\\
	&\le c_2x_d^{-\beta}\sum_{n=0}^\infty \int_{\bH \cap B(x, 2^nx_d) \setminus B(x,2^{n-1}x_d)}
	\frac{y_d^{q-\beta}}{|x-y|^{d+\alpha-2\beta}} dy \nn\\
		&\le 2^{d+\alpha-\beta}c_2x_d^{-d-\alpha+\beta}\sum_{n=0}^\infty 2^{-n(d+\alpha-2\beta)}\int_{\bH \cap B(x, 2^nx_d) }
y_d^{q-\beta} dy .
\end{align}
Since $q-\beta\ge -\beta>-1$, for every $n\ge 0$, we have
\begin{align}\label{e:int-F-3}
&	\int_{\bH \cap B(x, 2^nx_d) }
	y_d^{q-\beta} dy  \le \int_{\wt y \in \R^{d-1}, \, |\wt y|<2^nx_d} d\wt y \int_0^{(2^n +1)x_d} y_d^{q-\beta}dy_d\nn\\
	& = c_3 (2^nx_d)^{d-1} (2^n+1)^{q-\beta+1} x_d^{q-\beta+1} \le c_4 2^{n(d+q-\beta)} x_d^{d +q-\beta},
\end{align}
where we used $2^n+1\le 2^{n+1}$ in the last inequality. Combining \eqref{e:int-F-2} with \eqref{e:int-F-3}, since $q<\alpha-\beta$, we arrive at
\begin{align*}
	I_2\le  c_5 x_d^{q-\alpha} \sum_{n=0}^\infty 2^{-n(\alpha-q-\beta)} = c_6x_d^{q-\alpha}<\infty.
\end{align*}
This completes the proof.
\end{proof}

\begin{proof}[Proof of Proposition \ref{p:I2}] 
	Let 
	 $x =(\wt 0,x_d) \in \bH$ and $\delta \in (0,(x_d \wedge 1)/2)$. We set
	$$
	I(q,\delta):=\int_{\bH, \, |x-y|>\delta}  (y_d^q - x_d^q)\frac{ \F((y-x)/x_d)}{|x-y|^{d+\alpha}}dy.
	$$
	Using Lemma \ref{l:int-F} twice (with $q$ and $q=0$), one sees that the integrand in the above integral is absolutely integrable. Hence $I(q,\delta)$ is well-defined. Using the change of variables $z=y/x_d$ in the first equality below and  the change of  
	variables $\wt u = \wt z/(z_{d}-1)$ in the second, we get
	\begin{align}\label{e:generator-F-0}
		I(q,\delta)	&	=x_d^{q-\alpha} \int_{\bH,\, |\wt{z}|^2+(z_d-1)^2> (\delta/x_d)^2}
		\frac{(z_d^q-1)\,\F(z- \e_d)}{| (\wt{z}, z_d)-\mathbf{e}_d|^{d+\alpha}} d \wt{z}\,	dz_d \nn\\
		&=x_d^{q-\alpha}\int_{\bH,\, (|\wt{u}|^2 +1)(z_d-1)^2> (\delta/x_d)^2}
		\frac{(z_d^q-1)\,\F(((z_d-1)\wt u, z_d-1))}{(|\wt u|^2+1)^{(d+\alpha)/2} |z_d-1|^{1+\alpha}} d \wt u\,	dz_d.
	\end{align}
	Define  $\xi(\delta, \wt u):=(\delta/x_d)(|\wt{u}|^2+1)^{-1/2}$ for $\wt u \in \R^{d-1}$.  Then
	\begin{align}\label{e:generator-F-1}
		\xi(\delta, \wt u) \le \delta/x_d < 1/2 \quad \text{for all $\wt u \in \R^{d-1}$}.
	\end{align}
	By Fubini's theorem, we obtain from \eqref{e:generator-F-0} that
	\begin{align}\label{e:generator-F-2}
		I(q,\delta)	=x_d^{q-\alpha}\int_{\R^{d-1}} ( I_1(q,\delta,\wt u)  + I_2(q,\delta,\wt u) )
		\frac{d\wt u}{(|\wt u|^2+1)^{(d+\alpha)/2}},
	\end{align}
	where
	\begin{align*}
		I_1(q,\delta, \wt u)&:= \int_0^{1-	\xi(\delta, \wt u) } 
		\frac{(z_d^q -1)\, \F(((z_d-1)\wt u,z_d-1))}{|z_d-1|^{1+\alpha}}  dz_d,\\
		I_2(q,\delta, \wt u)&:= \int_{1+	\xi(\delta, \wt u) }^\infty 
		\frac{(z_d^q -1)\, \F(((z_d-1)\wt u,z_d-1))}{|z_d-1|^{1+\alpha}}  dz_d.
	\end{align*}
	By \eqref{e:profile-1}, for all $\wt u \in \R^{d-1}$ and $s>0$,
	\begin{align}\label{e:generator-F-3}
	 \F(((s-1)\wt u,s-1))  = \F((  (1/s-1)\wt u,  1/s-1)).
	\end{align}
	Using the change of the variables $s=1/z_d$ and \eqref{e:generator-F-3}, 
	we see that
	\begin{align*}
		I_2(q,\delta, \wt u)&= \int_0^{(1+\xi(\delta, \wt u))^{-1}}  
		\frac{(s^{-q}-1)}{(s^{-1}-1)^{1+\alpha}} \F(((1/s-1)\wt u, 1/s-1))\, \frac{ds}{s^2}\\
		&=\bigg(\int_0^{1-\xi(\delta, \wt u) }+ \int_{1-\xi(\delta, \wt u) }^{(1+\xi(\delta, \wt u))^{-1}} \bigg) 
		\frac{s^{\alpha-1-q}(1-s^q)}{(1-s)^{1+\alpha}} \F(((s-1)\wt u, s-1)) ds\\
		&=:I_{2,1}(q,\delta,\wt u) + I_{2,2}(q,\delta,\wt u).
	\end{align*}

	Let $\beta\in (\overline{\beta}, 1\wedge \alpha)$. 
Since $\xi(\delta, \wt u) < 1/2$, we have $(s-1)^2/s \le 1$ for all $s \in (1-\xi(\delta, \wt u), 1)$. Thus,
 using \eqref{e:profile-2}, the monotonicity of $\Phi$ and \eqref{e:Phi-scaling}, we get that for all $s \in (1-\xi(\delta, \wt u), 1)$,
	\begin{align} \label{e:Fep}
		\F(((s-1)\wt u, s-1)) &\le c_1 \Phi \left( \frac{(s-1)^2 (|\wt u|^2 + 1)}{s} \right)  \le c_1 \Phi ( |\wt u|^2 + 1 ) \le c_2 \Phi(1)  ( |\wt u|^2 + 1 ) ^\beta.
	\end{align}	
	Using  \eqref{e:Fep},
	\eqref{e:integrand-estimate-2},  \eqref{e:integrand-estimate-3} (with $\beta <1\wedge \alpha \le (\alpha+1)/2$) and  \eqref{e:generator-F-1}, we obtain
	\begin{align}\label{e:generator-F-4}	
		& \left| \int_{\R^{d-1}}(I_1(q,\delta, \wt u) + I_{2,1}(q,\delta, \wt u)) \frac{d\wt u}{(|\wt u|^2+1)^{(d+\alpha)/2}} 
		- C(\alpha,q,\F)\right| \nn\\	
		&= \int_{\R^{d-1}}\frac{1}{(|\wt u|^2+1)^{(d+\alpha)/2}}\int_{1-\xi(\delta, \wt u)}^1  
		\frac{(s^q -1)(1-s^{\alpha-1-q})}{(1-s)^{1+\alpha}} \F(((s-1)\wt u, s-1)) ds \, d \wt u\nn\\
		&\le c_3\int_{\R^{d-1}}\frac{d\wt u}{(|\wt u|^2+1)^{(d+\alpha)/2 -\beta}}\int_{1-\xi(\delta, \wt u)}^1   (1-s)^{1-\alpha} ds  
		  \le  c_4(\delta/x_d)^{2-\alpha}.
	\end{align}
	On the other hand, we note that, by \eqref{e:generator-F-1},  for every $\wt u \in \R^{d-1}$,
	\begin{align}\label{e:generator-F-5-1}	  (1+\xi(\delta, \wt u))^{-1}-(1-\xi(\delta, \wt u) )
	= 	\xi(\delta, \wt u)^2 (1+\xi(\delta, \wt u))^{-1} \le \xi(\delta, \wt u)^2 \le (\delta/x_d)^2.
	\end{align} 
	Further, using \eqref{e:Fep}, the mean value theorem and the fact that $\xi(\delta, \wt u)\le 1/2$, we see that for every $\wt u \in \R^{d-1}$ and $s\in (1-\xi(\delta,\wt u),1)$,
	\begin{align}\label{e:generator-F-5-2}	 
	\frac{s^{\alpha-1-q}(1-s^q)}{(1-s)^{1+\alpha}} \F(((s-1)\wt u, s-1))  \le \frac{qc_5 s^{\alpha-1-q} (s^{q-1} \vee 1)(|\wt u|^2+1)^\beta }{(1-s)^{\alpha}} \le \frac{c_6(|\wt u|^2+1)^\beta }{(1-s)^\alpha}.
	\end{align}
	Using \eqref{e:generator-F-5-1} and \eqref{e:generator-F-5-2} in the first inequality below, \eqref{e:generator-F-1} in the second and third, and  \eqref{e:integrand-estimate-3} in the fourth, we obtain that 
	\begin{align}\label{e:generator-F-6}
		  \int_{\R^{d-1}}  \frac{I_{2,2}(q,\delta, \wt u)}{(|\wt u|^2+1)^{(d+\alpha)/2}}d\wt u
		 &\le c_6\int_{\R^{d-1}}  \frac{1}{(|\wt u|^2+1)^{(d+\alpha)/2-\beta}} \int_{1-\xi(\delta, \wt u)}^{1-\xi(\delta, \wt u)+\xi(\delta, \wt u)^2} 
		 \frac{ds}{(1-s)^{\alpha}}\,d\wt u\nn\\
		 &	\le c_6 \int_{\R^{d-1}}  \frac{1}{(|\wt u|^2+1)^{(d+\alpha)/2-\beta}} \frac{\xi(\delta, \wt u)^2}{(\xi(\delta, \wt u)/2)^\alpha} d\wt u \nn\\
		 &\le  2^\alpha c_6  (\delta/x_d)^{2-\alpha} \int_{\R^{d-1}}  \frac{ d\wt u}{(|\wt u|^2+1)^{(d+\alpha)/2-\beta}}\le  c_7 (\delta/x_d)^{2-\alpha}.
		 \end{align}
	Combining \eqref{e:generator-F-2} with \eqref{e:generator-F-4} and \eqref{e:generator-F-6}, we arrive at the result.
\end{proof}

%%%%%%%%%%%%%%%%%%%%%%%%%%%%%%%%%%%%%%%%%%%%%%%%%%%%%%%%%%%%%%%%%%%%%%%%%%%%%%%%%%%%%%%%%%%%%%%%%%%%%%%%%%%%%%%%%%%%%%%%%%%%%%%%%%%%%
%%%%%%%%%%%%                  Dynkin-type formula and analysis of $X$                 %%%%%%%%%%%%%%%%%%%%%%%%%%%%%%%%%%%%%%%%%%%%%%%
%%%%%%%%%%%%%%%%%%%%%%%%%%%%%%%%%%%%%%%%%%%%%%%%%%%%%%%%%%%%%%%%%%%%%%%%%%%%%%%%%%%%%%%%%%%%%%%%%%%%%%%%%%%%%%%%%%%%%%%%%%%%%%%%%%%%%
\section{Analysis of $\overline X$ and $X$}\label{s:dynkin-X}
 
In this section, we assume that $D$ is a $C^{1,\eta}$ open set,  and that \hyperlink{A}{{\bf (A)}} and \hyperlink{K}{{\bf (K)}} hold. Recall that $X$  denotes either $X^0$ or $X^{\kappa}$.  

We first recall the L\'evy system formula together with its consequences, and  establish several results for  $\overline X$ based on its heat kernel estimates. We then  derive 
a priori  bounds for the heat kernel and Green function of $X$. We remark that all the results in this section remain valid under   
the weaker assumption that $D$ is a Lipschitz open set.
However, for the sake of brevity, we work throughout in the $C^{1,\eta}$ setting.

We recall that $m_d$ denotes the $d$-dimensional Lebesgue measure.
  Since the jump measure $\sB(x,y)|x-y|^{-d-\alpha}dxdy$ of $(\sE^0, \overline \sF)$ is absolutely continuous with respect to $m_d \otimes m_d$,  by using \cite[(5.3.15)]{FOT} and repeating the arguments in \cite[p.~40]{CK03}, one sees that  $\overline X$ satisfies the following L\'evy system formula: 
 For any $x \in  \overline D$, any non-negative Borel function $f$ on 
$\overline D \times \overline D$
vanishing on the diagonal, and any stopping time $\tau$,
\begin{align}\label{e:Levysystem-X-bar}
	\E_x\bigg[  \sum_{s \le \tau} f( \overline X_{s-},  \overline X_s)  \bigg] 
	= \E_x \left[ \int_0^\tau   \int_{D} \frac{f( \overline X_s, y) 
		\sB( \overline X_s, y)}{| \overline X_s-y|^{d+\alpha}} dy  \, ds \right].
\end{align}
Moreover, since 	the jump measure $\sB(x,y)|x-y|^{-d-\alpha}dxdy$ and the killing measure $\kappa(x)dx$ of  $(\sE^\kappa,  \sF^\kappa)$ are absolutely continuous with respect to  $m_d \otimes m_d$ and  $m_d$ respectively,  
$X$ satisfies the following L\'evy system formula (cf. \eqref{e:Levysystem-X-bar}): For any $x \in  D $, any non-negative Borel function $f$ on $D \times D_\partial$ vanishing on the diagonal, and any stopping time $\tau$,	
\begin{align}\label{e:Levysystem-X-kappa}	
\E_x\bigg[  \sum_{s \le \tau} f( X_{s-},  X_s)  \bigg] 	= \E_x \left[ \int_0^\tau \bigg(  \int_{D} 	\frac{f( X_s, y) \sB( X_s, y)}{| X_s-y|^{d+\alpha}} dy  	+ \kappa(X_s)f(X_s, \partial) \bigg) \,ds  \right].	
\end{align}
 When $X$ stands for $X^0$, the function $\kappa$ in \eqref{e:Levysystem-X-kappa} is understood to be zero.

For an open set $U \subset  D$,  let $\tau_U:=\inf\{t>0: X_t \notin U\}$.   We  denote by $X^{U}$ 
the part of the process $X$ killed upon exiting $U$.

The following is a standard consequence of \eqref{e:Levysystem-X-kappa} and \eqref{e:J-Phi}.

\begin{lemma}\label{l:levy-consequence}
(i) There exists $C>0$ such that	for any open set $U\subset D$ and Borel set $V \subset D$ with $\inf_{z\in U, \, y \in V} |z-y|\ge r >0$,
	\begin{align*}
		\P_x \left( X_{\tau_U} \in V \right) \le \frac{C}{r^\alpha} \E_x \left[\int_0^{\tau_U} \Phi\left(\frac{r}{\delta_D(X_s)} \right)ds \right] \quad \text{for all $x\in U$}.
	\end{align*}
	
\noindent	(ii) There exists $C>0$ such that	for any bounded open set $U\subset D$ and  bounded Borel set $V \subset D$ with  $r:=\inf_{z\in U, \, y \in V} |z-y|>0$,
	\begin{align*}
		\P_x \left( X_{\tau_U} \in V \right) \ge \frac{Cm_d(V)}{\diam(U \cup V)^{d+\alpha}} \E_x \left[\int_0^{\tau_U} \Phi\left(\frac{r \wedge A_0}{\delta_D(X_s)} \right)ds 
		\right], \quad \text{for all $x\in U$}.
	\end{align*}
\end{lemma}
\begin{proof}
By \eqref{e:Levysystem-X-kappa},  we get that for all $x\in U$,
\begin{align}\label{e:levy-1}
	&	\P_x \left( X_{\tau_U} \in V \right) =\E_x \left[\int_0^{\tau_U} \int_V \frac{\sB(X_s,y)}{|X_s-y|^{d+\alpha}}dyds \right].
\end{align}

\noindent (i) 	Combining \eqref{e:levy-1}  with \eqref{e:jump-tail-estimate}, we obtain
	\begin{align*}
	&	\P_x \left( X_{\tau_U} \in V \right) \le \E_x \left[\int_0^{\tau_U} \int_{D\setminus B(X_s, r)} \frac{\sB(X_s,y)}{|X_s-y|^{d+\alpha}}dyds \right] \le \frac{c_1}{r^\alpha} \E_x \left[\int_0^{\tau_U} \Phi\left(\frac{r}{\delta_D(X_s)} \right)ds \right].
	\end{align*}

\noindent (ii) For all $z\in U$ and $y\in V$,  by \eqref{e:Phi-scaling}, we have
\begin{align*}
	 \Phi\left( \frac{(|z-y|\wedge A_0)^2}{(\delta_D(z)\wedge A_0) (\delta_D(y) \wedge A_0)}\right) &\ge  \Phi\left( \frac{(|z-y|\wedge A_0)^2}{(\delta_D(z)\wedge A_0) ((\delta_D(z)+|z-y|) \wedge A_0)}\right)\\
	 &\ge c_1\begin{cases}
	 	\Phi\left( \frac{|z-y|\wedge A_0}{\delta_D(z)\wedge A_0}\right) &\mbox{ if $\delta_D(z)\le |z-y|$},\\
	 	1&\mbox{ if $\delta_D(z)> |z-y|$}
	 \end{cases}\\
	 &\ge c_2\Phi\left( \frac{|z-y|\wedge A_0}{\delta_D(z)\wedge A_0}\right) \ge c_2\Phi\left( \frac{r\wedge A_0}{\delta_D(z)}\right). \end{align*}
	 Hence, by \eqref{e:J-Phi}, for all $z\in U$ and $y\in V$,
\begin{align*}
	\frac{\sB(z,y)}{|z-y|^{d+\alpha}} \ge \frac{C_0^{-1}c_2}{\diam(U \cup V)^{d+\alpha}}
	\Phi\left( \frac{r\wedge A_0}{\delta_D(z)}\right).
\end{align*}

Combining this with \eqref{e:levy-1} yields the result.
\end{proof}

The following lemma follows from Lemma \ref{l:levy-consequence} and will be used to establish the interior parabolic H\"older regularity for $X$.
\begin{lemma}\label{l:phr-1}
	There exists $C>0$ such that for all $x \in \overline D$, $R\in (0, \delta_D(x)/2]$ and $r\in (0,R/2]$,
	\begin{align*}
		\sup_{y\in B(x,r/2)} \P_y \big( X_{\tau_{B(x,r)}} \in D\setminus B(x,R) \big) \le C (r/R)^\alpha.
	\end{align*} 
\end{lemma}
\begin{proof}
	Let $x\in \overline D$, $R\in (0,\delta_D(x)/2]$, $r\in (0,R/2]$ and  $y \in B(x,r/2)$. By Lemma \ref{l:levy-consequence}(i),  since $\delta_D(w) \ge \delta_D(x) - r \ge R$ for all $w\in B(x,r)$, we obtain
	\begin{align}\label{e:phr-1-1}
		&\P_y \big( X_{\tau_{B(x,r)}} \in D\setminus B(x,R) \big)  \nn\\
		&\le \frac{c_1}{(R-r)^\alpha}\E_y \left[\int_0^{\tau_{B(x,r)}} \Phi\left(\frac{R-r}{\delta_D(X_s)} \right)ds \right]  \le \frac{c_1 \Phi(1) \E_y[\tau_{B(x,r)}]}{(R/2)^\alpha}.
	\end{align}
	 Choose any $z\in \partial B(x,2r)$.  Note that $B(z,r/2) \subset B(x,2R) \subset D$. Using Lemma \ref{l:levy-consequence}(ii) and $\Phi\ge 1$, we obtain
	\begin{align*}
		1 \ge \P_y ( X_{\tau_{B(x,r)}} \in B(z,r/2)) \ge \frac{c_2}{r^\alpha} \E_y \left[\int_0^{\tau_{B(x,r)}} \Phi\left(\frac{(r/2) \wedge A_0}{\delta_D(X_s)} \right)ds \right] \ge \frac{c_2}{r^\alpha} \E_y[\tau_{B(x,r)}].
	\end{align*}
	Combining this with \eqref{e:phr-1-1} yields the desired result.
\end{proof}

\subsection{Analysis of $\overline X$}\label{ss:analysis-X}

\begin{lemma}\label{l:p-estimate}
	For any $T>0$,  there exists $C=C(T)>0$  such that 
	 for all $t\in (0,T]$, $x\in \overline D$ and $r>0$,
	\begin{align}\label{e:p-estimate-result}
		\int_{D\setminus B(x,r)}	\overline p(t, x,y)dy \le \frac{Ct}{r^{\alpha}} \Phi\left( \frac{r \wedge A_0}{(\delta_D(x) \vee t^{1/\alpha})\wedge A_0}\right) \le \frac{Ct}{r^{\alpha}} \Phi\left( \frac{r}{ t^{1/\alpha}}\right).
	\end{align}
\end{lemma}
\begin{proof}
Let $t\in (0,T]$, $x\in \overline D$ and $r>0$. Using \eqref{e:main-HKE1}, the monotonicity of $\Phi$  and 	Lemma \ref{l:p-estimate0}(ii) (with $\theta=\alpha>\ub$), we obtain
\begin{align*}
		\int_{D\setminus B(x,r)}	\overline p(t, x,y)dy&\le c_1t	\int_{D\setminus B(x,r)} \frac{ 1}{|x-y|^{d+\alpha}}\Phi\left( \frac{(|x-y| \wedge A_0)^2}{((\delta_ D(x) \vee t^{1/\alpha}) \wedge A_0)(\delta_D(y)  \wedge A_0) }\right)dy\\
	&\le \frac{c_2t}{r^\alpha} \Phi\left( \frac{r \wedge A_0}{(\delta_ D(x) \vee t^{1/\alpha}) \wedge A_0 }\right),
\end{align*}which proves the first inequality in \eqref{e:p-estimate-result}. The second inequality follows from \eqref{e:without-A0} and the 
monotonicity 
of $\Phi$.
\end{proof}

For any relatively open set $U$ of  $\overline D$, we let 
$
\overline \tau_U:=\inf\{t>0:\overline X_t \notin U\}.
$ 
From Lemma \ref{l:p-estimate}, we obtain the following tightness result for $\overline X_t$, which constitutes one of the key ingredients in the analysis of  $\overline X$.
\begin{lemma}\label{l:EP}
For any $T>0$,	there exists  $C=C(T)>0$ such that 
	\begin{align}\label{e:EP-1}
		\P_{x}( \overline \tau_{B_{\overline D}(x,r)} \le t ) \le \frac{Ct}{r^\alpha} \Phi\bigg(\frac{r}{t^{1/\alpha}}\bigg) \quad \text{for all $t\in (0,T]$, $x \in \overline D$ and $r>0$} .
	\end{align}
	Consequently, for any $T>0$ and $\lambda \in (0, 1-\ub/\alpha)$, 	there exists  $C=C(T,\lambda)>0$ such that 
	\begin{align}\label{e:EP-2}
		\P_{x}( \overline \tau_{B_{\overline D}(x,r)} \le t ) \le  C \bigg(\frac{t}{r^\alpha}\bigg)^{\lambda} \quad \text{for all $t\in (0,T]$, $x \in \overline D$ and $r>0$.} 
	\end{align}
\end{lemma}
\begin{proof}
		Since $(\sE^0,\overline \sF)$ is conservative, 
	\eqref{e:EP-1} follows from  
	Lemma \ref{l:p-estimate}   by a standard argument; see, e.g.,  the proof of \cite[Lemma 2.7]{CKW21}.
	We omit the details.  
	\eqref{e:EP-2} follows from \eqref{e:EP-1} and \eqref{e:Phi-scaling}. 
\end{proof}

The next proposition was established in \cite[Proposition 3.9(i)]{CKSV25}.

\begin{prop}\label{p:ndl0}
There exists a constant $R_1 \in (0, \wh R]$ such that for all $x_0 \in {\overline D}$ and $r\in (0, 	R_1]$, 
	the semigroup $(P^{B_{\overline D}(x_0,r)}_{t})_{t\ge 0}$ has a jointly continuous heat kernel $\overline p_{B_{\overline D}(x_0,r)}(t,x,y)$. Moreover, there exist constants $\delta_1,\delta_2\in (0,1)$ and $C>0$ independent of $x_0$ and $r$ such that
	\begin{align*}
		\overline  p_{B_{\overline D}(x_0,r)}(t,x,y) \ge C t^{-d/\alpha} \quad \text{for all $t\le  \delta_1r^\alpha$ and $x,y \in 
			B_{\overline D}(x_0, \delta_2 t^{1/\alpha})$
		}.
	\end{align*}
\end{prop}

Since $D$ is a $C^{1,\eta}$ open set, it satisfies \textit{the measure density condition}. That is, there exists  $C>0$ such that 
\begin{align}\label{e:VD}			m_d(	B_{ D}(x,r)) \ge Cr^d \quad \text{for all} \;\, x \in \overline D, \;\, r\in (0,\wh R].	\end{align}

In the next proposition, we extend the range of the Dirichlet heat kernel estimates established above.

	 \begin{prop}\label{p:ndl}
For any $R>0$ and $b\in (0,1)$, there exists  $C=C(R,b)>0$  such that for any $x_0 \in \overline D$, $r\in (0,R]$ and 
 $t\in (0, (b r)^\alpha]$, it holds that
	\begin{equation*}
		\overline	p_{	B_{\overline D}(x_0, r)}(t,z,y)  \ge  C t^{-d/\alpha} 
		\quad  \text{for all } y,z \in 	B_{\overline D}(x_0, b t^{1/\alpha}).
	\end{equation*}
\end{prop}
	\begin{proof}
Using 
 \eqref{e:Levysystem-X-bar}, the fact that $\sB(x,y) \ge C_0^{-1}$, 
 \eqref{e:VD},  
\eqref{e:EP-2} (instead of \cite[Lemma 4.3]{CKSV23}) and Proposition \ref{p:ndl0} (instead of \cite[(4.20)]{CKSV23}), 
 one can follow the proof of \cite[Proposition 4.5]{CKSV23} line by line to get the result. We omit the details.
	\end{proof}

For $\lambda\ge 0$, let $\overline G_\lambda$ denote the $\lambda$-resolvent of $(\sE^0,\overline \sF)$. Then $\overline G_\lambda$ admits a kernel $$\overline G_\lambda(x,y):=\int_0^\infty e^{-\lambda t} \overline p(t,x,y) dt, \quad x,y \in \overline D.$$

Define  for $x,y \in \overline D$,
\begin{align}\label{e:definition-green-function}
\overline	u_{d,\alpha}(x,y):= \begin{cases}
		|x-y|^{\alpha-d} &\mbox{ if $d>\alpha$},\\[3pt]
	\displaystyle	\log\left(e+|x-y|^{-1}\right) &\mbox{ if $d=1=\alpha$},\\[3pt]
1 &\mbox{ if $d=1<\alpha$}
	\end{cases}
\end{align}
and
\begin{align*}	\overline	g_{d,\alpha,\Phi}(x,y):=	\int_0^1  \left[t^{-d/\alpha} \wedge  \left( \frac{ t }{|x-y|^{d+\alpha}}\Phi\left( \frac{(|x-y| \wedge A_0)^2}{((\delta_D(x) \vee t^{1/\alpha}) \wedge A_0)((\delta_D(y) \vee t^{1/\alpha}) \wedge A_0) }\right) \right)\right] dt.
\end{align*}
Note that, for any dimension $d$, it holds that \begin{align}\label{e:green-lower}\overline	u_{d,\alpha}(x,y) \ge  |x-y|^{\alpha-d} \ind_{|x-y|\le 1} + (|x-y|^{\alpha-d} \wedge 1)\ind_{|x-y|>1}   \quad \text{for all $x,y \in \overline D$}.\end{align}
Moreover, for any bounded subset $K$ of $\overline D$, 
\begin{align}\label{e:green-lower-inf}
	\inf_{x,y \in K} \overline u_{d,\alpha}(x,y)  \ge \diam(K)^{\alpha-d}\ind_{d>\alpha} + \ind_{d\le \alpha} >0.
\end{align}

The following lemma will be used to derive estimates for  
$\overline G_1(x,y)$ and the Green functions of $X$.
\begin{lemma}\label{l:green-prelim}
	There exists $C>0$ such that
	\begin{align}\label{e:green-prelim-upper}
	\overline	g_{d,\alpha,\Phi}(x,y) \le C \overline u_{d,\alpha}(x,y) \quad \text{for all $x,y \in \overline D$}.
	\end{align}
	Moreover, for any  bounded subset $K$ of $\overline D$, there exists $C'=C'(\diam(K))>0$ such that
	\begin{align}\label{e:green-prelim-lower}
	\overline	g_{d,\alpha,\Phi}(x,y) \ge C'\overline u_{d,\alpha}(x,y)  \quad \text{for all $x,y \in K$}.
	\end{align}
\end{lemma}
\begin{proof}
	Using the monotonicity of $\Phi$ and \eqref{e:without-A0}, 	we see that for all $x,y \in \overline D$,
	\begin{align}\label{e:green-prelim-1}
\overline	g_{d,\alpha,\Phi}(x,y) 	& \le 
		\int_{|x-y|^\alpha \wedge 1}^{1}  t^{-d/\alpha} dt\nn\\
		&\quad +
		\int_0^{|x-y|^\alpha  \wedge 1} 
		\frac{ t }{|x-y|^{d+\alpha}}\Phi\left( \frac{(|x-y| \wedge A_0)^2}{((\delta_D(x) \vee t^{1/\alpha}) \wedge A_0)((\delta_D(y) \vee t^{1/\alpha}) \wedge A_0) }\right) dt\nn\\
		&\le 
	 c_1 \overline u_{d,\alpha}(x,y)+
		\int_0^{|x-y|^\alpha  \wedge 1} 
		\frac{ t }{|x-y|^{d+\alpha}}\Phi\left( \frac{(|x-y| \wedge A_0)^2}{( t^{1/\alpha} \wedge A_0)^2 }\right) dt\nn\\
		&\le  c_1\overline  u_{d,\alpha}(x,y)+
		\int_0^{ |x-y|^\alpha \wedge 1} 
		\frac{ t }{|x-y|^{d+\alpha}}\Phi\left( \frac{|x-y|^2}{t^{2/\alpha}}\right) dt.
	\end{align}
	Let $\beta \in (\ub,\alpha  \wedge 1 )$.  Since $\beta<\alpha$ and  $2\beta-d-\alpha \le 2\beta-1-\alpha<0$, 
	using \eqref{e:Phi-scaling}  and \eqref{e:green-lower},   we obtain that
	\begin{align*}
		&	\int_0^{ |x-y|^\alpha \wedge 1} 
		\frac{ t }{|x-y|^{d+\alpha}}\Phi\left( \frac{|x-y|^2}{t^{2/\alpha}}\right) dt\le  \frac{c_2}{|x-y|^{d+\alpha-2\beta}}	\int_0^{|x-y|^\alpha \wedge 1}  t^{1-2\beta/\alpha} dt  \nn\\
		&=   c_3 \left( |x-y|^{\alpha-d} \ind_{|x-y|\le 1} + |x-y|^{2\beta-d-\alpha} \ind_{|x-y|>1}\right) \nn\\
		&\le   c_3 \left( |x-y|^{\alpha-d} \ind_{|x-y|\le 1} + (|x-y|^{\alpha-d} \wedge 1)\ind_{|x-y|>1}\right) \le c_3 \overline u_{d,\alpha}(x,y). 
	\end{align*}
	Combining this with \eqref{e:green-prelim-1},
	we obtain  \eqref{e:green-prelim-upper}.
	
	To prove \eqref{e:green-prelim-lower}, we first note that, for all $x,y \in K$ and  $t\in [(|x-y|/(2\diam(K)))^\alpha,1]$, 
	\begin{align*}
		&t^{-d/\alpha} \wedge  \left( 	\frac{ t }{|x-y|^{d+\alpha}}\Phi\left( \frac{(|x-y| \wedge A_0)^2}{((\delta_D(x) \vee t^{1/\alpha}) \wedge A_0)((\delta_D(y) \vee t^{1/\alpha}) \wedge A_0) }\right)\right)\\
		& \ge t^{-d/\alpha} \wedge  \frac{t}{|x-y|^{d+\alpha}} \ge  
		\frac{t^{-d/\alpha}}{(1\vee (2\diam(K)))^{d+\alpha}}.
	\end{align*}
	It follows that for all $x,y \in K$,
	\begin{align*}
	\overline	g_{d,\alpha,\Phi}(x,y)  \ge  \frac{1}{(1\vee (2\diam(K)))^{d+\alpha}} \int_{(|x-y|/(2\diam(K)))^\alpha}^1 t^{-d/\alpha} dt \ge c_4 \overline u_{d,\alpha}(x,y).
	\end{align*}
	The proof is complete.
\end{proof}

\begin{lemma}\label{l:resolvent-estimate}
For any bounded subset $K$ of $\overline D$, there exists $C=C(\diam(K))>1$ such that
	\begin{align*}
	 C^{-1}\overline u_{d,\alpha}(x,y) \le 	\overline G_1(x,y) \le C\overline u_{d,\alpha}(x,y) \quad \text{for all $x,y\in K$}.
	\end{align*}
\end{lemma}
\begin{proof}
Let $K\subset \overline D$ be a bounded subset. 
	Using \eqref{e:main-HKE1}, \eqref{e:green-prelim-upper} 
 and \eqref{e:green-prelim-lower}, we see that for all  $x,y \in K$,
\begin{align*}
	\overline G_{1}(x,y) 
	\le c_1\overline  g_{d,\alpha,\Phi}(x,y) + \int_1^\infty e^{- t}\overline p(t,x,y)dt\le c_2\overline u_{d,\alpha}(x,y)+ \int_1^\infty e^{- t}\overline p(t,x,y)dt
\end{align*}
and $
\overline G_{1}(x,y) 
\ge c_3 g_{d,\alpha,\Phi}(x,y) \ge c_4\overline u_{d,\alpha}(x,y)$. Hence, to get the result, it suffices to show that
\begin{align}\label{e:resolvent-claim}
	\int_1^\infty e^{- t}\overline p(t,x,y)dt \le c_5\overline u_{d,\alpha}(x,y)\quad \text{for all $x,y\in K$}.
\end{align}

By the $L^\infty$-contraction property, the map $t\mapsto \sup_{x,y\in \overline D}\overline p(t,x,y)$ is non-increasing. Thus, by \eqref{e:main-HKE1}, we have $\overline p(t,x,y) \le c_6$ for all $t\ge 1$ and $x,y \in \overline D$. It follows that for all $x,y \in K$,
\begin{align*}
		\int_1^\infty e^{- t}\overline p(t,x,y)dt \le c_6 \int_1^\infty e^{- t} dt  = c_6e^{-1}.
\end{align*}
Combining this with \eqref{e:green-lower-inf}, we arrive at  \eqref{e:resolvent-claim}. The proof is complete. 
\end{proof}

We conclude this subsection with  the following proposition,  which provides a criterion for determining whether  $X^0$ and $\overline X$ coincide when started from $D$. Although the proposition is not needed in the proofs of our main theorems, we include it because of its independent interest. In the case $\Phi=1$, the result is well known; see \cite[Theorem 1.1]{BBC03}.

\begin{prop}\label{p:censored}
	(i) If $\alpha\le 1$, then 
	$\partial D$ is polar for $\overline X$. In particular,  $\overline \sF=\sF^0$, and  when started from $D$, $X^0$ and $\overline X$ have the same law.
	
	\noindent (ii) If $\alpha>1$, then $\overline \sF\neq \sF^0$ and  $X^0$ is a proper subprocess of $\overline X$. 
\end{prop}

Before presenting the proof of Proposition \ref{p:censored}, we recall several notions from potential theory that will be used in the proof.

Recall that the \textit{capacity} associated with  $( \sE^0,\overline \sF)$ is defined by
\begin{align*}
	\text{Cap}_{(\sE^0,\overline \sF)}(A):= \inf_{\text{$U$:open}, \, A\subset U} \,\inf \left\{ \sE^0_1(u,u): u \in \overline \sF, \, u \ge 1 \text{ a.e. on $U$} \right\}.
\end{align*}

A Radon measure $\mu$ on $\overline D$ is said to be of \textit{finite energy integral} with respect to $(\sE^0,\overline \sF)$ if there exists $C>0$ such that
\begin{align*}
	\int_{\overline D} |u(x)| \,\mu(dx) \le C \sE^0_1(u,u)^{1/2} \quad   \text{for all $u \in \overline \sF \cap C_c(\overline D)$.} 
\end{align*}
For $\lambda>0$ and  a Radon measure $\mu$ of finite energy integral, define the function $\overline G_\lambda \mu$ on $\overline D$ by $$\overline G_\lambda \mu(x):=\int_{\overline D} \overline G_\lambda(x,y)\mu(dy), \quad x\in \overline D.$$
It is known that $\overline G_\lambda \mu\in \overline \sF$ and that 
\begin{align*}
	\sE^0_\lambda(\overline G_\lambda \mu, u) = \int_{\overline D} u(x) \, \mu(dx) \quad   \text{for all $u \in \overline \sF \cap C_c(\overline D)$.} 
\end{align*}
See \cite[(2.2.2) and Exercise 4.2.2]{FOT}. Moreover, for any compact set $K\subset \overline D$, there exists a unique Radon measure $\mu_K$ of finite energy integral, supported on $K$, such that
\begin{align}\label{e:capacity-equilibrium-measure}
	\begin{split}	\text{Cap}_{(\sE^0,\overline \sF)}(K)&=\mu_K(K)\\
		&=  \sup\left\{ \mu(K):\begin{array}{l}
			\mu \text{ is a Radon measure of finite energy integral}\\
		\text{supported on $K$ such that $\overline G_1\mu \le 1$ q.e. on $\overline D$}
		\end{array} \right\}.
		\end{split}
\end{align}
The measure $\mu_K$ is called the 1-\textit{equilibrium measure} of $K$ associated with  $(\sE^0,\overline \sF)$; see \cite[(2.2.13)]{FOT}.

We now present the proof of Proposition \ref{p:censored}.

\begin{proof}[Proof of Proposition \ref{p:censored}]
	To prove the proposition, since $(\sE^0,\sF^0)$ has no killing part, it suffices to show that $\partial D$ is polar for $\overline X$ if and only if $\alpha\le 1$.  Since $\overline X$ admits a heat kernel, by \cite[Theorems 4.1.2 and 4.2.1]{FOT}, this is equivalent to proving  $	\text{Cap}_{(\sE^0,\overline \sF)}(\partial D)=0$ if and only if $\alpha \le 1$.
	
	 Define a bilinear form  $(\sC,\sF^\sC)$  on $L^2(D)$  by \begin{align*}
		\sC(u,v)&:= \frac12 \int_{D\times D} \frac{(u(x)-u(y))(v(x)-v(y))}{|x-y|^{d+\alpha}}  dxdy,\\
		\sF^\sC&:=\text{the closure of  $\text{Lip}_c(\overline D)$ in $L^2(D)$ under the norm $\sC_1^{1/2}$}.
	\end{align*} 
	Since $\Phi\ge 1$, we have $\overline \sF \subset \sF^\sC$ and $\sC(u,u) \le C_0 \sE^0(u,u)$ for all $u \in \overline \sF$. Thus, for any compact set $K\subset \overline D$, a standard approximation argument (see, e.g., \cite[Lemma 2.2.7(ii)]{FOT}), together with the fact that  $\text{Lip}_c(\overline D)$ is a common core for both $(\sE^0,\overline \sF)$ and $(\sC,\sF^\sC)$,  yields 
	\begin{align}\label{e:capacity-compare}
		\text{Cap}_{(\sC,\sF^\sC)}(K)&=\inf \left\{ \sC_1(u,u): u \in  \text{Lip}_c(\overline D), \, u \ge 1 \text{ on $K$} \right\}\nn\\
		&\le C_0\inf \left\{ \sE^0_1(u,u): u \in  \text{Lip}_c(\overline D), \, u \ge 1 \text{ on $K$} \right\}  = C_0	\text{Cap}_{(\sE^0,\overline \sF)}(K).
	\end{align}Here, $	\text{Cap}_{(\sC,\sF^\sC)}$ is the capacity associated with $(\sC,\sF^\sC)$.

If $\alpha > 1$, then, since \eqref{e:VD} holds and the Hausdorff dimension of $\partial D$ is $d-1>d-\alpha$,  \cite[Theorem 2.5(i) and Remark 2.2(i)]{BBC03} imply that $	\text{Cap}_{(\sC,\sF^\sC)}(\partial D)>0$. Hence, by approximating $\partial D$ by compact sets,  \eqref{e:capacity-compare} yields $	\text{Cap}_{(\sE^0,\overline \sF)}(\partial D)>0$, which proves the assertion in this case.
	
	Suppose  that $\alpha\le 1$. In this case, we get from  \cite[Theorem 2.5(i) and Remark 2.2(i)]{BBC03} that $	\text{Cap}_{(\sC,\sF^\sC)}(\partial D)=0$. 
		 If  $	\text{Cap}_{(\sE^0,\overline \sF)}(\partial D)>0$, then there exists a compact set $K\subset \partial D$ such that  $	\text{Cap}_{(\sE^0,\overline \sF)}(K)>0$. 
		 Let $\mu_K$ be the 1-equilibrium measure of $K$ associated with $(\sE^0,\overline \sF)$. 	 
		 By  \eqref{e:capacity-equilibrium-measure}, $\mu_K$ is  supported on $K$, $\mu_K(K)>0$ and  $\overline G_1\mu_K \le 1$ q.e. on $\overline D$.
		 
		  Since \eqref{e:VD} holds, it follows from  \cite[Theorem 1.1]{CK03}, together with a  standard semigroup argument,  that $(\sC,\sF^\sC)$ admits a jointly continuous heat kernel $p^\sC(t,x,y)$ on $(0,\infty)\times \overline D \times \overline D$ satisfying the following estimates: for every $T\in (0,\infty)$, there exists $C=C(T)>1$ such that
		 \begin{align*}
		C^{-1}\left(t^{-d/\alpha}  \wedge \frac{t}{|x-y|^{d+\alpha}}\right)\le  	p^\sC(t,x,y) \le C\left(t^{-d/\alpha}  \wedge \frac{t}{|x-y|^{d+\alpha}}\right),
		 \end{align*}
		  for all $t\in (0,T]$ and $x,y \in \overline D$. Using these estimates, one can follow the proof of Lemma \ref{l:resolvent-estimate} to show that the 1-resolvent of $(\sC,\sF^\sC)$ admits a density $G^\sC_1(x,y)$ and that there exists a constant $c_1>1$,  depending on $K$, such that 
		   $c_1^{-1} \overline u_{d,\alpha}(x,y) \le G^\sC_1(x,y) \le c_1 \overline u_{d,\alpha}(x,y)$ for all $x,y \in K$. By Lemma \ref{l:resolvent-estimate}, it follows that  $G^\sC_1(x,y)$ and $\overline G_1(x,y)$ are comparable on $K\times K$. Consequently, using the fact that $\overline G_1\mu_K\le 1$ q.e. on $K$, we obtain
		  \begin{align*}
		  &\int_{\overline D \times \overline D} G_1^\sC(x,y)\mu_K(dx)\mu_K(dy) =  \int_{K \times K} G_1^\sC(x,y)\mu_K(dx)\mu_K(dy)\\
		  & \le c_2 	  \int_{K \times K} \overline G_1(x,y)\mu_K(dx)\mu_K(dy) \le c_2 \mu_K(K)<\infty.
		  \end{align*}
		  This implies that $\mu_K$ is of finite energy integral with respect to $(\sC,\sF^\sC)$; see \cite[Exercise 4.2.2]{FOT}. Since $\text{Cap}_{(\sC,\sF^\sC)}(K)\le 	\text{Cap}_{(\sC,\sF^\sC)}(\partial D)=0$, 
		  it follows from \cite[Theorem 2.2.3]{FOT} that $\mu_K(K)=0$, which is a contradiction. Therefore, we obtain  $	\text{Cap}_{(\sE^0,\overline \sF)}(\partial D)=0$ and the proof is complete.
\end{proof}

\subsection{A priori  bounds for the heat kernel and Green function of $X$}\label{ss:a-priori-bounds}
Since $X$ is a subprocess of $\overline X$, \eqref{e:main-HKE1} yields the following upper bound for the transition density of $X$.
\begin{prop}\label{p:upper-heatkernel-2}
	The process $X$ has a transition density  $p(t,x,y)$  on $(0,\infty) \times D \times D$. 
Moreover,	for any $T \in (0, \infty)$, there exists a constant $C=C(T)>0$ such that for all $(t,x,y) \in (0, T] \times  D\times D$,
	\begin{align*}
		p(t,x,y) \le  
		C \left[t^{-d/\alpha} \wedge  \left( \frac{ t }{|x-y|^{d+\alpha}}\Phi\left( \frac{(|x-y| \wedge A_0)^2}{((\delta_D(x) \vee t^{1/\alpha}) \wedge A_0)((\delta_D(y) \vee t^{1/\alpha}) \wedge A_0) }\right) \right)\right].	
	\end{align*} 
\end{prop}
For an open set $U\subset D$, let $p_U(t,x,y)$ denote the transition density of $X^U$.	
%PK%no space
By \cite[Theorem 6.1.1]{FOT}, the semigroup $(P^{U}_t)_{t \ge 0}$ of $X^U$ can be represented by 
\begin{align}\label{e:Feymann-Kac}	P^{U}_t f(x)= \E_x \left[\,\exp\bigg( - \int_0^t \kappa(\overline X^{U}_s) ds \bigg)\, f(\overline X^{U}_t) \right].\end{align} 
%PK%no space
Observe that, by \hyperlink{K}{{\bf (K)}}, there exists $C>0$ such that $\kappa(y) \le C (\delta_D(y)\wedge 1)^{-\alpha}$ for all $y\in D$. Hence, by \eqref{e:Feymann-Kac}, for any open set $U\subset D$ with dist$(U,D^c)>0$,  we have
\begin{align*}
	p_U(t,x,y) \ge \overline p_U(t,x,y) \exp \left( - Ct ( \text{dist}(U,D^c) \wedge 1)^{-\alpha} \right)  \quad \text{for all $t>0$, all $x\in U$ and a.e. $y\in U$}.
\end{align*} 
Combining this with Proposition \ref{p:ndl}, we obtain the following proposition.
\begin{prop}\label{p:ndl-killed}
	For any $R>0$ and $b\in (0,1)$, there exists  $C=C(R,b)>0$  such that for any $x_0 \in  D$, $r\in (0,(b\delta_D(x_0)) \wedge R]$ and 
	$t\in (0, (b r)^\alpha]$, it holds that
	\begin{equation*}
			p_{	B(x_0, r)}(t,z,y)  \ge  C t^{-d/\alpha} 
		\quad  \text{for all $z\in B(x_0, bt^{1/\alpha})$ and a.e. $y\in B(x_0, bt^{1/\alpha})$}.
	\end{equation*}
\end{prop}

The following lemma is a standard consequence of Proposition \ref{p:ndl-killed}; see, e.g.,  \cite[Proposition 3.5(ii)]{CKW20}. We omit the proof.

\begin{lemma}\label{l:exit-time-lower-killed}
For any $R>0$, there exists $C=C(R)>0$ such that $\E_x \tau_{B(x,r)} \ge Cr^\alpha$ for all $x\in D$ and $r\in (0,\delta_D(x)\wedge R]$.
\end{lemma}

Let $Z := (V_s,  X_s)_{s\ge 0}$ be the  time-space process  where $V_s = V_0 -s$. 
The law of the time-space process $s\mapsto 
Z_s$ starting from $(t, x)$ will be denoted by 
$\P_{(t,x)}$. 
For an open subset $U$ of $[0,\infty) \times \R^d$, define $
\tau^Z_{U} = \inf\{s > 0:\, 
Z_s\notin U\}$.

For $x_0 \in \overline D$ and $0\le a<b<\infty$, a Borel  function $u:[0,\infty)\times D\to \R$ is said to be \textit{caloric} in $(a,b] \times B_D(x_0,r)$ 
with respect to $X$ if there is a properly exceptional set $\sN$ such that for every relatively compact open set $U\subset 
(a,b] \times B_{ D}(x_0,r)$, it holds that 
$u(t,z)=\E_{(t,z)}u(
Z_{
	\tau^Z_U})$ for all 
$(t,z)\in U$ with $z\notin \sN$.

Using Proposition \ref{p:ndl-killed} and Lemma \ref{l:phr-1} (in place of \cite[Proposition 12.1]{CK25+}), one can follow the argument in 
\cite[Proposition 12.2]{CK25+} to obtain the following interior parabolic H\"older regularity for $X$. We omit the details.

\begin{theorem}\label{t:phr}
	Let  $b \in (0,1)$. 
	There exist constants 
$\lambda=\lambda(b) \in (0,1)$ and $C=C(b)>0$ 
	such that for all  $x \in  D$,  
	$0<r\le \delta_D(x) \wedge R_0$,  $t_0\ge 0$, and bounded caloric function $u$ in $(t_0, 
	t_0+ r^\alpha]\times B(x, r)$ with respect to $X$, there is a properly exceptional set $\sN_u$ such that 
	\begin{equation*}
		|u(s,y) - u(t,z)| \le C \bigg(\frac{|s-t|^{1/\alpha} + |y-z|}{r}\bigg)^{\lambda } \esssup_{[t_0, t_0 + r^\alpha] \times D}\,|u|,
	\end{equation*}
	for all $s,t \in (t_0 + (1-b^\alpha)r^\alpha, 
	t_0+ r^\alpha]$ and $y,z \in B(x, b r)  \setminus \sN_u$.
\end{theorem}

As a standard consequence of Theorem \ref{t:phr}, the heat kernel $p(t,x,y)$ and the Dirichlet heat kernel $p_U(t,x,y)$ for an open set $U\subset D$ admit versions that are jointly continuous on $(0,\infty) \times D \times D$ and $(0,\infty) \times U \times U$, respectively; see \cite[Section 5]{GHH}. Throughout this paper, we always work with these jointly continuous versions.

Observe that by  \hyperlink{K}{{\bf (K)}}, 
 if $\alpha \le 1$, then  there exist constants $r_0 \in (0,1)$ and $C>0$ such that $\kappa(x) \ge C$  for any bounded connected component $D_0$ of $D$ and  $x\in D_0$ with $\delta_{D_0}(x)<r_0$. This shows that condition {\bf (K2)} in \cite{CKSV23} holds. Hence, since $\Phi\ge 1$, the next lemma follows from 
 the part of the proof of \cite[Lemma 4.23]{CKSV23} after \cite[(4.56)]{CKSV23} and the 
 proof of \cite[Lemma 4.25]{CKSV23}.
\begin{lemma}\label{l:eigenvalue-lower}
		For any $R>0$, there exists  a constant  $\mu_0=\mu_0(R)>0$  such that if $D_0$ is a bounded connected component of $D$ with $\text{\rm diam}(D_0) \le R$, then for all $u \in$ \text{\rm Lip}$_c(D_0)$,
	\begin{align*}
		\sE^\kappa(u,u)\ge \mu_0 \lVert u \rVert_{L^2(D_0)}^2.
	\end{align*}
\end{lemma}

\begin{lemma}\label{l:upper-heatkernel-3}
	For any $x_0\in \overline D$, $R>0$  and $T>0$,
	there exist constants   
	$C=C(R, T)>0$ and $\mu=\mu(R, T)>0$ 
	independent of $x_0$ such that
	\begin{equation*}
		p_{B_D(x_0, R)}(t,x,y) \le Ce^{-\mu t}, \quad\;\; 
		 t\ge T, 
		\; x,y \in B_D(x_0,R).
	\end{equation*}
\end{lemma}
\begin{proof}
	Since $X$ is not conservative, 
	using Propositions \ref{p:upper-heatkernel-2} and \ref{p:ndl},  \eqref{e:VD} and Lemma \ref{l:eigenvalue-lower}, 
	one can follow the proof of
\cite[Propositions 4.22 and  4.26]{CKSV23} 	line by line to get the result. We omit the details.			\end{proof}

For an open set $U\subset D$, let 
$$
G_U(x,y):=\int_0^\infty p_U(t,x,y)dt, \quad x,y \in U,
$$
be the Green function of $X$ on $U$.
	
Recall that the function $\overline u_{d,\alpha}(x,y)$ is defined in \eqref{e:definition-green-function}.

	\begin{prop}\label{p:green-upper-bound} 
	For any  $R>0$, there exists  $C=C(R)>0$ such that 		for all $x_0 \in \overline D$ and $x,y \in B_D(x_0,R)$,
		\begin{align*}
   	   G_{B_D(x_0,R)}(x,y) \le C\overline u_{d,\alpha}(x,y).
		\end{align*}\end{prop}
	\begin{proof}
We decompose the integral as
	\begin{align*}
		G_{B_D(x_0,R)}(x,y)= \int_0^{1} p_{B_D(x_0,R)}(t,x,y) dt +  \int_{1}^\infty p_{B_D(x_0,R)}(t,x,y) dt =:I_1+I_2.
	\end{align*}
	By Proposition \ref{p:upper-heatkernel-2}  and Lemma \ref{l:green-prelim}, we obtain $I_1\le \int_0^{1} p(t,x,y) dt \le c_1\overline  g_{d,\alpha,\Phi} (x,y)\le c_2\overline u_{d,\alpha}(x,y).$ 
	On the other hand, using  Proposition  \ref{l:upper-heatkernel-3} and \eqref{e:green-lower-inf}, we get that $	I_2\le c_3 \int_1^\infty e^{-\mu t}dt = c_3 \mu^{-1} e^{-\mu}  \le c_4\overline u_{d,\alpha}(x,y)$. 
		The proof is complete.
	\end{proof}

%%%%%%%%%%%%%%%%%%%%%%%%%%%%%%%%%%%%%%%%%%%%%%%%%%%%%%%%%%%%%%%%%%%%%%%%%%%%%%%%%%%%%%%%%%%%%%%%%%%%%%%%%%%%%%%%%%%%%%%%%%%%%%%%%%%%%
%%%%%%%%%%%%%%%%%                     Key estimates on cut-off distance function                 %%%%%%%%%%%%%%%%%%%%%%%%%%%%%%%%%%%%
%%%%%%%%%%%%%%%%%%%%%%%%%%%%%%%%%%%%%%%%%%%%%%%%%%%%%%%%%%%%%%%%%%%%%%%%%%%%%%%%%%%%%%%%%%%%%%%%%%%%%%%%%%%%%%%%%%%%%%%%%%%%%%%%%%%%%
\section{Barriers and preliminary boundary estimates}\label{s:barriers}

 In this section, we construct barrier functions using the  analytic results established in Sections \ref{s:est-geometric} and \ref{s:F}. Together with the Dynkin-type formula (see Subsection \ref{ss:dynkin}), these yield upper and lower bounds for  probabilities associated with the boundary behavior of $X$, with explicit boundary decay rates as $x\to \partial D$; see Proposition \ref{p:explicit-1}.

Throughout this section, we assume that $D$ is a  $C^{1,\eta}$ open set, and that  \hyperlink{A}{{\bf (A)}}, \hyperlink{B1}{{\bf (B1)}}, \hyperlink{B2}{{\bf (B2)}} and \hyperlink{K}{{\bf (K)}} hold. Recall that $\F$ is defined by \eqref{e:F-symmetrization}. By Lemma \ref{l:F-basic}(ii), $\F \in \mathfrak P_\Phi$.

\subsection{Key estimate}\label{ss:key-estimate}
Recall that
\begin{align*}
L^0_\sB u(x)=
\text{p.v.}\int_D (u(y)-u(x))J(x,y)\, dy 
=\lim_{\eps\to 0}\int_{D,\, |y-x|>\eps}(u(y)-u(x))\frac{\sB(x,y)}{|x-y|^{d+\alpha}}\, dy,
\end{align*}
 $D(r)=\{x\in D: \delta_D(x)<r\}$,  
$\rho$ denotes the regularized distance function.

For $q \in  [(\alpha-1)_+,\alpha-\overline{\beta}) \cap (0,\infty)$ and $R\in (0,R_0]$,  define 
	\begin{align}\label{e:def-hq}
		h_{q,R}(x)=\rho(x)^q  \ind_{D(R)}(x).
	\end{align} 
 
The next proposition is the main result of this subsection  and one of the key estimates of this work.

\begin{prop}\label{p:barrier}
For any $q\in [(\alpha-1)_+,\alpha-\overline{\beta}) \cap (0,\infty)$, there 
exist constants $C=C(q)>0$ and  
$ \theta_4(q)>0$  
 such that for any  $R\in (0,  R_0 ]$ and  $x \in D(R/8)$,
	\begin{align*}
		\left| 
		L^0_\sB
		h_{q,R}(x)-C(\alpha,q,\F) \sB(x,x)  \delta_D(x)^{-\alpha} h_{q,R}(x)\right| 
		\le C (\delta_D(x)/R)^{\theta_4(q)} \delta_D(x)^{q-\alpha}.
	\end{align*}
\end{prop}
\begin{proof}  
Let 
$\nu=\nu_0$ where $\nu_0$ is the constant in \hyperlink{B2}{{\bf (B2)}}, $R\in (0,R_0]$ and $x \in D(R/8)$. 
 Recall that $Q_x\in \partial D$ is the point such that 
$\delta_D(x)=|x-Q_x|$.  
We will work in the coordinate system CS$_{Q_x}$, so that $\bH^{Q_x}=\bH$. For simplicity, we use the notation 
$$
E_\nu:= E_{\nu}^{Q_x}(R), \quad r:=\delta_D(x) \quad \text{and} \quad k:=\rho(x)/r.
$$  
Note that  $C_1^{-1}\le k\le C_1$ by \eqref{e:regular-distance-1}.

 Since $q<\alpha-\ub$, we can choose $\eta'\in (0,\eta]$ and $\gamma \in (0,1)$  such that 
\begin{equation}
	\label{e:gamma}
q(1+\eta')<\alpha -\ub, \quad 
  (\alpha +(1-q)_+)\gamma\le (q \wedge 1)\eta', 
\quad \text{and} \quad (\alpha+\theta_1)\gamma \le \frac{\theta_2}{2},
\end{equation}
where  $\theta_1$ and $\theta_2$ are the constants in \eqref{e:ass-B5-2}. Define
\begin{align*}
	\sO:=B(x, 5^{-1} R^{-\gamma} r^{1+\gamma}).
\end{align*}
Observe that
$$	\delta_{E_\nu}(x) \ge \delta_{\{(\wt y,y_d):4|\wt y|<y_d<R/2\}}((\wt 0,r)) = r/\sqrt{17}>5^{-1}R^{-\gamma}r^{1+\gamma}.$$ Thus, $\sO\subset  E_{\nu}$.
Since $x_d=r$ and  $\rho(x)=kr =kx_d$, we have
\begin{align*}
&    L^0_\sB
h_{q,R}(x)	=\text{p.v.}\int_{\sO}  \frac{(\rho(y)^q-\rho(x)^q)\sB(x,y)}{|x-y|^{d+\alpha}}dy + \int_{E_\nu\setminus \sO}  \frac{(\rho(y)^q-\rho(x)^q)\sB(x,y)}{|x-y|^{d+\alpha}}dy\\
		&\qquad \qquad\quad \;\;+ \int_{D \setminus E_\nu}  \frac{(h_{q,R}(y)-\rho(x)^q) \sB(x,y)}{|x-y|^{d+\alpha}}dy \\
		&=\text{p.v.}\int_{\sO}  \frac{(\rho(y)^q-\rho(x)^q)\sB(x,y)}{|x-y|^{d+\alpha}}dy  \\
		&\quad + \int_{ E_\nu \setminus \sO }   \frac{(\rho(y)^q-(k
			y_d)^q) \sB(x,x)\F((y-x)/r)}{|x-y|^{d+\alpha}}dy \\ 
		&\quad +  \int_{E_\nu\setminus \sO}   \frac{((k
			y_d)^q-(k
			x_d)^q)\sB(x,x)\F( (y-x)/r)}{|x-y|^{d+\alpha}}dy  \\	
		&\quad +  \int_{ E_\nu \setminus \sO }   \frac{(\rho(y)^q-\rho(x)^q) (\sB(x,y)-\sB(x,x)\F((y-x)/r)) }{|x-y|^{d+\alpha}}dy\\
			&\quad +\int_{D \setminus E_\nu}  \frac{(h_{q,R}(y)-\rho(x)^q) \sB(x,y)}{|x-y|^{d+\alpha}}dy  \\
		&=: I_1+I_2+I_3+I_4+I_5.
\end{align*}
 Estimates of $I_2$ and $I_4$, postponed to Lemmas \ref{l:I4} and   \ref{l:I5} respectively, give that
	\begin{equation}\label{e:I4-I5}
|I_2|+	|I_4|  \le  c_1r^{q-\alpha+b_1} \le c_1(r/R)^{b_1}r^{q-\alpha}
	\end{equation}
for some constants $c_1,b_1>0$ independent of $R$ and $x$. In the rest of the proof, we estimate $I_1$, $I_3$ and $I_5$.

We further split $I_1$   as follows:
\begin{align*}
	I_1&=\text{p.v.} \int_{\sO} \frac{( \nabla h_{q,R}(x) \cdot (y-x))\sB(x,x)}{|x-y|^{d+\alpha}}dy + \int_{\sO}  \frac{(\rho(y)^q-\rho(x)^q-\nabla h_{q,R}(x) \cdot(y-x)) \sB(x,x) }{|x-y|^{d+\alpha}}dy \\
	& \quad	+\int_{\sO}(\rho(y)^q-\rho(x)^q) \frac{\sB(x,y)-\sB(x,x)}{|x-y|^{d+\alpha}}dy \\
	& =:I_{1,1}+I_{1,2}+I_{1,3}.
\end{align*}
By symmetry, $I_{1,1}=0$. Since $r<R$, by the triangle inequality,  
$(4/5)r \le \delta_D(y) \le (6/5)r$ for all $y\in \sO$. Hence, 	using 
\eqref{e:regular-distance-1}, \eqref{e:regular-distance-3} and \eqref{e:regular-distance-4},  we obtain that  for all $y\in \sO$ and  $1\le i,j\le d$,
\begin{align*}
	|\partial_{ij}  h_{q,R} (y)| &\le | q(q-1) \rho(y)^{q-2} \partial_{i} \rho(y) \partial_{j}\rho(y)| + |q\rho(y)^{q-1} \partial_{ij}\rho(y)|\le c_2 r^{q-2} .
\end{align*}
Thus, by Taylor's theorem, 
\begin{align*}
		|\rho(y)^q-\rho(x)^q-\nabla h_{q,R}(x) \cdot(y-x)|\le c_3 r^{q-2}|x-y|^2 \quad \text{for all $y\in \sO$}.
\end{align*}
Using this and  $\sB(x,x)\asymp 1$,  we get that
	\begin{align*}
		|I_{1,2}| &\le c_4 r^{q-2}\int_{\sO} \frac{dy}{|x-y|^{d+\alpha-2}}  = c_5r^{q-2} (5^{-1}r^{1+\gamma}/R^\gamma)^{2-\alpha} =  c_6(r/R)^{(2-\alpha)\gamma} r^{q-\alpha}.
	\end{align*}
	For $I_{1,3}$, we note that, by \eqref{e:regular-distance-1} and \eqref{e:regular-distance-4},   $|\rho(y)^q-\rho(x)^q| \le c_7r^{q-1} |x-y|$ for all $y\in \sO$. Thus, when $\alpha\ge 1$, by   \hyperlink{B1}{{\bf (B1)}},
	\begin{equation}\label{e:I13}
		|I_{1,3}| \le c_8 r^{q-1} \int_\sO \frac{(|x-y|/r)^{\theta_0}}{|x-y|^{d+\alpha-1}} dy  = c_9 r^{q-1-\theta_0} (5^{-1}r^{1+\gamma}/R^\gamma)^{\theta_0-\alpha+1}=  c_{10} (r/R)^{(\theta_0-\alpha+1)\gamma} r^{q-\alpha}.
	\end{equation}
		In the case  $\alpha<1$, observe that $|x-y|^2/(\delta_D(x)\delta_D(y)) \le (r/2)^2/(r\delta_D(y))\le c_{11}$ for all $y\in \sO$. Hence, $\sB(x,y) \le c_{12}$ for all $y\in \sO$. Using this, 
		we obtain \eqref{e:I13}  with $\theta_0=0$. This completes the estimate of $I_1$. 
		
		We next estimate $I_3$. 
		 By Proposition \ref{p:I2}, 
	\begin{align*}
		&\bigg| \int_{ 
		\bH \setminus \sO}   \frac{(
		y_d^q - x_d^q)\F( (y-x)/r)}{|x-y|^{d+\alpha}}dy -C(\alpha, q, \F)r^{q-\alpha} \bigg| 
		\le  c_{13}( r/R)^{(2-\alpha)\gamma} r^{q-\alpha}.
	\end{align*}	
	Thus, since  $E_\nu \subset \bH$, $x_d=r$  and $\rho(x)=kx_d$, using $\sB(x,x) \le C_0$ and $k\le C_1$, we obtain that
	\begin{align}\label{e:I2-1}
			&|I_3 - C(\alpha,q,\F) \sB(x,x) r^{-\alpha}  \rho(x)^q| = |I_3 - k^qC(\alpha,q,\F) \sB(x,x) r^{q-\alpha} |\nn\\
			& \le  k^q \sB(x,x) \bigg| \int_{\bH
		 \setminus \sO}   \frac{(y_d^q - x_d^q)\F( (y-x)/r)}{|x-y|^{d+\alpha}}dy  - C(\alpha,q,\F) r^{q-\alpha} \bigg|  \nn\\
		&\quad + k^q \sB(x,x) \bigg| \int_{\bH \setminus E_\nu}    \frac{(y_d^q - x_d^q)\F( (y-x)/r)}{|x-y|^{d+\alpha}}   dy   \bigg|\nn\\
		& \le c_{14} (r/R)^{(2-\alpha)\gamma} r^{q-\alpha} + c_{15} \int_{\bH \setminus E_\nu}    \frac{(y_d^q+r^q)\F( (y-x)/r)}{|x-y|^{d+\alpha}}  dy .
	\end{align}
Recall that $\F \in \mathfrak P_\Phi$.  	Note that $y_d^q\le (r+|x-y|)^q \le 2^q(r^q+|x-y|^q)$ for all $y\in\bH$.  Hence, using  \eqref{e:profile-2} in the second inequality below, and  
\eqref{e:Enu-difference-result-2} 
with $\theta =\alpha $ and $\theta=\alpha-q$ in the third, we get that
	\begin{align}\label{e:I2-2}
		&\int_{\bH \setminus E_\nu}    \frac{(y_d^q+r^q)\F( (y-x)/r)}{|x-y|^{d+\alpha}}  dy   \le (2^q+1) r^q\int_{\bH \setminus E_\nu}   \frac{\F( (y-x)/r)}{|x-y|^{d+\alpha}}   dy + 2^q\int_{\bH \setminus E_\nu}    \frac{\F( (y-x)/r)}{|x-y|^{d+\alpha-q}}   dy \nn\\
			&\le c_{16} r^q\int_{\bH \setminus E_\nu}   \frac{1}{|x-y|^{d+\alpha}} \Phi \bigg( \frac{|x-y|^2}{ry_d} \bigg)   dy + c_{16} \int_{\bH \setminus E_\nu}    \frac{1}{|x-y|^{d+\alpha-q}}  \Phi \bigg( \frac{|x-y|^2}{ry_d} \bigg)    dy \nn\\
		&\le c_{17}(r/R)^{b_2}r^{q-\alpha},
	\end{align}
	for some constants $c_{17},b_2>0$ independent of $R$ and $x$. Combining \eqref{e:I2-1} with \eqref{e:I2-2},   we arrive at
	\begin{align*}
		|I_3 - C(\alpha,q,\F) \sB(x,x) r^{-\alpha}  \rho(x)^q|
		& \le c_{18} (r/R)^{((2-\alpha)\gamma) \wedge b_2} r^{q-\alpha}.
	\end{align*}
	
	Finally, we estimate $I_5$.
	 	By definition of $h_{q,R}$, 
	 	and  \eqref{e:regular-distance-1}, 
	we see that for all $y\in D$,
	\begin{align*}
	h_{q,R}(y)+\rho(x)^q \le c_{19}\delta_D(y)^q + c_{19}r^q \le c_{19}(r+|x-y|)^q + c_{19}r^q \le c_{20}r^q + c_{21}|x-y|^q.
	\end{align*} 
	Hence, using \eqref{e:J-Phi}, \eqref{e:regular-distance-1}, and  \eqref{e:Enu-difference-result-1}
	with $\theta=\alpha$ and $\theta=\alpha-q$,  we obtain
	\begin{align*}
		|I_5|&\le c_{22}\int_{D \setminus E_\nu}  \frac{(h_{q,R}(y)+\rho(x)^q) \Phi(|x-y|^2/(r\delta_D(y)))}{|x-y|^{d+\alpha}}dy\\
		&\le  c_{23}r^q \int_{D \setminus E_\nu}  \frac{\Phi(|x-y|^2/(r\delta_D(y)))}{|x-y|^{d+\alpha}}dy   +   c_{24} \int_{D \setminus E_\nu}  \frac{\Phi(|x-y|^2/(r\delta_D(y)))}{|x-y|^{d+\alpha-q}}dy \\
		&\le c_{25} (r/R)^{b_3}r^{q-\alpha},
	\end{align*}
	for some $b_3>0$. 	
	Combining the estimates for $I_1$, $I_3$ and $I_5$ with \eqref{e:I4-I5} completes the proof.
\end{proof}

In the following two lemmas, we continue to  use   the notation of the proof of Proposition \ref{p:barrier}. In particular,  $\eta'\in (0,\eta]$ and $\gamma \in (0,1)$ denote the constants satisfying \eqref{e:gamma}.
 
\begin{lemma}\label{l:I4}
There exists  $C>0$ independent of $R$ and $x$  such that
	\begin{align*}
		|I_2| \le C \delta_D(x)^{q-\alpha+\gamma}.
	\end{align*}
\end{lemma} 
\begin{proof}  
Recall from the proof of Proposition \ref{p:barrier} that $r=\delta_D(x)$, $k=\rho(x)/r$,  $E_\nu=E_\nu^{Q_x}(R)$,   $\sO=B(x,5^{-1}R^{-\gamma}r^{1+\gamma})$, and that 
	$$
		I_2=\int_{ E_\nu \setminus \sO }   \frac{(\rho(y)^q-(k
			y_d)^q) \sB(x,x)\F((y-x)/r)}{|x-y|^{d+\alpha}}dy .
	$$
We continue to work in the coordinate system CS$_{Q_x}$. 
For all $y\in E_\nu$, using Lemma \ref{l:geometry}(ii) and the fact that $|x-y|\vee r\le1$ and $\eta'\le \eta$,  we get that 
	if $q\le 1$, then
	\begin{align}\label{e:I4-q<1}
		|\rho(y)^q-(ky_d)^q|  \le 	|\rho(y)-ky_d|^q  \le c_4(|x-y|\vee r)^{q\eta}|x-y|^q \le c_4(|x-y|\vee r)^{q\eta'}|x-y|^q,
	\end{align}
	         and that if $q>1$, then 
	\begin{equation}\label{e:I4-q>1}
			|\rho(y)^q-(ky_d)^q|  \le 	q(\rho(y) + ky_d)^{q-1}|\rho(y)-ky_d| \le c_5  (|x-y| \vee r)^{q-1+\eta'} |x-y|.
	\end{equation}
Further,	 by \eqref{e:profile-2}, we have 
\begin{align}\label{e:I4-F}
	\F((y-x)/r) \le c_6 \Phi(|x-y|^2/(r\delta_D(y))) \quad\text{for all $y \in E_\nu$}.
\end{align}

	Note that $\alpha-q(1+\eta')>\ub$ by \eqref{e:gamma}. If $q\le 1$, then using 
	\eqref{e:I4-q<1}, $\sB(x,x)\asymp 1$ and \eqref{e:I4-F} in the first inequality below,  Lemma \ref{l:tail-F} with $\theta = \alpha-q(1+\eta')$ and $\theta = \alpha-q$ in the third, we obtain
\begin{align*}
	|I_2|& \le c_7 \int_{E_\nu \setminus \sO} \frac{(|x-y| \vee r)^{q\eta'}}{|x-y|^{d+\alpha-q}} \Phi\bigg( \frac{|x-y|^2}{r\delta_D(y)}\bigg) dy \\
	&\le  c_7 \int_{E_\nu \setminus \sO} \frac{1}{|x-y|^{d+\alpha-q(1+\eta')}} \Phi\bigg( \frac{|x-y|^2}{r\delta_D(y)}\bigg) dy  +  c_7 r^{q\eta'} \int_{E_\nu \setminus \sO} \frac{1}{|x-y|^{d+\alpha-q}} \Phi\bigg( \frac{|x-y|^2}{r\delta_D(y)}\bigg) dy\\
	&\le   \frac{c_8\Phi((r/R)^\gamma)}{( r^{1+\gamma}/R^\gamma)^{\alpha-q(1+\eta')}}  +  \frac{c_8r^{q\eta'}\Phi((r/R)^\gamma)}{( r^{1+\gamma}/R^\gamma)^{\alpha-q}} \\
	&\le \frac{c_8\Phi(1)}{r^{(1+\gamma)(\alpha-q(1+\eta'))}}  +  \frac{c_8\Phi(1)}{ r^{(1+\gamma)(\alpha-q)-q \eta'}}\le \frac{2c_8\Phi(1)}{r^{\alpha-q-\gamma}}.
\end{align*}
In the last inequality, we used $(1+\gamma)(\alpha-q(1+\eta'))\le (1+\gamma)(\alpha-q)-q \eta'\le \alpha-q-\gamma$, which follows from \eqref{e:gamma}, together with the fact that $r<1$.

	Now, suppose that $q>1$. Then $\alpha-1>\alpha-q>\ub$ and  $\alpha-q-\eta'>\alpha-q(1+\eta')>\ub$  by \eqref{e:gamma}. Moreover, we see from \eqref{e:gamma} that
	$$	
	(1+\gamma)(\alpha-q-\eta') - (\alpha-q-\gamma) =  \gamma(\alpha+1-q-\eta')  - \eta' < \alpha \gamma -\eta'\le 0
	$$
	and
	$$
	 (1+\gamma)(\alpha-1) - q+1-\eta' - (\alpha-q-\gamma)  = \alpha \gamma - \eta' \le 0.
	$$
	 Hence, using  \eqref{e:I4-q>1}, $\sB(x,x)\asymp 1$, \eqref{e:I4-F}, Lemma \ref{l:tail-F} and $r\le R\le 1$,  we obtain
	\begin{align*}
		|I_2| &\le c_9 \int_{E_\nu \setminus \sO} \frac{(|x-y| \vee r)^{q-1+\eta'}}{|x-y|^{d+\alpha-1}} \Phi\bigg( \frac{|x-y|^2}{r\delta_D(y)}\bigg) dy\\
		&\le c_9 \int_{E_\nu \setminus \sO} \frac{1}{|x-y|^{d+\alpha-q-\eta'}} \Phi\bigg( \frac{|x-y|^2}{r\delta_D(y)}\bigg) dy+ c_9r^{q-1+\eta'} \int_{E_\nu \setminus \sO} \frac{1}{|x-y|^{d+\alpha-1}} \Phi\bigg( \frac{|x-y|^2}{r\delta_D(y)}\bigg) dy\\
		&\le \frac{c_{10}\Phi((r/R)^\gamma)}{( r^{1+\gamma}/R^\gamma)^{\alpha-q-\eta'}}  +  \frac{c_{10}r^{q-1+\eta'}\Phi((r/R)^\gamma)}{( r^{1+\gamma}/R^\gamma)^{\alpha-1}} 
		 \le  \frac{2c_{10}\Phi(1)}{r^{\alpha-q-\gamma}} .
	\end{align*}
	The proof is complete.
\end{proof}

\begin{lemma}\label{l:I5}
There exist constants $C>0$ and $\lambda \in (0,1/2)$ independent of $R$ and $x$  such that
\begin{align*}
	|I_4| \le C \delta_D(x)^{q-\alpha + \lambda \theta_2}.
\end{align*}
\end{lemma} 
\begin{proof} 
	Recall that  $I_4=\int_{E_\nu\setminus \sO} g(y)dy$, where 
	$$
		g(y):=  \frac{(\rho(y)^q-\rho(x)^q) (\sB(x,y)-\sB(x,x)\F((y-x)/r)) }{|x-y|^{d+\alpha}}.
	$$
	We use the coordinate system CS$_{Q_x}$. For all $y\in E_\nu$,
 since $\delta_D(y)\le |x-y|+r \le 2(|x-y| \vee r)$,  we get from Lemma \ref{l:F-basic}(i) that 
	\begin{align}\label{e:I5-main}
	 \big|\sB(x,y)-\sB(x,x) \F((y-x)/r) \big|&\le c_1\left(\frac{|x-y|\vee r \vee \delta_D(y)}{|x-y|\wedge r\wedge \delta_D(y)}\right)^{\theta_1} \, 
	(|x-y|\vee r\vee \delta_D(y))^{\theta_2}\nn\\
	&\le c_2\left(\frac{|x-y|\vee r}{|x-y|\wedge r\wedge \delta_D(y)}\right)^{\theta_1} \, 
	(|x-y|\vee r)^{\theta_2}.
\end{align}

		Define 
	\begin{align*}
		S_1&=\{y\in E_\nu\setminus \sO: 3|x-y|<y_d \vee (2r) \},\\
		S_2&=\{y \in E_\nu\setminus \sO: 3|x-y| \ge 2r \ge y_d\},\\
		S_3&=\{y \in E_\nu \setminus \sO: 3|x-y| \ge y_d > 2r\}.
	\end{align*}
	We first note that $S_1\subset B(x, (2/3)r)$. Indeed,  if $y \in S_1\setminus B(x,(2/3)r)$, then $y_d>3|x-y|\ge |x-y| + (4/3)r$. This leads to a contradiction since $y_d\le  |x-y|+r$. By \eqref{e:regular-distance-1}, it follows that $\delta_D(y)\asymp \rho(y)\asymp r > |x-y|$ for all $y\in S_1$. Hence, using \eqref{e:I5-main} and \eqref{e:gamma}, we obtain
	\begin{align}
		\int_{S_1 }  
	 |g(y)| 
		dy &\le \int_{S_1}\frac{(\rho(y)^q+\rho(x)^q) \big|\sB(x,y)-\sB(x,x)\F((y-x)/r)\big| }{|x-y|^{d+\alpha}}dy\nn\\
		& \le  c_3 r^{q+\theta_1+\theta_2} \int_{S_1} \frac{dy}{|x-y|^{d+\alpha+\theta_1}} \le  c_3 r^{q+\theta_1+\theta_2} \int_{B(x, 5^{-1}R^{-\gamma}r^{1+\gamma})^c} \frac{dy}{|x-y|^{d+\alpha+\theta_1}} \nn\\
		&= c_4 R^{\gamma(\alpha+\theta_1)} r^{q+\theta_1+\theta_2 -(\alpha+\theta_1) (1+\gamma)} \le c_4 r^{q-\alpha + \theta_2-(\alpha+\theta_1)\gamma} \le c_4 r^{q-\alpha + \theta_2/2} .
	\end{align}

	Since $\ub<1 \wedge (\alpha-q)$, we can	choose $\beta\in (\ub, 1\wedge (\alpha-q))$ and  $\lambda \in(0,1/2)$  so that 
	\begin{align}\label{e:lambda}
	\beta(1-\lambda) + \lambda \theta_1<1 \quad \text{and} \quad 	\beta(1-\lambda)+\lambda (\theta_1+\theta_2) <\alpha-q.
	\end{align}
	By  Lemma  \ref{l:distance-in-Enu}(i), $\delta_D(y) \asymp y_d$ for all $y\in E_\nu$. Observe that for all $y \in S_2\cup S_3$, $9|x-y|^2 \ge 2ry_d \ge c_5 r\delta_D(y)$. 	
	Hence,  using \eqref{e:J-Phi}, \eqref{e:without-A0}, $\sB(x,x)\asymp 1$,   \eqref{e:profile-2} and \eqref{e:Phi-scaling}, we obtain
	\begin{align}\label{e:I5-main+}
		\sB(x,y)+\sB(x,x)\F((y-x)/r)) &\le c_6\Phi\bigg( \frac{|x-y|^2}{r\delta_D(y)}\bigg) + c_7\Phi\bigg( \frac{|x-y|^2}{ry_d}\bigg) \nn\\
		& \le c_8	\Phi\left(\frac{9|x-y|^2}{c_5r \delta_D(y)}\right)\le c_9 \left(\frac{|x-y|^2}{r\delta_D(y)}\right)^{\beta}.
	\end{align}
	Using \eqref{e:I5-main} and \eqref{e:I5-main+}, we get that for all $y\in S_2\cup S_3$, 
	\begin{align}\label{e:I5-main2}
&	\left| \sB(x,y)-\sB(x,x)\F((y-x)/r)) \right|\nn\\
	&\le \big(\sB(x,y)+\sB(x,x)\F((y-x)/r))) \big)^{1-\lambda}	\big| \sB(x,y)-\sB(x,x)\F((y-x)/r)) \big|^\lambda \nn\\
		&\le c_{10} \left(\frac{|x-y|^2}{r\delta_D(y)}\right)^{\beta(1-\lambda)} \left(\frac{|x-y|\vee r}{|x-y|\wedge r\wedge \delta_D(y)}\right)^{\lambda\theta_1} \, 
		(|x-y|\vee r)^{\lambda\theta_2}.
	\end{align}
	
	For $y\in S_2$, by \eqref{e:regular-distance-1} and  Lemma  \ref{l:distance-in-Enu}(i), it holds that $\rho(y) \asymp \delta_D(y) \asymp  y_d \le 2r \le 3|x-y|$. Hence, using $\rho(x)\asymp r$ and \eqref{e:I5-main2}, we get that for all $y\in S_2$,
	\begin{align*}
		|g(y)| &\le  \frac{c_{11} r^{q-\beta(1-\lambda)} }{|x-y|^{d+\alpha + \lambda\theta_1-\lambda\theta_2}} \left(\frac{|x-y|^2}{\delta_D(y)}\right)^{\beta(1-\lambda)+\lambda \theta_1}.
	\end{align*}
Since $\beta(1-\lambda)+\lambda\theta_1<(\alpha+\lambda \theta_1-\lambda \theta_2) \wedge 1$ by \eqref{e:lambda}, it follows from Lemma \ref{l:p-estimate0}(i) that
	\begin{align*}
		\int_{S_2} 
		 |g(y)| 
		dy&\le c_{11} r^{q-\beta(1-\lambda)} \int_{D\setminus B(x,(2/3)r)} \frac{1}{|x-y|^{d+\alpha + \lambda\theta_1-\lambda\theta_2}} \left(\frac{|x-y|^2}{\delta_D(y)}\right)^{\beta(1-\lambda)+\lambda \theta_1} dy\nn\\
		&\le c_{12} r^{q-\alpha + \lambda \theta_2 } .
	\end{align*}

	For  $y\in S_3$, we have $\delta_D(y) \ge c_{13}r$ and  $|\rho(y)^q-\rho(x)^q|\le \rho(y)^q+\rho(x)^q \le c_{14}(y_d^q + r^q) \le c_{15}|x-y|^q$. Thus,  by \eqref{e:I5-main2}, we have for all $y\in S_3$,
	\begin{align*}
		|g(y)|\le \frac{c_{16}}{r^{\beta(1-\lambda) + \lambda \theta_1} |x-y|^{d+\alpha-q-\lambda (\theta_1+\theta_2)}} \bigg(\frac{|x-y|^2}{\delta_D(y)}\bigg)^{\beta(1-\lambda)}.
	\end{align*}
	Since $\beta(1-\lambda)<(\alpha-q-\lambda(\theta_1+\theta_2)) \wedge 1$ by \eqref{e:lambda}, we get from Lemma \ref{l:p-estimate0}(i) that
	\begin{align*}
		\int_{S_3} 
		 |g(y)| 
		dy &\le \frac{c_{16}}{r^{\beta(1-\lambda) + \lambda \theta_1} }  \int_{D\setminus B(x,(2/3)r)}\frac{1}{ |x-y|^{d+\alpha-q-\lambda (\theta_1+\theta_2)}} \bigg(\frac{|x-y|^2}{\delta_D(y)}\bigg)^{\beta(1-\lambda)} dy\\
		&\le  c_{13} r^{q-\alpha+\lambda\theta_2}.
	\end{align*}The proof is complete.
\end{proof}

\subsection{Preliminary boundary estimates}\label{ss:prelim-ub}

Since $D$ is a $C^{1,\eta}$ open set,  there exists a constant $a \in (0,1/4]$ such that for any $Q\in \partial D$ and $R\in (0,\wh R]$, one can choose a point  $w_{Q,R} \in B_D(Q,R)$ satisfying 
\begin{align}\label{e:W(Q,R)}
	B(w_{Q,R}, a R) \subset D(R)\setminus  D(R/4) . 
\end{align}
We fix such a point $w_{Q,R}$ and set
\begin{align}\label{e:def-W(Q,R)}
	W(Q,R):=B(w_{Q,R}, a R) . 
\end{align}

Recall that $p\in [(\alpha-1)_+, \alpha-\ub)$ is the unique positive constant satisfying \eqref{e:C(alpha,p,F)}.

The goal of this subsection is to establish the following proposition. 

\begin{prop}\label{p:explicit-1}
	There exist constants $\eps_1 \in (0,1/8]$ and $C>1$ such that for any $Q\in \partial D$,    $R\in (0, R_0]$ and  $x\in B_D(Q,\eps_1R/2)$, the following hold:
	\begin{gather*}
		\P_x\big(X_{\tau_{B_D(Q,\eps_1R)}} \in W(Q,R) \big) \le C(\delta_D(x)/R)^p,\\
		\P_x\big(X_{\tau_{B_D(Q,\eps_1R)}} \in D \big) \ge C^{-1}(\delta_D(x)/R)^p.
	\end{gather*}
\end{prop}

To prove Proposition \ref{p:explicit-1}, we first construct barrier functions using Proposition \ref{p:barrier}.   Throughout this subsection, we let  $\theta_3$ and 
$\theta_4(p)$
 be the  constants in \eqref{e:def-K} and Proposition \ref{p:barrier} respectively. Set
\begin{align*}
	q:=p +  2^{-1} \left(  (\alpha -\overline{\beta} -p) \wedge 
	\theta_3 \wedge 	 \theta_4(p)	\right).
\end{align*}
For $R\in (0,R_0]$, define
\begin{align}\label{e:def-phipsi}
	\phi_{R}(x):= \frac{3h_{p,R}(x)}{R^p} - \frac{h_{q,R}(x)}{(C_1R)^q} \quad \text{and} \quad 
	\psi_{R}(x):= \frac{h_{p,R}(x)}{R^p} + \frac{h_{q,R}(x)}{(C_1R)^q},
\end{align}
where $C_1>1$ is the constant in \eqref{e:regular-distance-1}. 
Since $\rho\in C^2(\R^d \setminus \partial D)$,  we have $\phi_R, \psi_R \in C^2(D(R))$.
From \eqref{e:def-hq} and \eqref{e:def-phipsi}, we have 
	\begin{align}\label{e:(a)}
 (\rho(x)/R)^p \ind_{D(R)}(x)\le \psi_R(x) \le \phi_R(x) \le 3(\rho(x)/R)^p  \ind_{D(R)}(x) \quad \text{for all }x \in D,
	\end{align}

\begin{lemma}\label{l:Dynkin1}
	There  exists a constant $\eps_1 \in (0,1/8]$  such that  for all $R \in (0, R_0]$,		\begin{align}\label{e:(c)}
		L^\kappa_\sB
		\phi_R(x) \le 0 \le 
		L^\kappa_\sB
		\psi_R(x)  \quad 
		 \text{for all } x \in  D(\eps_1R).
	\end{align}
	\end{lemma}
\begin{proof}  
Let $R\in (0,R_0]$.	
	By  Proposition \ref{p:barrier}, \eqref{e:def-K} and  \eqref{e:C(alpha,p,F)},  for all $x \in D(R/8)$,
	\begin{align}\label{e:generator-hp}
		\left|	
		L^\kappa_\sB
		h_{p,R}(x) \right|&=	\left|	
		L^0_\sB
		h_{p,R}(x) - \kappa(x) h_{p,R}(x)  \right| \nn\\
		& \le \left|  
		L^0_\sB
		h_{p, R}(x) - C(\alpha,p, \F) \sB(x,x)\delta_D(x)^{-\alpha} h_{p, R}(x) \right| \nn\\
		&\quad +    \left| \kappa(x) - \kappa_0 \sB(x,x)\delta_D(x)^{-\alpha}  \right| h_{p, R}(x)\nn\\
		&\le c_1 (\delta_D(x)/R)^{ \theta_4(p)} \delta_D(x)^{p-\alpha}  + c_2 \delta_D(x)^{p-\alpha+\theta_3}\nn\\
		&\le  (c_1+c_2)(\delta_D(x)/R)^{\theta_3 \wedge 
	 \theta_4(p)} \delta_D(x)^{p-\alpha}.
	\end{align}
	Moreover, using Proposition \ref{p:barrier},  $\sB\ge C_0^{-1}$, \eqref{e:def-K} and \eqref{e:regular-distance-1},  since $C(\alpha,q,\F)>C(\alpha,p,\F)$ by Lemma \ref{l:constant}, we also get that for all $x \in D(R/8)$,
	\begin{align}\label{e:generator-hq}
		L^\kappa_\sB
		h_{q, R}(x) &=	
		L^0_\sB
		h_{q, R}(x) - \kappa(x) h_{q, R}(x) \nn\\
		&=    
		L^0_\sB
		h_{q, R}(x) - C(\alpha,q, \F)\sB(x,x)\delta_D(x)^{-\alpha} h_{q, R}(x)   \nn\\
		&\quad  +  (C(\alpha,q,\F) - C(\alpha,p,\F))  \sB(x,x) \delta_D(x)^{-\alpha} h_{q, R}(x)  \nn\\
		&\quad -   (\kappa(x) - C(\alpha,p,\F) \sB(x,x)\delta_D(x)^{-\alpha})h_{q, R}(x) \nn\\
		&\ge - c_3(\delta_D(x)/R)^{ \theta_4(q)} \delta_D(x)^{q-\alpha}\nn\\
		&\quad  + 
		 C_0^{-1} 
		(C(\alpha,q,\F) - C(\alpha,p,\F))\delta_D(x)^{-\alpha}  \rho(x)^q   -  c_4 \delta_D(x)^{q-\alpha+\theta_3}\nn\\
		&\ge c_5\delta_D(x)^{q-\alpha}-(c_3+c_4)(\delta_D(x)/R)^{\theta_3 \wedge 
	 \theta_4(q)} \delta_D(x)^{q-\alpha}.
	\end{align}

	Since $p<q<p+\theta_3 \wedge \theta_4(p)$, there exists $\eps_1 \in (0,1/8]$ such that
	\begin{align}\label{e:property-eps1}
	 		c_5 s^{q-p-\theta_3 \wedge 
	 \theta_4(p)} - (c_3+c_4)s^{q-p-\theta_3 \wedge 
	 \theta_4(p)+\theta_3 \wedge 
	 \theta_4(q)} \ge  3 C_1^{q}(c_1+c_2) 
	\quad \text{for all} \;\, s \in (0,\epsilon_1].
	\end{align}
 Using \eqref{e:generator-hp}, \eqref{e:generator-hq} and \eqref{e:property-eps1}, we get that  for all $x \in D(\epsilon_1 R)$,
	\begin{align*}
		&-
		L^\kappa_\sB
		\phi_R(x) =   (C_1R)^{-q}
		L^\kappa_\sB
		h_{q, R}(x) - 3 R^{-p}	
		L^\kappa_\sB
		h_{p, R}(x)  \nn\\
			& \ge
		 C_1^{-q}   R^{-q} [c_5\delta_D(x)^{q-\alpha}-(c_3+c_4)(\delta_D(x)/R)^{\theta_3 \wedge 
	 \theta_4(q)} \delta_D(x)^{q-\alpha}] \\
	 &\qquad - 3 R^{-p}	(c_1+c_2)(\delta_D(x)/R)^{\theta_3 \wedge 
	 \theta_4(p)} \delta_D(x)^{p-\alpha}\\
	 & \ge
		    \Big( C_1^{-q}c_5(\delta_D(x)/R)^{q} - C_1^{-q}(c_3+c_4)(\delta_D(x)/R)^{q+\theta_3 \wedge 
	 \theta_4(q)}\\
	 &\qquad - 3    (c_1+c_2)(\delta_D(x)/R)^{\theta_3 \wedge 
	 \theta_4(p)+p} \Big) \delta_D(x)^{-\alpha}\\
		& \ge   \left(   C_1^{-q}c_5 (\delta_D(x)/R)^{q-p-\theta_3 \wedge 	 \theta_4(p)}   
		- 
		 C_1^{-q}(c_3+c_4) (\delta_D(x)/R)^{q-p-\theta_3 \wedge 
	 \theta_4(p)+\theta_3 \wedge 
	 \theta_4(q)}- 3 (c_1+c_2) \right) \nn\\
		&\qquad \times (\delta_D(x)/R)^{p+\theta_3 \wedge 
	 \theta_4(p)}  \delta_D(x)^{-\alpha} \ge 0,
	\end{align*}
	and
	\begin{align*}
		&
		L^\kappa_\sB
		\psi_R(x) = (C_1R)^{-q} 
		L^\kappa_\sB
		h_{q, R}(x)  + R^{-p}	
		L^\kappa_\sB
		h_{p, R}(x) \nn\\
		& \ge  \left(  C_1^{-q}c_5  (\delta_D(x)/R)^{q-p-\theta_3 \wedge 
	 \theta_4(p) }  - C_1^{-q}(c_3+c_4) (\delta_D(x)/R)^{q-p-\theta_3 \wedge 
	 \theta_4(p)+\theta_3 \wedge 
	 \theta_4(q)}
		- (c_1+c_2) \right)\nn\\
		&\quad \times  (\delta_D(x)/R)^{p+\theta_3 \wedge 
	 \theta_4(p)}  \delta_D(x)^{-\alpha} \ge 0.
	\end{align*}
	
	The proof is complete. 
\end{proof}

\begin{proof}[Proof of Proposition \ref{p:explicit-1}]
	Let $Q\in \partial D$ and $R\in (0,R_0]$.  Let  $\eps_1\in (0,1/8]$ be the constant in Lemma \ref{l:Dynkin1} and $x\in B_D(Q,\eps_1R/2)$.
	Define $\phi_R$ and $\psi_R$ by \eqref{e:def-phipsi}. Then $\phi_R$ and $\psi_R$ are 
	bounded 
	and $\phi_R,\psi_R \in C^2(D(\eps_1R))$. 
 By \eqref{e:regular-distance-1}, 
\eqref{e:W(Q,R)}
and
\eqref{e:(a)}, 
we get that
\begin{align}\label{e:explicit-phi}
\phi_R(x) \le 3(C_1\delta_D(x)/R)^p,  \qquad
 \inf_{y\in W(Q,R)} \phi_R(y) \ge \inf_{y\in W(Q,R)} 
(C_1^{-1}\delta_D(y)/R)^p \ge (4C_1)^{-p},
\end{align}
and 
\begin{align}\label{e:explicit-psi}
	\begin{split} 
		&\psi_R(x) \ge (C_1^{-1}\delta_D(x)/R)^p, \qquad \sup_{y\in D} \psi_R(y) \le 3C_1^p.
	\end{split} 
\end{align}
On the other hand, 
by Lemma \ref{l:Dynkin1}, 
we obtain from Corollary \ref{c:Dynkin-local}  that
\begin{align*}
	&\phi_R(x) = \E_x\left[ \phi_R( X_{\tau_{B_D(Q,\eps_1R)}})\right]  - \E_x \left[ \int_0^{\tau_{B_D(Q,\eps_1R)}} 
	L^\kappa_\sB
	\phi_R(X_s) ds\right]\ge \E_x\big[ \phi_R( X_{\tau_{B_D(Q,\eps_1R)}})\big] ,\\
		&\psi_R(x) = \E_x\left[ \psi_R( X_{\tau_{B_D(Q,\eps_1R)}})\right]  - \E_x \left[ \int_0^{\tau_{B_D(Q,\eps_1R)}}  
		L^\kappa_\sB
		\psi_R(X_s) ds\right]\le \E_x\big[ \psi_R( X_{\tau_{B_D(Q,\eps_1R)}})\big] .
\end{align*}
Combining these with \eqref{e:explicit-phi} and \eqref{e:explicit-psi}, we arrive at
\begin{gather*}
	3(C_1\delta_D(x)/R)^p \ge \phi_R(x) \ge \E_x\big[ \phi_R( X_{\tau_{B_D(Q,\eps_1R)}})\big]  \ge 	
	(4C_1)^{-p}
	\P_x\big(X_{\tau_{B_D(Q,\eps_1R)}} \in W(Q,R) \big),\\
	 (C_1^{-1}\delta_D(x)/R)^p \le \psi_R(x) \le \E_x\big[ \psi_R( X_{\tau_{B_D(Q,\eps_1R)}})\big] \le 
	 3 C_1^p \P_x\big(X_{\tau_{B_D(Q,\eps_1R)}} \in D \big). 
\end{gather*}
These yield the desired results.
\end{proof}

%%%%%%%%%%%%%%%%%%%%%%%%%%%%%%%%%%%%%%%%%%%%%%%%%%%%%%%%%%%%%%%%%%%%%%%%%%%%%%%%%%%%%%%%%%%%%%%%%%%%%%%%%%%%%%%%%%%%%%%%%%%%%%%%%%%%%
%%%%%%%%%%%%%%%%%%%%%%%%                               Barriers                         %%%%%%%%%%%%%%%%%%%%%%%%%%%%%%%%%%%%%%%%%%%%%
%%%%%%%%%%%%%%%%%%%%%%%%%%%%%%%%%%%%%%%%%%%%%%%%%%%%%%%%%%%%%%%%%%%%%%%%%%%%%%%%%%%%%%%%%%%%%%%%%%%%%%%%%%%%%%%%%%%%%%%%%%%%%%%%%%%%%

\section{Box method}\label{s:box}

In this section, we continue to assume that $D$ is a  $C^{1,\eta}$ open set, and that  \hyperlink{A}{{\bf (A)}}, \hyperlink{B1}{{\bf (B1)}}, \hyperlink{B2}{{\bf (B2)}} and \hyperlink{K}{{\bf (K)}} hold. Let $\eps_1\in (0,1/8]$ be the constant in Proposition \ref{p:explicit-1}. Recall that $W(Q,R)$ is defined in \eqref{e:def-W(Q,R)}.

The following proposition is the main result of this section, showing that the survival probability $	\P_x(X_{\tau_{B_D(Q,\eps_1R)}} \in D )$ is dominated by the survival probability  with vertical escape, $	\P_x(X_{\tau_{B_D(Q,\eps_1R)}} \in W(Q,R) )$.
\begin{prop}\label{p:explicit-2}
	There exists $C>1$ such that for any $Q\in \partial D$,    $R\in (0, \wh R]$ and  $x\in B_D(Q,\eps_1R/2)$, 
	\begin{align*}
		\P_x\big(X_{\tau_{B_D(Q,\eps_1R)}} \in D \big)\le C	\P_x\big(X_{\tau_{B_D(Q,\eps_1R)}} \in W(Q,R) \big).
	\end{align*}
\end{prop} 
To prove Proposition \ref{p:explicit-2}, we adapt the \textit{box method},  whose origin dates back to \cite{BB89}. We will use some ideas from a  version of this method, adapted to jump processes, recently developed in \cite{CC24}. The proof of Proposition \ref{p:explicit-2} will be given at the end of this section.

For an open set $U\subset D$,  $x\in U$ and $R>0$, define
\begin{align*}	\Theta^U(x,R)=	\Theta^{U}(\Phi,x,R):= \E_x \left[\int_0^{\tau_{U}} \Phi\left(\frac{R}{\delta_D(X_s)} \right)ds \right].
\end{align*}
Note that there exist constants $C'>C>0$ such that, for 
any $Q\in \partial D$,  $R \in (0,\wh R]$, $z\in B_D(Q,\eps_1R)$ and $w\in W(Q,R)$, $CR\le |z-w| \le C'R$. 
 Hence,	by Lemma \ref{l:levy-consequence} and \eqref{e:Phi-scaling},  for any $Q\in \partial D$,    $R\in (0, \wh R]$ and $x\in B_D(Q,\eps_1R/2)$,
\begin{align}\label{e:box-equiv}
	\P_x\big(X_{\tau_{B_D(Q,\eps_1R)}} \in W(Q,R) \big) \asymp R^{-\alpha}  \Theta^{B_D(Q,\eps_1R)}(x,\eps_1R) \asymp R^{-\alpha}  \Theta^{B_D(Q,\eps_1R)}(x,R).
\end{align}

Before proving Proposition \ref{p:explicit-2}, we first record a consequence of Propositions \ref{p:explicit-2} and  \ref{p:explicit-1},  which provides explicit estimates for survival probabilities.
 \begin{thm} \label{t:explicit-equiv} 
	There are comparison constants such that  for any $Q\in \partial D$,    $R\in (0, \wh R]$ and $x\in B_D(Q,\eps_1R/2)$, 
	\begin{align*}
		&\P_x\big(X_{\tau_{B_D(Q,\eps_1R)}} \in D \big)\asymp \P_x\big(X_{\tau_{B_D(Q,\eps_1R)}} \in W(Q,R) \big) \asymp R^{-\alpha}  \Theta^{B_D(Q,\eps_1R)}(x,R) \asymp  (\delta_D(x)/R)^p.
	\end{align*}
\end{thm}
\begin{proof}
Let $Q\in \partial D$, $R\in (0,\wh R]$ and $x\in B_D(Q,\eps_1R/2)$.  Since $\P_x(X_{\tau_{B_D(Q,\eps_1R)}} \in W(Q,R) )\le \P_x(X_{\tau_{B_D(Q,\eps_1R)}} \in D)$, by  \eqref{e:box-equiv},  it suffices to show that
\begin{align}\label{e:box-claim}
	\begin{split} 
	\P_x\big(X_{\tau_{B_D(Q,\eps_1R)}} \in D \big) &\le c_1(\delta_D(x)/R)^p,\\
	\P_x\big(X_{\tau_{B_D(Q,\eps_1R)}} \in W(Q,R) \big) &\ge c_2(\delta_D(x)/R)^p,
	\end{split}
\end{align}
for some constants $c_1,c_2>0$ independent of $Q,R$ and $x$.

	If $R\le R_0$, then \eqref{e:box-claim} follows from Propositions \ref{p:explicit-1} and \ref{p:explicit-2}. 
	
	Suppose $R_0<R\le \wh R$. Set $R':=(R_0/\wh R)R$, so that $R'\le R\wedge R_0$. Recall that $Q_x\in \partial D$ is the point satisfying $\delta_D(x)=|x-Q_x|$. 	
		If $\delta_D(x) \ge \eps_1 R'/8$, then, since $B(x, \eps_1R'/8) \subset B_D(Q,\eps_1R)$, using \eqref{e:box-equiv}, $\Phi\ge 1$ and Lemma \ref{l:exit-time-lower-killed}, we obtain that
		\begin{align*}
				&\P_x\big(X_{\tau_{B_D(Q,\eps_1R)}} \in W(Q,R) \big) \ge c_3 R^{-\alpha} \E_x [\tau_{B_D(Q,\eps_1R)}]\\
				& \ge c_3 R^{-\alpha} \E_x [\tau_{B(x,\eps_1R'/8)}] \ge c_4 R^{-\alpha} (\eps_1R'/8)^\alpha \ge c_5.
		\end{align*}
		Thus, $\P_x\big(X_{\tau_{B_D(Q,\eps_1R)}} \in W(Q,R) \big)  \asymp \P_x\big(X_{\tau_{B_D(Q,\eps_1R)}} \in D \big)  \asymp (\delta_D(x)/R)^p \asymp 1$ and \eqref{e:box-claim} holds.
		If $\delta_D(x) \le \eps_1 R'/8$, then, since $B_D(Q_x, \eps_1R'/4) \subset B_D(x, \eps_1R'/2) \subset B_D(Q, \eps_1R)$, we get from Propositions \ref{p:explicit-2} and  \ref{p:explicit-1} that 
	\begin{align*}
	&	\P_x\big(X_{\tau_{B_D(Q,\eps_1R)}} \in D \big) \le 
		\P_x\big(X_{\tau_{B_D(Q_x,\eps_1R'/4)}} \in D \big) \\
		&\le 
				 C\P_x\big(X_{\tau_{B_D(Q_x,\eps_1R'/4)}} \in W(Q_x, R'/4) \big)  \le c_6(4\delta_D(x)/R')^p = c_7(\delta_D(x)/R)^p.
	\end{align*}
On the other hand, since $|z-w| \asymp R$ for $z\in  B_D(Q_x,\eps_1R'/4)$ and $w\in W(Q,R) $, by Lemma \ref{l:levy-consequence} and \eqref{e:Phi-scaling}, we have
	\begin{align*}
			\P_x\big(X_{\tau_{B_D(Q_x,\eps_1R'/4)}} \in W(Q,R) \big) \asymp  R^{-\alpha}  \Theta^{B_D(Q_x,\eps_1R'/4)}(x,R) \asymp  R^{-\alpha}  \Theta^{B_D(Q_x,\eps_1R'/4)}(x,R'/4).
	\end{align*}
Therefore, by \eqref{e:box-equiv} and Propositions \ref{p:explicit-2} and  \ref{p:explicit-1}, 
		\begin{align*}
		&\P_x\big(X_{\tau_{B_D(Q,\eps_1R)}} \in W(Q,R) \big) \ge 
		\P_x\big(X_{\tau_{B_D(Q_x,\eps_1R'/4)}} \in W(Q,R) \big) \\
		&\ge c_8R^{-\alpha}  \Theta^{B_D(Q_x,\eps_1R'/4)}(x,R'/4)\ge c_9
		\P_x\big(X_{\tau_{B_D(Q_x,\eps_1R'/4)}} \in W(Q_x, R'/4) \big)  \\
		&\ge c_{10}
		\P_x\big(X_{\tau_{B_D(Q_x,\eps_1R'/4)}} \in D \big)   \ge c_{11}(4\delta_D(x)/R')^p = c_{12}(\delta_D(x)/R)^p,
	\end{align*}
	completing the proof of \eqref{e:box-claim}, and thus of the theorem.
\end{proof}

We begin the proof of Proposition \ref{p:explicit-2} with the following lemma. In the remainder of this section, we fix a constant $\lambda \in (0, 1- \ub/\alpha)$.

\begin{lemma}\label{l:box-EP}	There exists $C\ge 1$ such that for any $z\in \overline D$, 
	$r\in (0,1]$,  $x\in B_D(z,3r/4)$ and $s\in (0,r/8]$,
\begin{align}\label{e:box-EP-1}
	\P_x\big( X_{\tau_{B_D(x,s)}} \in D \big) \le C\bigg(1\wedge \frac{\Theta^{B_D(z,r)} (x,r) }{s^\alpha} \bigg)^{\lambda/(1+\lambda)}.
\end{align}
Consequently, there exists $C'\ge 1$ such that for any $z\in \overline D$, $r\in (0,1]$ and  $x\in B_D(z,3r/4)$,
\begin{align}\label{e:box-EP-2}
	\P_x\big( X_{\tau_{B_D(z,r)}} \in D \big) \le  C'\bigg(1\wedge \frac{\Theta^{B_D(z,r)} (x,r) }{r^\alpha} \bigg)^{\lambda/(1+\lambda)}.
\end{align}
\end{lemma}
\begin{proof}
Let $z\in \overline D$, $r\in (0,1]$, $x\in B_D(z,3r/4)$ and $s\in (0,r/8]$. Since $B_D(x,r/8) \subset B_D(z,r)$, $	\P_x
( X_{\tau_{B_D(z,r)}} \in D 
) \le 	\P_x
( X_{\tau_{B_D(x,r/8)}} \in D 
) $ for all $x\in B_D(z,3r/4)$. Thus,
	\eqref{e:box-EP-2} follows by taking $s=r/8$ in \eqref{e:box-EP-1}. 
	Now we prove \eqref{e:box-EP-1}.
	
	If $\Theta^{B_D(z,r)}(x,r) \ge s^\alpha$, then,  since $	\P_x( X_{\tau_{B_D(x,s)}} \in D )\le 1$,  \eqref{e:box-EP-1}  holds. Suppose  $\Theta^{B_D(z,r)}(x,r) <s^\alpha$. 
	By Lemma \ref{l:levy-consequence}(i), since $B_D(x,s) \subset B_D(z,r)$,  we have
	\begin{align}\label{e:box-EP-term1}
		&\P_x\big( X_{\tau_{B_D(x,s)}} \in D\setminus  B(x,2s) \big) \le  \frac{c_1 \Theta^{B_D(x,s)}(x,s)}{s^\alpha} \le  \frac{c_1 \Theta^{B_D(z,r)}(x,r)}{s^\alpha}.
	\end{align}
	On the other hand, for any $t>0$, it holds that
	\begin{align*}
		\P_x\big( X_{\tau_{B_D(x,s)}} \in  B_D(x,2s) \big) \le \P_x( \tau_{B_D(x,2s)}\ge t ) + \P_x( \tau_{B_D(x,s)} < t \wedge \zeta )=:I_1+I_2.  
	\end{align*}
	Using $B_D(x,2s) \subset B_D(z,r)$,  Markov's inequality and $\Phi\ge 1$, we obtain
	\begin{align*}
		I_1\le \P_x( \tau_{B_D(z,r)}\ge t )  \le t^{-1} \E_x[ \tau_{B_D(z,r)}] \le t^{-1} \Theta^{B_D(z,r)}(x,r).
	\end{align*}
	For $I_2$,  note that on the event $\{\tau_{B_D(x,s)} < t \wedge \zeta\}$, the processes $(X_u)_{u\le \tau_{B_D(x,s)}}$ and $(\overline X_u)_{u\le \tau_{B_D(x,s)}}$ have the same law. Thus, by Lemma \ref{l:EP}, we obtain that for all $t\in (0,1]$,
	\begin{align*}
		I_2 \le \P_x( \overline \tau_{B_{\overline D}(x,s)} < t) \le c_2 (t/s^\alpha)^\lambda.
	\end{align*}
	 Note that $(s^{\alpha \lambda}\Theta^{B_D(z,r)}(x,r) )^{1/(1+\lambda)} <s^\alpha\le (1/8)^\alpha$.  
	Taking $t=(s^{\alpha \lambda}\Theta^{B_D(z,r)}(x,r) )^{1/(1+\lambda)}$, we arrive at
	\begin{align}\label{e:box-EP-term2}
			\P_x\big( X_{\tau_{B_D(x,s)}} \in  B_D(x,2s) \big)  &\le (1+c_2) \bigg(\frac{\Theta^{B_D(z,r)}(x,r) }{s^\alpha}\bigg)^{\lambda/(1+\lambda)}.
	\end{align}
Combining \eqref{e:box-EP-term1} with \eqref{e:box-EP-term2} and using $\P_x( X_{\tau_{B_D(x,s)}} \in  D ) \le1$, we obtain \eqref{e:box-EP-1}.
\end{proof}

\begin{lemma}\label{l:box-cor}
	There exists $C>0$ such that for any $z\in \overline D$, 
	$r\in (0,1]$ and $x\in B_D(z,3r/4)$,
	\begin{align*}
		\frac{\Theta^{B_D(z,r)}(x,r)}{r^\alpha \P_x\big( X_{\tau_{B_D(z,r) }  } \in D \big) }\ge C \bigg(\P_x\big( X_{\tau_{B_D(z,r) }  } \in D \big) +   1\wedge \frac{\Theta^{B_D(z,r)}(x,r)}{r^\alpha}  \bigg)^{1/\lambda}. 
	\end{align*}
\end{lemma}
\begin{proof} Let $z\in \overline D$, $r\in (0,1]$ and $x\in B_D(z,3r/4)$. 
If $\Theta^{B_D(z,r)}(x,r)/r^\alpha \ge 1$, then
\begin{align*}
		\frac{\Theta^{B_D(z,r)}(x,r)}{r^\alpha \P_x\big( X_{\tau_{B_D(z,r) }  } \in D \big) } \ge 	\frac{\Theta^{B_D(z,r)}(x,r)}{r^\alpha }  \ge 1 \ge \frac1{2^{1/\lambda}} \bigg(\P_x\big( X_{\tau_{B_D(z,r) }  } \in D \big) +   1 \bigg)^{1/\lambda}.
	\end{align*} If $\Theta^{B_D(z,r)}(x,r) / r^\alpha<1$, then  
	applying \eqref{e:box-EP-2} twice, we obtain 
	\begin{align*}
	\frac{\Theta^{B_D(z,r)}(x,r)}{r^\alpha \P_x\big( X_{\tau_{B_D(z,r) }  } \in D \big) }&\ge c_1 \bigg(\frac{\Theta^{B_D(z,r)} (x,r) }{r^\alpha} \bigg)^{1/(1+\lambda)} \ge c_2\P_x\big( X_{\tau_{B_D(z,r) }  } \in D \big)^{1/\lambda}.
	\end{align*}
	Moreover, since $\lambda<1$, it holds that
	\begin{align*}
			\frac{\Theta^{B_D(z,r)}(x,r)}{r^\alpha \P_x\big( X_{\tau_{B_D(z,r) }  } \in D \big) } \ge 	\frac{\Theta^{B_D(z,r)}(x,r)}{r^\alpha } \ge  \bigg(	\frac{\Theta^{B_D(z,r)}(x,r)}{r^\alpha  }\bigg)^{1/\lambda}.
	\end{align*}
	Therefore, by the inequality   $a^{1/\lambda} \vee b^{1/\lambda}  \ge 2^{-1/\lambda}(a+b)^{1/\lambda}$ for  $a,b\ge0$,  the desired inequality also holds in this case.
\end{proof}

By Lemma \ref{l:box-cor} and \eqref{e:box-equiv}, the inequality in Proposition \ref{p:explicit-2} holds for any  $x\in B_D(Q,\eps_1R/2)$ for which
\begin{align}\label{e:def-Xi}
\Xi^{B_D(Q,\eps_1R)}(x):=\P_x( X_{\tau_{B_D(Q,\eps_1R) }  } \in D ) +   1\wedge \frac{\Theta^{B_D(Q,\eps_1R)}(x,\eps_1R)}{(\eps_1R)^\alpha}
\end{align}
 satisfies  $
 \Xi^{B_D(Q,\eps_1R)}(x) \ge c$ for some constant $c>0$. 
It therefore remains to consider the case where $\Xi^{B_D(Q,\eps_1R)}(x)$ is small. 
To this end, we use a box method to compare
 $$\Theta^{B_D(Q,\eps_1R)}(x,\eps_1R)/ \P_x( X_{\tau_{B_D(Q,\eps_1R) }  } \in D ) $$ 
 with  the corresponding quantity at a point $y\in B_D(Q,\eps_1R/2)$ for which $
 \Xi^{B_D(Q,\eps_1R)}(y)$ is bounded away from zero.  
The desired conclusion then follows from Lemma \ref{l:box-cor} and \eqref{e:box-equiv}. 
The proof is adapted from   that of  \cite[Proposition 3.1]{CC24}.

\begin{proof}[Proof of Proposition \ref{p:explicit-2}]
Let $Q\in \partial D$, $R\in (0,\wh R]$ and $x\in B_D(Q,\eps_1R/2)$. Set $r:=\eps_1R$. For simplicity, we write $B:=B_D(Q,r)$. By \eqref{e:box-equiv}, it suffices to prove that
	\begin{align}\label{e:explicit-claim}		\P_x(X_{\tau_{B}} \in D ) \le c_1r^{-\alpha} \Theta^{B} (x,r), \quad\;\; x \in B_D(Q,r/2),
	\end{align}
	for some constant $c_1>0$ independent of $Q$ and $R$.
	Define
	\begin{align*}
		U_0:= \bigg\{ y\in B_D\Big(Q, \frac{3r}{4}\Big): 
		\Xi^B(y) \ge 2^{-1}\bigg\},
	\end{align*}
	and for $j\ge 1$,
	\begin{align*}
			U_j&:= \bigg\{ y\in B_D\Big(Q, \frac{3r}{4} - \sum_{k=1}^j \frac{3r}{2\pi^2k^2}\Big):
			\Xi^B(y)  \in [2^{-j-1}, 2^{-j})\bigg\},\\
			V_j&:=\bigcup_{k=0}^{j-1}U_k, \qquad W_j:= B_D\Big(Q, \frac{3r}{4} - \sum_{k=1}^{j-1} \frac{3r}{2\pi^2k^2}\Big) \setminus V_j,
	\end{align*}
where $
	\Xi^{B}(y)$ is defined  in \eqref{e:def-Xi}. 
	Since $\sum_{k=1}^\infty k^{-2} = \pi^2/6$, we have 
	\begin{align}\label{e:box-decomp}
		B_D(Q, r/2) \subset \bigcup_{j=0}^\infty V_j = \bigcup_{j=0}^\infty U_j \subset B_D(Q,3r/4). 
	\end{align}
	Moreover,   for $j\ge 1$,
	\begin{align}\label{e:VW-j}
		\begin{split} 
		V_j \supset \bigg\{ y\in B_D\Big(Q, \frac{3r}{4} - \sum_{k=1}^{j-1} \frac{3r}{2\pi^2k^2}\Big):
		\Xi^B(y)\ge 2^{-j}\bigg\},\\
		W_j = \bigg\{ y\in B_D\Big(Q, \frac{3r}{4} - \sum_{k=1}^{j-1} \frac{3r}{2\pi^2k^2}\Big):
		\Xi^B(y) < 2^{-j}\bigg\}.
		\end{split}
	\end{align}
	
	For $j\ge 1$, define
	\begin{align*}
		a_j = \inf_{y\in V_j} 	\frac{\Theta^{B}(y,r)}{r^\alpha \P_y( X_{\tau_{B }  } \in D ) } \quad \text{if $V_j\neq \emptyset$} \quad \text{ and } \quad a_j=\infty  \quad \text{if $V_j= \emptyset$} .
	\end{align*}
	Since $(a_j)$ is non-increasing and $x\in V_j$ for some $j\ge 1$ by \eqref{e:box-decomp}, in order to prove \eqref{e:explicit-claim}, it suffices to show that there exists $c_2>0$, independent of $Q$ and $R$, such that 
	\begin{align}\label{e:explicit-claim-2}
		\lim_{j\to \infty} a_j \ge c_2.
	\end{align}

	We first claim that, for any $j\ge 1$, the sets  
	$V_j, W_j$ and $W_{i-1} \setminus (W_i\cup V_i)$, $2\le i\le j$, are pairwise disjoint, and that
	\begin{align}\label{e:explicit-decomp}
		 V_j \cup  W_j \cup \Big( \bigcup_{i=2}^j (W_{i-1} \setminus (W_i \cup V_i)) \Big)=	B_D(Q,3r/4).
	\end{align}
	For $j=1$, it is easy to see that the claim holds. For  $j\ge 2$, since $V_i$ is increasing, we have
	\begin{align*}
	B_D(Q,3r/4) &\supset 	V_j \cup W_j \cup \Big( \bigcup_{i=2}^j (W_{i-1} \setminus (W_i \cup V_i)) \Big)  = 	V_j \cup W_j \cup \Big( \bigcup_{i=2}^j (W_{i-1} \setminus W_i) \Big)  \\
	&= V_j \cup W_1 \supset V_1 \cup W_1 = B_D(Q,3r/4),
	\end{align*}
	proving that \eqref{e:explicit-decomp} holds. 
	Clearly, $V_j \cap W_j = \emptyset$. By \eqref{e:VW-j}, we obtain that  for any $2\le i\le j$ and  $y \in  W_{i-1} \setminus (W_i \cup V_i)$, 
	\begin{align} \label{e:explicit-decomp-1}
	y\in  B_D\Big(Q, \frac{3r}{4} - \sum_{k=1}^{i-2} \frac{3r}{2\pi^2k^2}\Big) \setminus  B_D\Big(Q, \frac{3r}{4} - \sum_{k=1}^{i-1} \frac{3r}{2\pi^2k^2}\Big)
	\end{align}
	and 
	\begin{align}\label{e:explicit-decomp-2} 
		\Xi^B(y) < 2^{-i+1}.
		\end{align} 
		By \eqref{e:explicit-decomp-1}, we deduce that the sets $W_j$ and  $W_{i-1} \setminus (W_i \cup V_i)$, $2\le i\le j$, are  pairwise disjoint.
	Moreover,	for any $2\le i\le j$ and  $y \in  W_{i-1} \setminus (W_i \cup V_i)$, 
	 by \eqref{e:explicit-decomp-1}, $y\notin U_k$ for $k\ge i$, and by \eqref{e:explicit-decomp-2}, $y\notin U_k$ for $k\le i-1$. Thus, $ V_j \cap (W_{i-1} \setminus (W_i \cup V_i))=\emptyset$. This completes the proof of the claim.

By the claim,	for any $j\ge 1$ and $ y \in U_j$, we have
	\begin{align*}
	 \P_y( X_{\tau_{B }  } \in D ) & =  \P_y\big( X_{\tau_{B_D(y, 3r/(8\pi^2 j^2))}} \in V_j, \,  X_{\tau_{B }  } \in D \big)\\
	  &\quad  + \sum_{i=2}^j \P_y\big( X_{\tau_{B_D(y, 3r/(8\pi^2 j^2))}} \in W_{i-1} \setminus (W_i \cup V_i), \,  X_{\tau_{B }  } \in D \big) \\
	   &\quad + \P_y\big( X_{\tau_{B_D(y, 3r/(8\pi^2 j^2))}} \in W_j, \,  X_{\tau_{B }  } \in D \big) \\
	  &\quad +  \P_y\big( X_{\tau_{B_D(y, 3r/(8\pi^2 j^2))}} \in D \setminus B_D(Q,3r/4), \,  X_{\tau_{B }  } \in D \big)\\
	  &=:I_1+I_2+I_3+I_4.
	\end{align*}
	For $I_1$, using the strong Markov property twice and the definition of $a_j$, we get that
	\begin{align}\label{e:box-I1}
		I_1 &= \E_y \left[ \P_{X_{\tau_{B_D(y, 3r/(8\pi^2 j^2))}}}( X_{\tau_{B}} \in D )  ; X_{\tau_{B_D(y, 3r/(8\pi^2 j^2))}} \in V_j \right]\nn\\
		&\le \frac{1}{a_j r^\alpha}\E_y \left[  \Theta^{B} (X_{\tau_{B_D(y, 3r/(8\pi^2 j^2))}}, r) ; X_{\tau_{B_D(y, 3r/(8\pi^2 j^2))}} \in V_j \right] \nn\\
		&= \frac{1}{a_j r^\alpha} \E_y \,\bigg[\int_{\tau_{B_D(y, 3r/(8\pi^2 j^2))}}^{\tau_B} \Phi \bigg( \frac{r}{\delta_D(X_s)} \bigg) 
		 ds 
		; X_{\tau_{B_D(y, 3r/(8\pi^2 j^2))}} \in V_j \bigg]\nn\\
		&\le \frac{ \Theta^{B} (y, r) - \Theta^{B_D(y, 3r/(8\pi^2 j^2))}(y,r)}{a_j r^\alpha}.
	\end{align}
	To estimate $I_2$, we observe that for any $2\le i\le j$, $z\in B_D(y, 3r/(8\pi^2 j^2))$ and $w \in W_{i-1} \setminus (W_i \cup V_i)$,
since $y\in U_j$ and $W_i \cup V_i \supset B_D(Q, 3r/4-\sum_{k=1}^{i-1} 3r/(2\pi^2k^2))$, 	
	\begin{align*}
		|z-w|& \ge  |Q-w| - |Q-y| - |y-z| \\
		&\ge \bigg( \frac{3r}{4} - \sum_{k=1}^{i-1} \frac{3r}{2\pi^2k^2}  \bigg) - \bigg( \frac{3r}{4} - \sum_{k=1}^j \frac{3r}{2\pi^2k^2}  \bigg)- \frac{3r}{8\pi^2j^2}\ge  \frac{3r}{2\pi^2i^2}-  \frac{3r}{8\pi^2j^2} \ge  \frac{r}{\pi^2i^2} .
	\end{align*} 
Hence, for any $2\le i\le j$, using the strong Markov property in the equality  below, \eqref{e:explicit-decomp-2} in the first inequality,  and Lemma \ref{l:levy-consequence}(i) in the second inequality,  we obtain
\begin{align*}
	&\P_y\big( X_{\tau_{B_D(y, 3r/(8\pi^2 j^2))}} \in W_{i-1} \setminus (W_i \cup V_i), \,  X_{\tau_{B }  } \in D \big)\\
	&=\E_y \left[ \P_{X_{\tau_{B_D(y, 3r/(8\pi^2 j^2))}}}( X_{\tau_{B}} \in D )  ; X_{\tau_{B_D(y, 3r/(8\pi^2 j^2))}} \in W_{i-1} \setminus (W_i \cup V_i) \right]\\
	&\le 2^{-i+1} \P_y \big(X_{\tau_{B_D(y, 3r/(8\pi^2 j^2))}} \in W_{i-1} \setminus (W_i \cup V_i) \big)  \\
	&\le \frac{c_3 2^{-i}  \Theta^{B_D(y, 3r/(8\pi^2 j^2))} (y, r/(\pi^2i^2))}{(r/(\pi i^2))^\alpha}\le  \frac{c_4 i^{2\alpha} 2^{-i}   \Theta^{B_D(y, 3r/(8\pi^2 j^2))} (y, r)    }{r^\alpha}.
\end{align*}	
It follows that
\begin{align}\label{e:box-I2}
	I_2 \le \frac{c_4  \Theta^{B_D(y, 3r/(8\pi^2 j^2))} (y, r)   }{r^\alpha} \sum_{i=2}^\infty i^{2\alpha} 2^{-i} = \frac{c_5  \Theta^{B_D(y, 3r/(8\pi^2 j^2))} (y, r)   }{r^\alpha}.
\end{align}
For $I_4$, note that for any $z\in B_D(y, 3r/(8\pi^2 j^2))$ and $w\in D\setminus B_D(Q,3r/4)$,
\begin{align*}
	|z-w|& \ge  |Q-w| - |Q-y| - |y-z| \\
&\ge \frac{3r}{4} - \bigg( \frac{3r}{4} - \sum_{k=1}^j \frac{3r}{2\pi^2k^2}  \bigg)- \frac{3r}{8\pi^2j^2}\ge  \frac{3r}{2\pi^2}-  \frac{3r}{8\pi^2j^2} \ge  \frac{r}{\pi^2} .
\end{align*}
Hence, using Lemma \ref{l:levy-consequence}(i), we obtain that
\begin{align}\label{e:box-I4}
	I_4&\le  \P_y( X_{\tau_{B_D(y, 3r/(8\pi^2 j^2))}} \in D \setminus B_D(Q,3r/4) )\nn\\
	& \le \frac{c_6\Theta^{B_D(y, 3r/(8\pi^2 j^2))} (y, r/\pi^2)}{(r/\pi^2)^\alpha}\le \frac{c_6   \Theta^{B_D(y, 3r/(8\pi^2 j^2))} (y, r)   }{r^\alpha}.
\end{align}
We now deal with $I_3$.  
By Lemma \ref{l:box-EP},  since $y\in U_j$, it holds that
\begin{align*}
	&\P_y (X_{\tau_{B_D(y, 3r/(8\pi^2 j^2))}} \in W_j)\le \P_y (X_{\tau_{B_D(y, 3r/(8\pi^2 j^2))}} \in D)\\
	&\le c_7\bigg(1\wedge \frac{\Theta^{B} (y,r) }{(3r/(8\pi^2 j^2))^\alpha} \bigg)^{\lambda/(1+\lambda)}\le c_7\bigg( \frac{ 2^{-j}}{(3/(8\pi^2 j^2))^\alpha} \bigg)^{\lambda/(1+\lambda)} = c_8 j^{2\alpha \lambda/(1+\lambda)} 2^{-j \lambda/(1+\lambda)}.
\end{align*}
Hence, using the strong Markov property,  \eqref{e:VW-j} and the fact that $y\in U_j$, we get that
\begin{align}\label{e:box-I3}
	I_3&= \E_y \left[ \P_{X_{\tau_{B_D(y, 3r/(8\pi^2 j^2))}}}( X_{\tau_{B}} \in D )  ; X_{\tau_{B_D(y, 3r/(8\pi^2 j^2))}} \in W_j \right]\nn\\
	&\le 2^{-j} \P_y \big(X_{\tau_{B_D(y, 3r/(8\pi^2 j^2))}} \in W_j \big) \le c_8 j^{2\alpha\lambda/(1+\lambda)}2^{ - j\lambda/(1+\lambda)+1}  \Xi^B(y).
\end{align}

Combining \eqref{e:box-I1}, \eqref{e:box-I2}, \eqref{e:box-I3} and \eqref{e:box-I4}, we conclude that for all $j\ge 1$ and  $y\in U_j$, 
\begin{align*}
&(1 -  c_8 j^{2\alpha\lambda/(1+\lambda)}2^{ - j\lambda/(1+\lambda)+1}  )  \P_y(X_{\tau_B} \in D) \\
&\le \frac{ \Theta^{B} (y, r) - \Theta^{B_D(y, 3r/(8\pi^2 j^2))}(y,r)}{a_j r^\alpha} + \frac{c_9   \Theta^{B_D(y, 3r/(8\pi^2 j^2))} (y, r)   }{r^\alpha}\\
&\quad  + c_8 j^{2\alpha\lambda/(1+\lambda)}2^{ - j\lambda/(1+\lambda)+1}\frac{\Theta^B(y,r)}{r^\alpha}\\
&\le ( (a_j^{-1} \vee c_9 )  + c_8 j^{2\alpha\lambda/(1+\lambda)}2^{ - j\lambda/(1+\lambda)+1})\frac{\Theta^B(y,r)}{r^\alpha}\\
&\le (a_j^{-1} \vee c_9)( 1 + c_8 j^{2\alpha\lambda/(1+\lambda)}2^{ - j\lambda/(1+\lambda)+1})\frac{\Theta^B(y,r)}{r^\alpha}, 
\end{align*}
where $c_9:=(c_5+c_6) \vee 1$. Since $a_{j+1} = \left(  \inf_{y\in U_j}  \Theta^B(y,r)/(r^\alpha \P_y(X_{\tau_B}\in 
 D 
))\right)\wedge a_j$, it follows that 
\begin{align}\label{e:box-recursive}
	a_{j+1} \ge \bigg( \frac{1 -  c_8 j^{2\alpha\lambda/(1+\lambda)}2^{ - j\lambda/(1+\lambda)+1} }{1 + c_8 j^{2\alpha\lambda/(1+\lambda)}2^{ - j\lambda/(1+\lambda)+1}} \bigg) (a_j \wedge c_9^{-1}) \quad \text{for all $j\ge 1$}.
\end{align}
Let $j_0:=\min\{j\ge 1:   c_8 j^{2\alpha\lambda/(1+\lambda)}2^{ - j\lambda/(1+\lambda)+1} \le 1/2 \}$. By Lemma \ref{l:box-cor} and the definition of $U_k$, if $y\in U_k$ for $1\le j\le j_0$ and $0\le k\le j-1$, then 
\begin{align*}
	\frac{\Theta^B(y,r)}{r^\alpha \P_y(X_{\tau_B} \in D)} \ge c_{10} 2^{-(k+1)/\lambda} \ge c_{10} 2^{-j_0/\lambda}.
\end{align*}
Thus, $ a_{j_0} \ge c_{10} 2^{-j_0/\lambda}$. Combining this with \eqref{e:box-recursive}, we obtain that
\begin{align*}
	\inf_{j>j_0} a_j \ge ((c_{10} 2^{-j_0/\lambda}) \wedge c_9^{-1}) \prod_{j=j_0}^\infty \bigg( \frac{1 -  c_8 j^{2\alpha\lambda/(1+\lambda)}2^{ - j\lambda/(1+\lambda)+1} }{1 + c_8 j^{2\alpha\lambda/(1+\lambda)}2^{ - j\lambda/(1+\lambda)+1}} \bigg).
\end{align*}
Since $\sum_{j=j_0}^\infty c_8j^{2\alpha \lambda/(1+\lambda)}2^{-j\lambda/(1+\lambda)+1}$  converges,  the right-hand side of the above equation is a positive constant. This yields \eqref{e:explicit-claim-2}, which completes the proof of the proposition.
\end{proof}

%%%%%%%%%%%%%%%%%%%%%%%%%%%%%%%%%%%%%%%%%%%%%%%%%%%%%%%%%%%%%%%%%%%%%%%%%%%%%%%%%%%%%%%%%%%%%%%%%%%%%%%%%%%%%%%%%%%%%%%%%%%%%%%%%%%%%
%%%%%%%%%%%%%%%%%%%%%                        Estimate of the Green potential                         %%%%%%%%%%%%%%%%%%%%%%%%%%%%%%%%
%%%%%%%%%%%%%%%%%%%%%%%%%%%%%%%%%%%%%%%%%%%%%%%%%%%%%%%%%%%%%%%%%%%%%%%%%%%%%%%%%%%%%%%%%%%%%%%%%%%%%%%%%%%%%%%%%%%%%%%%%%%%%%%%%%%%%

\section{Estimate of the Green potentials}\label{s:green-pot}

In this section, we continue to assume that $D$ is a  $C^{1,\eta}$ open set, and that  \hyperlink{A}{{\bf (A)}}, \hyperlink{B1}{{\bf (B1)}}, \hyperlink{B2}{{\bf (B2)}} and \hyperlink{K}{{\bf (K)}} hold. Recall that $p\in [(\alpha-1)_+,\alpha-\ub)$ is the positive constant satisfying \eqref{e:C(alpha,p,F)}.

The goal of this section is to obtain the upper bound for the Green potentials  in Proposition \ref{p:bound-for-integral-new}. To this end, we establish suitable upper bounds for the Green functions on balls. 
We begin by establishing survival estimates and small time near-diagonal upper heat kernel estimates. Although these results will be used only to derive upper bounds for the Green functions in the case $d=1\le \alpha$, we 
state and prove them for general dimension $d$ since they may be of independent interest.

\begin{lemma}\label{l:survival}
For any $T>0$,  there exists  $C=C(T) >0$ such that 
\begin{align*}
	\P_x(\zeta>t) \le C \left(1 \wedge \frac{\delta_D(x)}{t^{1/\alpha}}\right)^{p}  \quad \text{for all $t\in (0,T]$ and $x \in D$}.
\end{align*}
\end{lemma}
\begin{proof}
Let $t\in (0,T]$ and $x\in D$. Set  $r_t:=(t/T)^{1/\alpha} \wh R \in (0,\wh R]$. If $\delta_D(x)\ge \eps_1r_t/2 $, where $\eps_1$ is the constant in Proposition \ref{p:explicit-1}, then the desired inequality holds since $\P_x(\zeta>t)\le 1$. If $\delta_D(x)<\eps_1r_t/2$, then  using the fact that $\Phi \ge 1$ in the third inequality below, and Theorem \ref{t:explicit-equiv} in the fourth, we obtain \begin{align*} 	  	
	\P_x (\zeta>t) &\le \P_x (\tau_{B_D(Q_x,\eps_1 r_t)}>t) +\P_x (\tau_{B_D(Q_x,\eps_1 r_t)}<\zeta )\\
	&\le t^{-1}\E_x [\tau_{B_D(Q_x,\eps_1 r_t)}] +\P_x \big(X_{\tau_{B_D(Q_x,\eps_1 r_t)}}\in D \big) \\
	&\le t^{-1} \Theta^{B_D(Q_x,\eps_1 r_t)}(x,r_t)+\P_x \big(X_{\tau_{B_D(Q_x,\eps_1 r_t)}}\in D \big)\\
	&\le c_1 t^{-1} (\delta_D(x)/r_t)^p  r_t^\alpha  +c_3 (\delta_D(x)/r_t)^p= c_4(T)
	(\delta_D(x)/t^{1/\alpha})^p. 
\end{align*}
\end{proof}

\begin{lemma}\label{l:UHKD}
	For any $T>0$,  there exists  $C=C(T) >0$ such that 
	\begin{equation*}
			p(t,x,y) \le C \left(1 \wedge \frac{\delta_D(x)}{t^{1/\alpha}}\right)^{p} \left(1 \wedge \frac{\delta_D(y)}{t^{1/\alpha}}\right)^{p}t^{-d/\alpha} \quad \text{for all 
			$t\in (0,T]$ and $x,y \in D$}.
	\end{equation*}
\end{lemma}
\begin{proof}
	Using  the semigroup property, symmetry,  Proposition \ref{p:upper-heatkernel-2} and Lemma \ref{l:survival},  we obtain that for all $t\in (0,T]$ and $x,y\in D$,
	\begin{align*}	p(t,x,y)& = \int_D \int_D p(t/3,x,z) p(t/3,z,w) p(t/3,w,y) dz dw
		\\
		&\le c_1t^{-d/\alpha} \int_D p(t/3,x,z) dz \int_D p(t/3,y,w)dw  \\
		&= c_1t^{-d/\alpha}\P_x(\zeta>t/3)\P_y(\zeta>t/3)\le c_2\left(1 \wedge \frac{\delta_D(x)}{t^{1/\alpha}}\right)^{p} \left(1 \wedge \frac{\delta_D(y)}{t^{1/\alpha}}\right)^{p}t^{-d/\alpha}  .\end{align*}
\end{proof}

Define for $x,y \in  D$,
\begin{align*}
	u_{d,\alpha}^D(x,y):= \begin{cases}
		|x-y|^{\alpha-d} &\mbox{ if $d>\alpha$},\\[3pt]
		\displaystyle\log\left(e+ \frac{ \delta_D(x) \vee \delta_D(y)}	{|x-y|}\right) &\mbox{ if $d=1=\alpha$},\\[8pt]
		\left( \delta_D(x) \vee \delta_D(y) \vee |x-y|\right)^{\alpha-1} &\mbox{ if $d=1<\alpha$}.
	\end{cases}
\end{align*} 
As in  \eqref{e:green-lower}, we have, in every dimension $d$,  
\begin{align}\label{e:green-lower-2}
	u_{d,\alpha}^D(x,y) \ge  |x-y|^{\alpha-d} \ind_{|x-y|\le 1} + (|x-y|^{\alpha-d} \wedge 1) \ind_{|x-y|>1} \quad \text{for all $x,y \in D$}.
\end{align}
Moreover, since $\sup_{r>0} \log(e+r)/(1+r^{1/2}) <\infty$ and $\sup_{r>0} \ind_{\alpha>1}(1+r^{\alpha-1})/(1+r^{\alpha-1/2})<\infty$, it holds that for all $x,y \in D$,
\begin{align}\label{e:green-upper}
	u_{d,\alpha}^D(x,y)& \le |x-y|^{\alpha-d} \bigg(1 + \ind_{d=1=\alpha} \log\left(e+ \frac{ \delta_D(x) \vee \delta_D(y)}	{|x-y|}\right)+ \ind_{d=1<\alpha} \bigg( \frac{\delta_D(x) \vee \delta_D(y)}{|x-y|}\bigg)^{\alpha-1} \bigg) \nn\\
	&\le C |x-y|^{\alpha-d} \bigg( 1 +  \ind_{d=1\le \alpha}\bigg( \frac{\delta_D(x) \vee \delta_D(y)}{|x-y|}\bigg)^{\alpha-1/2} \bigg) .
\end{align}

\begin{lemma}\label{l:prelub-priori} 
	For any $x_0 \in \overline D$ and $R>0$, there exists $C=C(R)>0$, independent of $x_0$, such that  for all $x,y \in B_D(x_0,R)$,
	\begin{align*}
		G_{B_D(x_0, R)}(x, y) \le  Cu_{d,\alpha}^D(x,y).
	\end{align*}
\end{lemma}
\begin{proof}
	By Proposition \ref{p:green-upper-bound}, it suffices to consider the case $d=1\le \alpha$. 
	Let $x_0 \in \overline D$, $R>0$ and $x,y \in B_D(x_0,R)$. For simplicity,  write $B:=B_D(x_0,R)$.	 By  the symmetry of $G_{B}$, we assume without loss of generality that $\delta_D(x) \le \delta_D(y)$.  

	If $\delta_D(y) > 2R$, then 
	by  Proposition \ref{p:green-upper-bound}, we obtain
	\begin{align*}
		G_{B}(x,y) \le c_1\left( \ind_{\alpha=1} \log\left (e+\frac{1}{|x-y|}\right)  + \ind_{\alpha>1} \right) \le c_2 u^D_{1,\alpha}(x,y) 
	\end{align*}
	for some constant $c_2=c_2(R)>0$. 
	
	Suppose $\delta_D(y)<2R$. 
	 Note that, since $p \ge \alpha-1$,
	\begin{align}\label{e:green-d=1-general}
		u_{1,\alpha}^D(x,y)  \ge (\delta_D(y) \vee |x-y|)^{\alpha-1}   \ge (2R)^{\alpha-1-p}(\delta_D(y) \vee |x-y|)^p.
	\end{align}
	Let $T=(2R)^\alpha$.
	By  the semigroup property,	Lemmas \ref{l:upper-heatkernel-3} and \ref{l:UHKD}, 
	\begin{align}\label{e:green-d=1-1}
&		\int_{2T}^\infty p_{B}(t,x,y) dt \le \int_{2T}^\infty \int_{B} p_{B}(t-T,x,z)  p(T,z,y)  dz dt \nn\\
		&\le c_3m_d(B) T^{-d/\alpha} (\delta_D(y)/T^{1/\alpha})^p \int_{2T}^\infty  
		e^{-\mu (t-T)}    dt = c_4 \delta_D(y)^p \le c_5 u^D_{1,\alpha}(x,y),
	\end{align}
	where we used \eqref{e:green-d=1-general} in the last inequality. 
Further, 	 since $1+2p>1+p \ge \alpha$,  we get from Lemma \ref{l:UHKD}  that 
	\begin{align}\label{e:green-d=1-2}
		& \int_{|x-y|^\alpha}^{2T} p_B(t,x,y)  dt 	\le c_6\int_{|x-y|^\alpha}^{(|x-y| \vee \delta_D(y))^\alpha} \frac{dt}{t^{1/\alpha}}+c_6\delta_D(x)^p\delta_D(y)^p\int_{(|x-y| \vee \delta_D(y))^\alpha}^{2T} \frac{dt}{t^{(1+2p)/\alpha}} \nn\\
		& 
		 \le c_7\left(  \ind_{\alpha=1}    
		  \log\bigg(e+ \frac{\delta_D(y)}{|x-y|}\bigg)+	 \ind _{\alpha>1}(|x-y| \vee \delta_D(y))^{\alpha-1} + \frac{\delta_D(y)^{2p}}{ (|x-y| \vee \delta_D(y))^{1+2p-\alpha}} \right)\nn\\
		 &\le c_7\left(  \ind_{\alpha=1} 
		 \log\bigg(e+ \frac{\delta_D(y)}{|x-y|}\bigg)+	2(|x-y| \vee \delta_D(y))^{\alpha-1} \right) \le c_8 u_{1,\alpha}^D(x,y).
	\end{align}
	Let $\beta \in (\ub,\alpha)$. Using Proposition \ref{p:upper-heatkernel-2}, \eqref{e:without-A0} and \eqref{e:Phi-scaling}, we obtain that
	\begin{align}\label{e:green-d=1-3}
		&	
		\int_0^{|x-y|^\alpha} p_B(t,x,y) dt  \le c_9
		\int_0^{|x-y|^\alpha} 
		\frac{ t \Phi ( {|x-y|^2}/{t^{2/\alpha}})}{|x-y|^{1+\alpha}} dt\nn\\
		&\le \frac{c_{10}}{|x-y|^{1+\alpha-2\beta}}  \int_0^{|x-y|^\alpha}  t^{1-2\beta/\alpha} dt	  =  c_{11} |x-y|^{\alpha-1}\le c_{11}  u_{1,\alpha}^D(x,y).
	\end{align}
	Combining \eqref{e:green-d=1-1}, \eqref{e:green-d=1-2} and \eqref{e:green-d=1-3} yields the desired result.
\end{proof}

\begin{prop}\label{p:prelub} 
	For any $x_0 \in \overline D$ and $R>0$, there exists $C=C(R)>0$, independent of $x_0$, such that  for all $x,y \in B_D(x_0,R)$,
	\begin{align*}
		G_{B_D(x_0,R)}(x, y) \le C
		\left(\frac{\delta_D(x)}{|x-y|} \wedge 1\right)^p	\left(\frac{\delta_D(y) }{|x-y|} \wedge 1\right)^p  u^D_{d,\alpha}(x,y).
	\end{align*}
\end{prop}

 We will use the following  lemma to prove Proposition \ref{p:prelub}.

\begin{lemma}\label{l:Green-iteration}
Let $x_0 \in \overline D$. Suppose that for any $R\ge \wh R$, there exist constants  $A=A(R)>0$ and $b\in [0,\alpha-\ub)$, with $b$ independent of $R$, such that
	\begin{align}\label{e:Green-iteration-ass}
	G_{B_D(x_0,R)}(x, y) \le A
	\left(\frac{\delta_D(x) \wedge \delta_D(y) }{|x-y|} \wedge 1\right)^b
u^D_{d,\alpha}(x,y) \quad \text{for all $x,y \in B_D(x_0, R)$}.
\end{align}
Then  for any $R\ge \wh R$, there exists 
$C=C(R,b)>0$, independent of $x_0$, such that
	\begin{align}\label{e:Green-iteration-result}
		G_{B_D(x_0,R)}(x, y) \le C
		\left(\frac{\delta_D(x) \wedge \delta_D(y) }{|x-y|} \wedge 1\right)^p 	\left(\frac{\delta_D(x) \vee \delta_D(y) }{|x-y|} \wedge 1\right)^b u^D_{d,\alpha}(x,y)
	\end{align}
	for all $x,y \in B_D(x_0, R)$.
\end{lemma}
\begin{proof}	 	
Let  $\eps_1\in (0,1/8]$ be the constant in Proposition \ref{p:explicit-1}.  
We write $B:=B_D(x_0,R)$ and $\wt B  :=B_D(x_0,2R)$. 
By  the symmetry of $G_B(\cdot,\cdot)$,
it suffices to prove \eqref{e:Green-iteration-result} for $x,y \in B$ with $\delta_D(x) \le \delta_D(y)$.   Let $x,y\in B$ with $\delta_D(x) \le \delta_D(y)$. If $\delta_D(x)\ge (2^{-3}\eps_1\wh R/R)|x-y|$, then by \eqref{e:Green-iteration-ass},
\begin{align*}
	G_B(x,y)& \le   (2^{-3}\eps_1\wh R/R)^{-p} A\left(\frac{\delta_D(x) }{|x-y|} \wedge 1\right)^{p+b}
	u_{d,\alpha}^D(x,y)\\
	&\le  (2^{-3}\eps_1\wh R/R)^{-p}A\left(\frac{\delta_D(x)  }{|x-y|} \wedge 1\right)^p 	\left(\frac{\delta_D(y) }{|x-y|} \wedge 1\right)^b 	u_{d,\alpha}^D(x,y),
\end{align*}
proving that \eqref{e:Green-iteration-result} holds.

Suppose $\delta_D(x)<(2^{-3}\eps_1\wh R/R)|x-y|$. Then $\delta_D(y)\le \delta_D(x)+|x-y|<(1+2^{-3}\eps_1\wh R/R)|x-y|$.  Hence, to get \eqref{e:Green-iteration-result}, it suffices to show that
\begin{align}\label{e:Green-iteration-result-2}
	G_B(x,y) \le c_1 \delta_D(x)^p \delta_D(y)^b |x-y|^{-p-b} u^D_{d,\alpha}(x,y),
\end{align}
for some constant $c_1=c_1(R,A,b)>0$ independent of $x_0,x$ and $y$.

 Recall that 
$Q_x\in \partial D$ is the point satisfying  $|x-Q_x|=\delta_D(x)$.  
Set 
$r:=(2^{-2}\wh R/R)|x-y|.$
 Then $\delta_D(x)<\eps_1r/2$, $\delta_D(y)<c_2r$  and  $r<(\wh R/2) \wedge (|x-y|/4)$. Thus, $x\in B_D(Q_x,\eps_1r/2)$ and 
 \begin{align}\label{e:prelub-1}
B_D(Q_x,\eps_1r)\subset 	B_D(Q_x,r) \subset B_D(x, 2r) \subset B_D(x_0, R+2r) \setminus B(y,2r )\subset  \wt B \setminus B(y,2r).
 \end{align}
 In particular, $G_{\wt B}(\cdot, y)$ is regular harmonic in $B_D(Q_x,\eps_1r)$. By the domain  
 monotonicity  and the harmonicity, 
	\begin{align*}
		G_{B}(x, y) \le G_{\wt B}(x, y)&=\E_x\left[G_{\wt B}(X_{\tau_{B_D(Q_x,\eps_1r)}}, y); X_{\tau_{B_D(Q_x,\eps_1r)}} \in 
		B_D(Q_x,r) 	\right] \\
		&\quad 
			+ \E_x\left[G_{\wt B}(X_{\tau_{B_D(Q_x,\eps_1r)}}, y); X_{\tau_{B_D(Q_x,\eps_1r)}}  \in \wt B \setminus B(Q_x,r)	\right]\\
		&=:I_1+I_2.
	\end{align*}

For all $z\in B_D(Q_x,r)$, we have   
 $\delta_D(z) \le \delta_D(x)+2r<c_3r$ and  $2r\le |z-y| \le |x-y|+2r \le c_4r$  by \eqref{e:prelub-1}. Hence, since $\delta_D(y)<c_2r$, it follows from \eqref{e:Green-iteration-ass},  \eqref{e:green-upper} and \eqref{e:green-lower-2} that 
\begin{align*}
	G_{\wt B}(z,y) &\le 
	c_5A(2R) \delta_D(y)^br^{-b} u_{d,\alpha}^D(z,y) \le  
	c_6 \delta_D(y)^br^{\alpha-d-b}\le 
	c_7\delta_D(y)^b r^{-b} u_{d,\alpha}^D(x,y).
\end{align*}
Combining this with Theorem \ref{t:explicit-equiv}, we obtain
\begin{align*}
	&I_1\le  c_7
	 \delta_D(y)^b r^{-b} u_{d,\alpha}^D(x,y)\,\P_x(X_{\tau_{B_D(Q_x,\eps_1r)}}\in B_D(Q_x,r))\le   c_8   \delta_D(x)^p \delta_D(y)^b r^{-p-b} u_{d,\alpha}^D(x,y).
\end{align*}

 For $I_2$, note that by \eqref{e:Levysystem-X-kappa}, 
\begin{align}\label{e:prelub-I2}
	I_2& = \E_x \bigg[ \int_0^{\tau_{B_D(Q_x,\eps_1 r)}}  \int_{\wt B\setminus B(Q_x,r)} G_{\wt B}(z,y)   \frac{\sB(X_s,z)}{|X_s-z|^{d+\alpha}}\,dz  ds\bigg].
\end{align}
We claim that for all $w\in B_D(Q_x,\eps_1r)$,
\begin{align}\label{e:prelub-claim}
	\int_{\wt B\setminus B(Q_x,r)} G_{\wt B}(z,y)   \frac{\sB(w,z)}{|w-z|^{d+\alpha}}\,dz \le \frac{c_{9}
		 \delta_D(y)^b}{ r^{d+b}} \Phi \left( \frac{r}{\delta_D(w)}\right).
\end{align}
Indeed, let $w\in B_D(Q_x,\eps_1r)$.  Note that $B(Q_x,r) \subset B(y, |x-y| + \delta_D(x) + r) \subset B(y, c_{10}r)$ for some constant $c_{10}\ge c_2$. Fix a constant $\beta \in (\ub,1 \wedge (\alpha-b))$. Using \eqref{e:Green-iteration-ass} and \eqref{e:green-upper}, we get that for all $m \in \Z$ and $z\in \wt B \cap (B(y, 2^{m+1}c_{10}r) \setminus B(y,2^{m}c_{10}r))$, 
\begin{align}\label{e:prelub-2}
G_{\wt B}(z,y) 
&\le \frac{c_{11}A(2R)}{2^{m(d-\alpha)}r^{d-\alpha}}\left( \frac{ \delta_D(y)}{2^{m}c_{10}r} \wedge 1 \right)^b   \bigg( 1 + \ind_{1=d\le \alpha} \left(\frac{ \delta_D(y) + 2^{m+1} c_{10}r}	{2^{m}c_{10}r}\right)^{\alpha-1/2}  \bigg)\nn\\
&\le \frac{c_{12}}{2^{m(d-\alpha)}r^{d-\alpha}}\left( \frac{ \delta_D(y)}{2^{m}c_{10}r} \wedge 1 \right)^b   \bigg( 1 + \ind_{1=d\le \alpha} \left(\frac{ \delta_D(y)}	{2^{m}c_{10}r}\right)^{\alpha-1/2}  \bigg).
\end{align}
Using \eqref{e:prelub-2},   \eqref{e:jump-tail-estimate} and \eqref{e:Phi-scaling}, since $\delta_D(y)<c_2r$,   $|w-y|<c_{10}r$ and $d+b-\beta\ge 1+b-\beta>0$, we obtain
\begin{align}\label{e:prelub-3}
&	\int_{\wt B\setminus B(y,2c_{10}r)} G_{\wt B}(z,y)   \frac{\sB(w,z)}{|w-z|^{d+\alpha}}\,dz\nn\\
&  \le  \frac{c_{13}
	 \delta_D(y)^b}{r^{d-\alpha+b}} \sum_{m=1}^\infty2^{-m(d-\alpha+b)}	\int_{\wt B \cap (B(y, 2^{m+1}c_{10}r) \setminus B(y,2^{m}c_{10}r))}   \frac{\sB(w,z)}{|w-z|^{d+\alpha}}\,dz \nn\\
& \le  \frac{c_{13} \delta_D(y)^b}{r^{d-\alpha+b}} \sum_{m=1}^\infty2^{-m(d-\alpha+b)}	\int_{D \setminus B(w,2^{m-1}c_{10}r))}   \frac{\sB(w,z)}{|w-z|^{d+\alpha}}\,dz \nn\\
& \le  \frac{c_{14} \delta_D(y)^b}{r^{d+b}} \sum_{m=1}^\infty2^{-m(d+b)} \Phi \left( \frac{2^{m-1}c_{10}r}{\delta_D(w)}\right)\nn\\
&\le  \frac{c_{15} \delta_D(y)^b}{ r^{d+b}} \Phi \left( \frac{r}{\delta_D(w)}\right) \sum_{m=1}^\infty2^{-m(d+b-\beta)} =  \frac{c_{16} \delta_D(y)^b}{ r^{d+b}} \Phi \left( \frac{r}{\delta_D(w)}\right).
\end{align}

For all  $z\in (\wt B \setminus B(Q_x,r))\cap B(y,2c_{10}r)$, since $|w-z| \ge |Q_x-z| - |Q_x-w| \ge (1-\eps_1)r$,  $|w-z| \le |w-Q_x|+|Q_x-x|+|x-y|+|y-z| \le c_{17}r$ and $\delta_D(z) \le \delta_D(y)+2c_{10}r<3c_{10}r$, using \eqref{e:J-Phi}, \eqref{e:without-A0} and \eqref{e:Phi-scaling},  we get that
\begin{align}\label{e:prelub-kernel-estimate}
	\frac{\sB(w,z)}{|w-z|^{d+\alpha}} \le \frac{c_{18}}{r^{d+\alpha}}\Phi \bigg( \frac{(c_{17}r)^2}{\delta_D(w)\delta_D(z)} \bigg) \le \frac{c_{19}}{\delta_D(z)^\beta r^{d+\alpha-\beta} }\Phi \bigg( \frac{r}{\delta_D(w)} \bigg) . 
\end{align}
Let $m_0\in \Z$ be such that $2^{m_0-1}c_{10}r\le \delta_D(y)<2^{m_0}c_{10}r$. Since $\delta_D(y)<c_2r\le c_{10}r$, we have $m_0\le 0$. For all $m\le m_0-3$ and  $z\in (\wt B \setminus B(Q_x,r))\cap (B(y, 2^{m+1}c_{10}r) \setminus B(y,2^{m}c_{10}r))$, since  $\delta_D(z) \ge \delta_D(y) - 2^{m_0-2}c_{10}r \ge \delta_D(y)/2$,    using \eqref{e:prelub-2} and \eqref{e:prelub-kernel-estimate}, we obtain that 
\begin{align*}
	G_{\wt B}(z,y)   \frac{\sB(w,z)}{|w-z|^{d+\alpha}}\le   \frac{c_{20}}{(\delta_D(y)/2)^\beta r^{2d-\beta}} \Phi \bigg( \frac{r}{\delta_D(w)} \bigg)   \bigg( \frac{1}{ 2^{m(d-\alpha)}}  + \ind_{1=d\le \alpha}\frac{\delta_D(y)^{\alpha-1/2}}{2^{m/2} r^{\alpha-1/2}} \bigg).
\end{align*}
Thus,
\begin{align}\label{e:prelub-4}
		&\int_{(\wt B \setminus B(Q_x,r))\cap B(y,2^{m_0-2}c_{10}r)} G_{\wt B}(z,y)   \frac{\sB(w,z)}{|w-z|^{d+\alpha}}\,dz \nn\\
		& =  \sum_{m=-\infty}^{m_0-3} \int_{(\wt B \setminus B(Q_x,r))\cap (B(y, 2^{m+1}c_{10}r) \setminus B(y,2^{m}c_{10}r))} G_{\wt B}(z,y)   \frac{\sB(w,z)}{|w-z|^{d+\alpha}}\,dz  \nn\\
		&\le  \frac{c_{20}}{(\delta_D(y)/2)^\beta r^{2d-\beta}} \Phi \bigg( \frac{r}{\delta_D(w)} \bigg)    \sum_{m=-\infty}^{m_0-3}   \left( \frac{1}{ 2^{m(d-\alpha)}}  + \ind_{d=1\le \alpha}\frac{\delta_D(y)^{\alpha-1/2}}{2^{m/2} r^{\alpha-1/2}} \right) \int_{B(y, 2^{m+1}c_{10}r) }  \,dz\nn\\
			& \le  \frac{c_{21}}{\delta_D(y)^\beta r^{d-\beta}}\Phi \bigg( \frac{r}{\delta_D(w)} \bigg)  \sum_{m=-\infty}^{m_0-3}  \bigg( 2^{m\alpha}  + \ind_{d=1\le \alpha}\frac{2^{m/2}  \delta_D(y)^{\alpha-1/2}}{r^{\alpha-1/2}} \bigg) \nn\\
			&\le \frac{  c_{22} }{\delta_D(y)^\beta r^{d-\beta}}\Phi \bigg( \frac{r}{\delta_D(w)} \bigg)  \bigg( 2^{m_0\alpha}  + \ind_{d=1\le \alpha}\frac{2^{m_0/2}  \delta_D(y)^{\alpha-1/2}}{r^{\alpha-1/2}} \bigg) \nn\\
			&\le \frac{  c_{23} (\delta_D(y)/r)^\alpha  }{\delta_D(y)^\beta r^{d-\beta}}\Phi \bigg( \frac{r}{\delta_D(w)} \bigg)  \le \frac{  c_{24} \delta_D(y)^b  }{ r^{d+b}}\Phi \bigg( \frac{r}{\delta_D(w)} \bigg) ,
\end{align}
where we used the facts that $2^{m_0} \le 2\delta_D(y)/(c_{10}r)$ in the fourth inequality, and that $b<\alpha-\beta$ and $\delta_D(y)<c_{10}r$ in the last inequality.   
On the other hand, for all $ m_0-2\le m\le 0$ and  $z\in (\wt B \setminus B(Q_x,r))\cap (B(y, 2^{m+1}c_{10}r) \setminus B(y,2^{m}c_{10}r))$,  since $\delta_D(y)<2^{m_0}c_{10}r$, by \eqref{e:prelub-2} and \eqref{e:prelub-kernel-estimate}, it holds  that 
\begin{align*}
	G_{\wt B}(z,y)   \frac{\sB(w,z)}{|w-z|^{d+\alpha}}\le   \frac{c_{25} \delta_D(y)^b}{2^{m(d-\alpha+b)}\delta_D(z)^\beta r^{2d-\beta+b}} \Phi \bigg( \frac{r}{\delta_D(w)} \bigg).  
\end{align*}
Hence, using Lemma \ref{l:Aikawa-integrals}(i) and $b<\alpha-\beta$, we get that
\begin{align}\label{e:prelub-5}
	&	\int_{(\wt B \setminus B(Q_x,r))\cap (B(y, 2c_{10}r) \setminus B(y,2^{m_0-2}c_{10}r))} G_{\wt B}(z,y)   \frac{\sB(w,z)}{|w-z|^{d+\alpha}}\,dz\nn\\
		& =  \sum_{m=m_0-2}^{0} \int_{(\wt B \setminus B(Q_x,r))\cap (B(y, 2^{m+1}c_{10}r) \setminus B(y,2^{m}c_{10}r))} G_{\wt B}(z,y)   \frac{\sB(w,z)}{|w-z|^{d+\alpha}}\,dz  \nn\\
	&  \le \frac{c_{25} \delta_D(y)^b}{ r^{2d-\beta+b}} \Phi \bigg( \frac{r}{\delta_D(w)} \bigg) \sum_{m=m_0-2}^{0} 2^{-m(d-\alpha+b)}	\int_{B_D(y, 2^{m+1}c_{10}r)}   \frac{dz}{\delta_D(z)^{\beta}} \nn\\
	&  \le   \frac{c_{26} \delta_D(y)^b}{ r^{d+b}} \Phi \bigg( \frac{r}{\delta_D(w)} \bigg) \sum_{m=m_0-2}^{0} 2^{m(\alpha-b-\beta)} \le \frac{c_{27} \delta_D(y)^b}{ r^{d+b}} \Phi \bigg( \frac{r}{\delta_D(w)} \bigg).
\end{align}
	Combining \eqref{e:prelub-3}, \eqref{e:prelub-4} and \eqref{e:prelub-5} completes the proof of the claim  \eqref{e:prelub-claim}.

	Now, combining \eqref{e:prelub-I2} with \eqref{e:prelub-claim}, we arrive at
	\begin{align*}
		I_2&\le\frac{c_{9} \delta_D(y)^b}{r^{d+b}} \E_x \bigg[ \int_0^{\tau_{B_D(Q_x,\eps_1 r)}}  \Phi\bigg( \frac{r}{\delta_D(X_s)} \bigg) ds\bigg]  \\
		&=\frac{c_{9}\delta_D(y)^b}{r^{d+b}}  \Theta^{B_D(Q_x,\eps_1r)}(x,r)\le \frac{c_{28}\delta_D(x)^p\delta_D(y)^b}{r^{d+p+b-\alpha}} \le \frac{c_{29}\delta_D(x)^p\delta_D(y)^b}{r^{p+b}}u^D_{d,\alpha}(x,y),
	\end{align*}
	where the second inequality follows from Theorem \ref{t:explicit-equiv} and the third inequality follows from \eqref{e:green-lower-2}.  This completes the proof of \eqref{e:Green-iteration-result-2} and hence the proposition.
\end{proof}

\begin{proof}[Proof of Proposition \ref{p:prelub}]
	By the domain monotonicity of the Green function, it suffices to prove the proposition for $R\ge \wh R$. By 
	Lemma \ref{l:prelub-priori},   condition \eqref{e:Green-iteration-ass} holds for $b=0$. Thus, by Lemma \ref{l:Green-iteration},  condition \eqref{e:Green-iteration-ass} holds for $b=p$. Applying Lemma \ref{l:Green-iteration} once again yields the desired result.
	
	 \end{proof}

Note that, since $D$ is a $C^{1,\eta}$ open set, there exists $C=C(d,\eta, \Lambda)>0$   
such that
\begin{align}\label{e:boundary-layer}
	m_d\left( B_D(Q,r) \cap D(s)\right) \le C s r^{d-1} \quad \text{for all $Q\in \partial D$ and $0<s\le r\le 4\wh R$}.
\end{align}

\begin{prop}\label{p:bound-for-integral-new}  
For  any $\gamma>-p-1$, there exists 
$C=C(\gamma)>0$ 
such that for any $R\in (0, \wh R/2]$, $Q\in \partial D$ and $x\in B_D(Q,R)$,
	\begin{align} \label{e:bound-for-integral-new} 
		\E_x \int_0^{\tau_{B_D(Q,R)}}\delta_D(X_t)^{\gamma }\, dt \le C
		\begin{cases} 
			R^{\alpha+  \gamma   -p} \delta_D(x)^p, &  \gamma >p-\alpha,\\[2pt]
			\delta_D(x)^p
			\log(e+R/\delta_D(x)), &  \gamma =p-\alpha,  
			\\[2pt]
			\delta_D(x)^{\alpha+ \gamma }, &\gamma <p-\alpha.  
		\end{cases}
		\end{align}
\end{prop}
\begin{proof} 
Let  $g_{R,\gamma}(x)$ be the function after the inequality in \eqref{e:bound-for-integral-new}.  
By considering each case separately, we see that
\begin{align}\label{e:integral-g-gamma}
	g_{R,\gamma}(x) \ge \delta_D(x)^{\alpha+\gamma} , \quad x\in B_D(Q,R).
\end{align} 
Let $x\in B_D(Q,R)$ and let $n_0\in \N$ be such that $2^{n_0-1}\delta_D(x)<R\le 2^{n_0}\delta_D(x)$.	 
By the domain monotonicity and  Proposition \ref{p:prelub},
	\begin{align*}
&\int_{B_D(Q,R)} G_{B_D(Q,R)}(x,y) \delta_D(y)^{ \gamma }\, dy\le \int_{B_D(x,2R)} 
G_{B_D(Q,\wh R)}(x,y) \delta_D(y)^{ \gamma }\, dy \\
		&\le c_1 \int_{ B_D(x,2R)}\left(\frac{\delta_D(x) }{|x-y|}\wedge 1\right)^p   \left(\frac{\delta_D(y)}{|x-y|}\wedge 1\right)^p  
	\delta_D(y)^\gamma u_{d,\alpha}^D(x,y) dy\\
		&\le c_1 \bigg(  \int_{ B(x,2^{-n_0-1}R)}
		\delta_D(y)^\gamma u_{d,\alpha}^D(x,y)   dy  \\
		&\qquad \quad   +  \sum_{n=1}^{n_0+2} \frac{\delta_D(x)^p}{(2^{1-n}R)^{2p}}\int_{  B_D(x,2^{2-n}R)\setminus  B(x, 2^{1-n}R)} \delta_D(y)^{p+\gamma} u_{d,\alpha}^D(x,y)   dy  \bigg)\\
		&=:c_1\bigg(I_1+ \sum_{n=1}^{n_0+2}I_{2,n}\bigg).
	\end{align*}
	For all $y \in B(x,\delta_D(x)/2)$, since $(1/2)\delta_D(x) \le \delta_D(y) \le (3/2)\delta_D(x)$, we get from \eqref{e:green-upper} that
	\begin{align*}
		\delta_D(y)^\gamma u_{d,\alpha}^D(x,y) &\le c_2 \delta_D(x)^\gamma |x-y|^{\alpha-d} \bigg( 1+ \ind_{d=1\le\alpha} \left(  \frac{\delta_D(x)}{|x-y|}\right)^{\alpha-1/2}\bigg).
	\end{align*}
	Since $2^{-n_0-1}R\le \delta_D(x)/2$, it follows that
	\begin{align*}
		I_1 &\le c_2\delta_D(x)^\gamma \int_{ B(x,\delta_D(x)/2)} \left( |x-y|^{\alpha-d} +  \ind_{d=1\le \alpha} \delta_D(x)^{\alpha-1/2} |x-y|^{-1/2}\right)  dy     \\
		&  \le c_3 \delta_D(x)^{\alpha+\gamma} \le c_3g_{R,\gamma}(x),
	\end{align*}
	where the last inequality follows from \eqref{e:integral-g-gamma}. Thus, it suffices to estimate $\sum_{n=1}^{n_0+2}I_{2,n}$.

 For all $1\le n \le n_0+2$ and  $y \in    B_D(x,2^{2-n}R)\setminus  B(x, 2^{1-n}R)$, since $\delta_D(x)<2^{1-n_0}R \le 2^{3-n}R\le 4|x-y|$ and $\delta_D(y) \le \delta_D(x)+|x-y|<5|x-y|$,  by
 \eqref{e:green-upper},
\begin{align}\label{e:bound-for-integral-new-1}
	  u^D_{d,\alpha}(x,y)& \le c_4  |x-y|^{\alpha-d} .
\end{align}
	Recall that $Q_x\in \partial D$ is the point satisfying $\delta_D(x)=|x-Q_x|$. 
	For any $1\le n\le n_0+2$, since $\delta_D(x) < 2^{1-n_0} R \le 2^{3-n}R$, we have $B_D(x,2^{2-n}R) \subset B_D(Q_x,2^{4-n}R) = B_D(Q_x,2^{4-n}R) \cap D(2^{4-n}R)$. 
	 Thus, 
 using \eqref{e:bound-for-integral-new-1},  \eqref{e:boundary-layer} and 
 $p+\gamma+1>0$,  
 we obtain that
	\begin{align*}
		I_{2,n} &\le \frac{c_5\delta_D(x)^p}{(2^{-n}R)^{d-\alpha+2p}}\int_{  B_D(x,2^{2-n}R)\setminus  B(x, 2^{1-n}R)} \delta_D(y)^{p+\gamma}   dy \\
		&\le  \frac{c_5\delta_D(x)^p}{(2^{-n}R)^{d-\alpha+2p}} \int_{   B_D(Q_x,2^{4-n}R) \cap D(2^{4-n}R)} \delta_D(y)^{p+\gamma}   dy  
		 \\
		&\le    \frac{c_5\delta_D(x)^p}{(2^{1-n}R)^{d-\alpha+2p}} \sum_{m=n}^{\infty} (2^{p+\gamma} \vee 1) (2^{3-m}R)^{p+\gamma} \int_{  B_D(Q_x,2^{4-n}R)  \cap (D(2^{4-m}R) \setminus D(2^{3-m}R))}  dy\\
			&\le   c_6\delta_D(x)^p (2^{-n}R)^{\alpha-2p-1} \sum_{m=n}^{\infty} (2^{-m}R)^{p+\gamma+1} =  c_7 \delta_D(x)^p (2^{-n}R)^{\alpha+\gamma-p}.
	\end{align*}
	It follows that
	\begin{align*}
		\sum_{n=1}^{n_0+2} I_{2,n} & \le c_7\delta_D(x)^p 	R^{\alpha +\gamma -p}   \sum_{n=1}^{n_0+2} 2^{-n(\alpha+\gamma-p)}  \\
		&\le c_8 \delta_D(x)^p R^{\alpha+\gamma-p} \left( \1_{\gamma>p-\alpha} + n_0 \1_{\gamma=p-\alpha} + 2^{n_0(p-\alpha-\gamma)} \1_{\gamma<p-\alpha} \right).
	\end{align*}
Using $2^{n_0-1}\delta_D(x)<R\le 2^{n_0}\delta_D(x)$, we obtain that the last term  above is comparable to $g_{R,\gamma}(x)$. The proof is complete.
\end{proof}

%%%%%%%%%%%%%%%%%%%%%%%%%%%%%%%%%%%%%%%%%%%%%%%%%%%%%%%%%%%%%%%%%%%%%%%%%%%%%%%%%%%%%%%%%%%%%%%%%%%%%%%%%%%%%%%%%%%%%%%%%%%%%%%%%%%%%
%%%%%%%%%%%%%%%%%%%%%               Carleson's estimate and  the  boundary Harnack principle         %%%%%%%%%%%%%%%%%%%%%%%%%%%%%%%%
%%%%%%%%%%%%%%%%%%%%%%%%%%%%%%%%%%%%%%%%%%%%%%%%%%%%%%%%%%%%%%%%%%%%%%%%%%%%%%%%%%%%%%%%%%%%%%%%%%%%%%%%%%%%%%%%%%%%%%%%%%%%%%%%%%%%%

\section{Carleson's estimate and  the  boundary Harnack principle}\label{s:bhp}

In this section we prove Carleson's estimate and  the  boundary Harnack principle. We start with a simple result following directly from Theorem \ref{t:explicit-equiv}. Recall that  $\eps_1\in (0,1/8]$ is the constant in Proposition \ref{p:explicit-1}.  

\begin{lemma}\label{l:Carleson-1}
There exists a constant $K_0>4$ such that for all $x\in D$ with $\delta_D(x)<\eps_1 \wh{R}/K_0$ it holds that
\begin{align*}
\P_x(\tau_{B_{D}\left(x, (K_0+1)\delta_D(x)\right)}<\zeta) < \frac{1}{2}.
\end{align*}
\end{lemma}

\begin{proof}
Let $K_0>4$ to be chosen later. 
For $x\in D$ with $\delta_D(x)<\eps_1 \wh{R}/K_0$, let $Q_x\in \partial D$ be a point such that $|x-Q_x|=\delta_D(x)$. 
Then it holds that $B_D(Q_x, K_0\delta_D(x))\subset B_D(x, (K_0+1)\delta_D(x))$. This inclusion and 
Theorem \ref{t:explicit-equiv} with $R=K_0\delta_D(x)/\eps_1<\wh{R}$ imply that
\begin{align*}
&\P_x(\tau_{B_{D}\left(x, (K_0+1)\delta_D(x)\right)}<\zeta)\le \P_x\left(X_{\tau_{B_D\left(Q_x, K_0\delta_D(x)\right)}}\in D\right)
\le c_1 \left(\frac{\delta_D(x)}{K_0\delta_D(x)/\eps_1}\right)^p = c_2 K_0^{-p}.
\end{align*}
Choosing $K_0:=(2c_2)^{1/p}\vee 5 >4$, we have proved the lemma.\end{proof}

We now give the proof of   Carleson's estimate. Proposition \ref{p:bound-for-integral-new} will be used in 
an essential way. 
\begin{thm}\label{t:Carleson}\!{\rm(Carleson's estimate)}	 
There exists  $C\ge1$ such that for any $Q \in \partial D$, $r \in (0, \wh R]$ and non-negative Borel function $f$ on $D$ which is harmonic in $B_D(Q,r)$ with respect to $X$ and vanishes continuously on $\partial D \cap B(Q,r)$, we have
	\begin{align}\label{e:Carleson}
		f(x) \le C f(z_0) \quad \text{for all} \;\, x \in B_D(Q,r/2),
	\end{align}
where $z_0 \in B_D(Q,8r/9)$ is any point with $\delta_D(z_0) \ge r/8$.
\end{thm}
\begin{proof} 
Let $\beta\in (\overline{\beta}, 1\wedge \alpha)$ be such that  $p \in [(\alpha-1)_+, \alpha-\beta)$. 
Let $Q \in \partial D$, $r \in (0, \wh R]$, $z_0 \in B_D(Q,8r/9)$ with $\delta_D(z_0) \ge r/8$ 
and let $f$ be a non-negative Borel function on $D$ which is harmonic in $B_D(Q,r)$ and vanishes continuously 
on $\partial D \cap B(Q,r)$.  
Recall that $\eps_1 \in (0,1/8]$ is  the constant in Proposition \ref{p:explicit-1}.  
By the Harnack inequality (\!\!\cite[Theorem 1.1]{KSV23})
and a standard chain argument, it suffices to prove \eqref{e:Carleson} for
$x\in B_D(Q, \eps_1 r / (48K_0))$ where 
$K_0 >4$ is the constant in Lemma \ref{l:Carleson-1}. 
Moreover, we also deduce from \cite[Theorem 1.1]{KSV23} and a standard chain argument that there exist  constants $c_1,\gamma>0$ independent of $Q,r,f$ and $z_0$ such that 
	\begin{align}\label{e:Carleson-chain}
 		f(x) \le  c_1 (\delta_D(x)/r)^{-\gamma} f(z_0) \quad \text{for all} \;\, x\in 
		B_D(Q, \eps_1 r / (24K_0)).
	\end{align} 
In the following, the constants  $c_i$ are always independent of $Q,r, f$ and $z_0$.

Choose $\lambda>0$ so that
\begin{align}\label{e:choice-lambda}
\lambda  \in (0,  \frac{\alpha-p}{d+\alpha}\wedge \frac{\alpha-\beta-p}{d+\alpha-2\beta}) \subset (0, 1).
\end{align}

Define
	\begin{align*}
		U_1:=B(z_0, \delta_D(z_0)/8), \qquad U_2:=B(z_0, \delta_D(z_0)/4)
	\end{align*}
and for $x\in B_D(Q, \eps_1 r / (12K_0))$,
	\begin{align*}
		V_1(x):=B_D(x, (K_0+1)\delta_D(x)), \qquad V_2(x):= B_D(x, (2K_0+2)r^{1-\lambda} \delta_D(x)^\lambda).
	\end{align*}
First note that,  since for all $y \in \overline U_1$, 
	$$
		|y-Q| \le |y-z_0|+|z_0-Q| \le \delta_D(z_0)/8+ |z_0-Q| \le 9 |z_0-Q| /8 <r,
	$$
we have $\overline{U_1} \subset B_D(Q,r)$.
We also have that $V_1(x) \subset V_2(x) \cap B_D(Q,r)$. Indeed, for $w\in V_1(x)$, 
\begin{align}
\label{e:V1}
|w-Q|\le |w-x|+|x-Q|\le (K_0+1)\delta_D(x)+|x-Q|\le (K_0+1)\frac{\eps_1 r}{12K_0}+\frac{\eps_1 r}{12 K_0}<r,
\end{align}
and since $\delta_D(x)< r^{1-\lambda}\delta_D(x)^{\lambda}$, also $V_1(x)\subset V_2(x)$.

By Lemma \ref{l:Carleson-1}, we see that 
	\begin{align}\label{e:Carleson-lifetime}
		\P_x (\tau_{V_1(x)}<\zeta) < 1/2
		 \quad \text{for all} \;\, x\in B_D(Q, \eps_1r / (24K_0)).
	\end{align}

Pick any $x \in B_D(Q, \eps_1r / (24K_0))$.  Since $V_1(x) \subset B_D(Q,r)$, by the harmonicity of $f$,  we have
	\begin{align}\label{e:split-f(x)}
		f(x) &= \E_x \left[ f(X_{\tau_{V_1(x)}}) ; X_{\tau_{V_1(x)}} \in V_2(x)\right] 
			+ \E_x \left[ f(X_{\tau_{V_1(x)}}) ; X_{\tau_{V_1(x)}} \in (D\setminus V_2(x))\cap U_2	\right]\nn\\
		& \quad  + \E_x \left[ f(X_{\tau_{V_1(x)}}) ; X_{\tau_{V_1(x)}} \in (D\setminus V_2(x))\cap U_2^c	\right]=:I_1+I_2+I_3.
	\end{align}
We divide the estimate of the right-hand side above into several steps.

\noindent \textbf{Step 1:} By \eqref{e:Carleson-lifetime}, 
	\begin{align}\label{e:Carleson-I1}
		I_1 \le  \Big( \sup_{y \in V_2(x)} f(y)\Big)\, 
		\P_x \big( X_{\tau_{V_1(x)}} \in V_2(x)\big) \le 2^{-1} \sup_{y \in V_2(x)} f(y).
	\end{align}
	
\medskip
\noindent \textbf{Step 2:} We show now that 
\begin{align}\label{e:step2-carleson}
f(z_0)\ge c_2 r^{\alpha}\int_{D\setminus U_2} f(y)J(z_0,y)\, dy.
\end{align}
Let $z\in U_1$ and $y\in D\setminus U_2$. Then $|z-y|\le |z-z_0|+|z_0-y|\le 2|z_0-y|$, and $\delta_D(z)\le (9/8)\delta_D(z_0)$. 
Therefore, by using \eqref{e:J-Phi} in the first, \eqref{e:Phi-scaling-2} and \eqref{e:without-A0} in the second, and $d+\alpha>2\beta$ in the fourth inequality, we see that
\begin{align*}
\frac{J(z_0,y)}{J(z,y)}&\le C_0^2\, \frac{|z_0-y|^{-d-\alpha}\Phi\left(\frac{(|z_0-y|\wedge A_0)^2}{(\delta_D(z_0)\wedge A_0)(\delta_D(y)\wedge A_0)}\right)}{|z-y|^{-d-\alpha}\Phi\left(\frac{(|z-y|\wedge A_0)^2}{(\delta_D(z)\wedge A_0)(\delta_D(y)\wedge A_0)}\right)}\\
&\le c_3 \frac{|z-y|^{d+\alpha}}{|z_0-y|^{d+\alpha}}\left(1+\frac{|z_0-y|^{2\beta}}{|z-y|^{2\beta}}\, 
\frac{\delta_D(z)^{\beta}}{\delta_D(z_0)^{\beta}}\right) \\
&\le c_3\left(\frac{|z-y|^{d+\alpha}}{|z_0-y|^{d+\alpha}}+\frac{|z-y|^{d+\alpha-2\beta}}{|z_0-y|^{d+\alpha-2\beta}}
\, \left(\frac{9}{8}\right)^{\beta}\right)\le c_4\left(2^{d+\alpha}+2^{d+\alpha-2\beta}\right)=:c_5.
\end{align*}
Further, since $f$ is non-negative and harmonic in $B_D(Q,r)$ and $\overline U_1 \subset B_D(Q,r)$,
the above estimate implies
\begin{align*}
f(z_0)&\ge \E_{z_0}\left[ f(X_{\tau_{U_1}}):\, X_{\tau_{U_1}}\in D\setminus U_2 \right]
	=\E_{z_0}\Big[\int_0^{\tau_{U_1}}\int_{D\setminus U_2} f(y)J(X_s, y)\, dy \, ds\Big]\\
	&\ge c_5^{-1}\E_{z_0}\Big[\int_0^{\tau_{U_1}}\int_{D\setminus U_2} f(y)J(z_0, y)\, dy \, ds\Big]
	=c_5^{-1} (\E_{z_0}\tau_{U_1})\int_{D\setminus U_2} f(y)J(z_0, y)\, dy\\
	&\ge c_6\delta_D(z_0)^{\alpha}\int_{D\setminus U_2} f(y)J(z_0, y)\, dy\ge c_7 r^{\alpha}\int_{D\setminus U_2} f(y)J(z_0, y)\, dy.
\end{align*}
In the last line we used 
Lemma \ref{l:exit-time-lower-killed}.

\medskip
\noindent \textbf{Step 3:} In this step we prove that
\begin{align*}
I_2&\le c_8 f(z_0)\left(\frac{\delta_D(x)}{r}\right)^{\alpha-\beta}.
\end{align*}
Let $w\in V_1(x)$ and $y\in U_2$. Then $|Q-y| >3\delta_D(z_0)/{4} \ge 3r/(32)$ and 
$$
|Q-w|\le |Q-x|+|x-w|\le \frac{\eps_1 r}{24K_0}+(K_0+1)\delta_D(x)\le \frac{K_0+2}{24K_0}\eps_1 r\le \frac{K_0+2}{24K_0}\frac{r
}{8} <\frac{r}{64}.
$$
Hence, by using this and $\delta_D(z_0)\ge r/8$, 
\begin{align}\label{e:pdist}
|w-y|&\ge |Q-y|-|Q-w|\ge \frac{3}{32}r-
\frac{1}{64}r= \frac{5}{64}r. 
\end{align} 
By the Harnack inequality (\!\!\cite[Theorem 1.1]{KSV23}), 
we have that $f(z)\le c_9 f(z_0)$ for every $z\in U_2$. Therefore, by using 
Lemma \ref{l:levy-consequence}(i) with \eqref{e:pdist} in the second inequality, the inclusion $V_1(x)\subset B_D(Q, (K_0+2)\delta_D(x))$
and \eqref{e:Phi-scaling}
 in the third, and Proposition \ref{p:bound-for-integral-new}  in the fourth  (note that $-\beta>p-\alpha$), we get
\begin{align*}
I_2 &= \E_x \int_0^{\tau_{V_1(x)}}\int_{(D\setminus V_2(x))\cap U_2} f(y)J(X_s, y)\, dy\, ds \le c_9 f(z_0)
\P_x (X_{\tau_{V_1(x)}} \in U_2)\\
&\le c_{10} f(z_0)r^{-\alpha} \E_x \int_0^{\tau_{V_1(x)}}\Phi\left(\frac{r}{\delta_D(X_s)}\right)\, ds \\
&\le c_{11} f(z_0)r^{-\alpha}r^{\beta}\E_x \int_0^{\tau_{B_D(Q, (K_0+2)\delta_D(x))}}\delta_D(X_s)^{-\beta}\, ds\\
&\le c_{12} f(z_0)r^{-\alpha+\beta} \delta_D(x)^p ((K_0+2)\delta_D(x))^{\alpha-\beta-p}=c_{13} f(z_0)\left(\frac{\delta_D(x)}{r}\right)^{\alpha-\beta}.
\end{align*} 

\medskip
\noindent \textbf{Step 4:} In this step we prove that
\begin{align*}
I_3\le c_{14}f(z_0)\left(\Big(\frac{\delta_D(x)}{r}\Big)^{\alpha-\lambda(d+\alpha)}
+\Big(\frac{\delta_D(x)}{r}\Big)^{\alpha-\lambda(d+\alpha-2\beta)-\beta}\right).
\end{align*}

Let $w\in V_1(x)$ and $y\in D\setminus V_2(x)$. Since 
$\delta_D(w)\le |x-w| + \delta_D(x) <(K_0+2)\delta_D(x) <(K_0+2)r^{1-\lambda}\delta_D(x)^{\lambda}$, we have
\begin{align*}
|w-y|\ge |x-y|-|w-x|\ge (2K_0+2)r^{1-\lambda}\delta_D(x)^{\lambda}-(K_0+2)r^{1-\lambda}\delta_D(x)^{\lambda}\ge K_0^{1-\lambda}r^{1-\lambda}\delta_D(w)^{\lambda}.
\end{align*}
Therefore, $r\le K_0^{\lambda-1}r^{\lambda}\delta_D(w)^{-\lambda}|w-y|\le r^{\lambda}\delta_D(w)^{-\lambda}|w-y|$ and \eqref{e:V1}, implying that
\begin{align}\label{e:z0y}
|z_0-y|&\le |w-y|+|w-Q|+|Q-z_0|\le |w-y|+2r\le (1+2r^{\lambda}\delta_D(w)^{-\lambda})|w-y| \nn\\
&\le 3r^{\lambda}\delta_D(w)^{-\lambda}|w-y|, \quad \text{ for all }w\in V_1(x) \text{ and }y\in D\setminus V_2(x).
\end{align}
Using \eqref{e:J-Phi} in the first,    \eqref{e:Phi-scaling-2} and \eqref{e:without-A0} in the second, 
and  \eqref{e:z0y}  and $d+\alpha>2\beta$ in the third inequality, we see that
\begin{align*}
\frac{J(w,y)}{J(z_0,y)}&\le C_0^2 \, \frac{|w-y|^{-d-\alpha}\Phi\left(\frac{(|w-y|\wedge A_0)^2}{(\delta_D(w)\wedge A_0)(\delta_D(y)\wedge A_0)}\right)}{|z_0-y|^{-d-\alpha}\Phi\left(\frac{(|z_0-y|\wedge A_0)^2}{(\delta_D(z_0)\wedge A_0)(\delta_D(y)\wedge A_0)}\right)}\\
&\le c_{15} \frac{|z_0-y|^{d+\alpha}}{|w-y|^{d+\alpha}}\left(1+\frac{|w-y|^{2\beta}}{|z_0-y|^{2\beta}}\, 
\frac{\delta_D(z_0)^{\beta}}{\delta_D(w)^{\beta}}\right) \\
&\le c_{15}\left((3r^{\lambda}\delta_D(w)^{-\lambda})^{d+\alpha}+(3r^{\lambda}\delta_D(w)^{-\lambda})^{d+\alpha-2\beta}\frac{r^{\beta}}{\delta_D(w)^{\beta}}\right)\\
&\le c_{16}\left(\Big(\frac{r}{\delta_D(w)}\Big)^{\lambda(d+\alpha)}+\Big(\frac{r}{\delta_D(w)}\Big)^{\lambda(d+\alpha-2\beta)+
\beta}\right).
\end{align*}
We now use the above estimate in $I_3$ and apply Proposition \ref{p:bound-for-integral-new} in a similar way as in Step 3. 
Note that the use of Proposition \ref{p:bound-for-integral-new} in the  third inequality is justified because of the choice of $\lambda$ in \eqref{e:choice-lambda}.  
Using \eqref{e:step2-carleson} in  the second inequality below we have
\begin{align*}
I_3&= \E_x \int_0^{\tau_{V_1(x)}}\int_{(D\setminus V_2(x))\cap U_2^c} f(y)J(X_s, y)\, dy\, ds\\
&\le c_{16} \E_x \int_0^{\tau_{V_1(x)}}\left(\Big(\frac{r}{\delta_D(X_s)}\Big)^{\lambda(d+\alpha)}+\Big(\frac{r}{\delta_D(X_s)}\Big)^{\lambda(d+\alpha-2\beta)+\beta}\right)\int_{D\setminus U_2}f(y)J(z_0,y)\, dy\\
&\le c_{16}c_2^{-1}f(z_0)r^{-\alpha}\left(r^{\lambda(d+\alpha)}\E_x \int_0^{B_D(Q,(K_0+2)\delta_D(x))} \delta_D(X_s)^{-\lambda(d+\alpha)}\, ds\right.\\
&\qquad \qquad \left. +r^{\lambda(d+\alpha-2\beta)+\beta}\E_x \int_0^{B_D(Q,(K_0+2)\delta_D(x))} 
\delta_D(X_s)^{-\lambda(d+\alpha-2\beta)-\beta}\, ds \right)\\
&\le c_{17}f(z_0)r^{-\alpha}\left(r^{\lambda(d+\alpha)}\delta_D(x)^p ((K_0+2)\delta_D(x))^{\alpha-p-\lambda(d+\alpha)}\right.\\
&\qquad \qquad \left. +r^{\lambda(d+\alpha-2\beta)+\beta} \delta_D(x)^p ((K_0+2)\delta_D(x))^{\alpha-p-\lambda(d+\alpha-2\beta)-\beta} \right)\\
&\le c_{18}f(z_0)\left(\Big(\frac{\delta_D(x)}{r}\Big)^{\alpha-\lambda(d+\alpha)}
+\Big(\frac{\delta_D(x)}{r}\Big)^{\alpha-\lambda(d+\alpha-2\beta)-\beta}\right).
\end{align*}
\medskip

\medskip
\noindent \textbf{Step 5:}
Recall that $\delta_D(x)/r\le \eps_1/(24K_0)<\eps_1$. 
Since $\alpha-\lambda(d+\alpha)>0$ and $\alpha-\lambda(d+\alpha-2\beta)-\beta>0$ by the choice of $\lambda$ in \eqref{e:choice-lambda}, 
\begin{align*}
&c_8 \Big(\frac{\delta_D(x)}{r}\Big)^{\alpha-\beta} +c_{14}\left(\Big(\frac{\delta_D(x)}{r}\Big)^{\alpha-\lambda(d+\alpha)}
+\Big(\frac{\delta_D(x)}{r}\Big)^{\alpha-\lambda(d+\alpha-2\beta)-\beta}\right)
\\
&\le c_8 \eps_1^{\alpha-\beta}+
c_{14}(
\eps_1^{\alpha-\lambda(d+\alpha)}+\eps_1^{\alpha-\lambda(d+\alpha-2\beta))
-\beta})=:c_{19}.
\end{align*}
Finally, by \eqref{e:split-f(x)}, \eqref{e:Carleson-I1}, steps 3 and 4, 
we finish the estimate of the right-hand side in \eqref{e:Carleson-induction} and obtain that
\begin{align}\label{e:Carleson-induction}
		f(x) \le 2^{-1} \sup_{y \in V_2(x)} f(y) + c_{19} f(z_0)  \quad \text{for all} \;\, x\in B_D(Q, \eps_1r / (24K_0)).
\end{align}

Now we prove that \eqref{e:Carleson} holds for all $x\in B_D(Q, \eps_1r / (48K_0))$ with 
	$$
		C=M:=3c_{19} + c_1 \bigg(\frac{24K_0(2K_0+2)}{a_0\eps_1}\bigg)^{\gamma/\lambda},
	$$
where
 	$$
 		a_0:=2^{-1}\bigg(\sum_{n=0}^\infty \bigg(\frac{3}{4}\bigg)^{n\lambda/ \gamma} \bigg)^{-1}.
 	$$
Suppose this fails. Then there exists $x_1 \in  B_D(Q, \eps_1r / (48K_0))$  such that  $f(x_1) >Mf(z_0)$.  In the following, we  construct a sequence  $(x_n)_{n \ge 2}$ in  $B_D(Q, \eps_1r / (24K_0))$  such that for all $n \ge 2$, 
	\begin{align}\label{e:Carleson-claim}
		& |x_{n}-x_{n-1}|<\frac{a_0\eps_1 r}{24K_0}\bigg(\frac{3}{4}\bigg)^{(n-2)\lambda/\gamma} \;\;\text{ and } \;\; f(x_{n}) 
			\ge f(x_1)\bigg(\frac{4}{3}\bigg)^{n-1}.
 	\end{align}
This leads to a contradiction since $f$ is bounded. 

By \eqref{e:Carleson-induction},  since 
$M>3c_{19}$, there exists $x_2 \in  \overline{V_2(x_1)}$ such that 
 	$$
 		f(x_2) \ge 2(f(x_1) - c_{19} f(z_0)) \ge (4/3) f(x_1).
 	$$
Note that 
$\delta_D(x_1) /r\le  (c_1 f(z_0)/f(x_1))^{1/\gamma}  <  (c_1/M)^{1/\gamma}$ 
by \eqref{e:Carleson-chain}. Thus, we have
 	\begin{align*}
	|x_2-x_1| 
		\le
		 (2K_0+2) r^{1-\lambda} \delta_D(x_1)^{\lambda} < (2K_0+2)(c_1/M)^{\lambda/\gamma}
 			r <a_0\eps_1 r/(24K_0)
 	\end{align*}
so that
	$$
 		|x_2-Q|\le |x_1-Q| + |x_1-x_2| < \eps_1 r/ (24K_0)
 	$$
by the triangle inequality. Hence, $x_2 \in B_D(Q,\eps_1 r/ (24K_0))$ and \eqref{e:Carleson-claim}  holds for $n=2$. Next, assume that $x_n \in B_D(Q,\eps_1 r/ (24K_0))$, $1\le n \le k$, are chosen to satisfy \eqref{e:Carleson-claim} for all $1\le n \le k$, for some $k \ge 2$. By \eqref{e:Carleson-induction}, since $f(x_k)\ge f(x_1) >Mf(z_0)$, there exists $x_{k+1} \in
\overline{ V_2(x_k)}
$ such that 
 	$$
 		f(x_{k+1}) \ge 2(f(x_k) - c_{14} f(z_0)) \ge (4/3) f(x_k).
 	$$
Since 
	$$
		\frac{\delta_D(x_k)}{r} \le \left(c_1\frac{f(z_0)}{f(x_k)}\right)^{1/\gamma} 
		\le    \left(c_1 \frac{f(z_0)}{f(x_1) (4/3)^{k-1} }\right)^{1/\gamma}   
		<  (c_1/M)^{1/\gamma} \bigg(\frac{3}{4}\bigg)^{(k-1)\lambda/\gamma}
	$$ 
by \eqref{e:Carleson-chain} and the induction hypothesis, we have
 	\begin{align*}
			&|x_{k+1}-x_k| \le (2K_0+2)r\bigg(\frac{\delta_D(x_k)}{r}\bigg)^{\lambda}\\
 		&< (2K_0+2) (c_1/M)^{\lambda/\gamma} r  \bigg(\frac{3}{4}\bigg)^{(k-1)\lambda/\gamma}
 			< \frac{a_0\eps_1 r}{24K_0}\bigg(\frac{3}{4}\bigg)^{(k-1)\lambda/\gamma}.
	 \end{align*}
Using this and the induction hypothesis, we get
	\begin{align*}
		|x_{k+1}-Q| \le
		 |x_1-Q| + \sum_{n=2}^{k+1} |x_{n}-x_{n-1}|< \frac{\eps_1 r}{24K_0}\bigg( \frac{1}{2} 
			+ a_0\sum_{n=2}^{k+1} \bigg(\frac{3}{4}\bigg)^{(n-2)\lambda/\gamma} \bigg)  <  \frac{\eps_1 r}{24K_0}.
	\end{align*}
Therefore, $x_{k+1} \in B_D(Q, \eps_1r/(24K_0))$ and we deduce that \eqref{e:Carleson-claim} holds for all $n$ by the induction. The proof is complete.

\end{proof}

Using Proposition \ref{p:bound-for-integral-new} and Carleson's estimate, we can prove Theorem \ref{t:BHPnew} by using a similar approach to the one in \cite{KSV25}.

\begin{proof}[Proof of Theorem \ref{t:BHPnew}]
Let $\beta\in (\overline{\beta}, 1\wedge \alpha)$ be such that  $p \in [(\alpha-1)_+, \alpha-\beta)$. 
Recall that $a \in (0,1/4]$ is the constant in 
\eqref{e:W(Q,R)} and that,  for any $Q\in \partial D$ and $R\in (0,\wh R]$, $w_{Q,R} \in B_D(Q,R)$  
is
the point in \eqref{e:W(Q,R)} so that $B(w_{Q,R}, a R) \subset (D(R)\setminus  D(R/4)) \cap B_D(Q,5R/4)$.

 Fix  $r \in (0, \wh R]$ and  $Q \in \partial D$, 
and let $f$ be a non-negative Borel function on $D$ which is harmonic in $B_D(Q,r)$ and vanishes continuously 
on $\partial D \cap B(Q,r)$.  
Recall that $\eps_1 \in (0,1/8]$ is  the constant in Proposition \ref{p:explicit-1}.  
By the Harnack inequality (\!\!\cite[Theorem 1.1]{KSV23}) 
and a standard chain argument, it suffices to prove the theorem for 
$x\in B_D(Q, \eps_1 r / 2)$.

Let 
$U={B_D(Q,\eps_1r)}$ and $V={B_D(Q,r/2)}$.
By the Harnack inequality (\!\!\cite[Theorem 1.1]{KSV23})
and 
Theorem \ref{t:explicit-equiv}, we have
for $x\in B_D(Q, \eps_1 r / 2)$,
\begin{align}\label{e:TAMSe6.37}
f(x)&=
\E_x[f(X_{\tau_{U}})]\ge \E_x[f(X_{\tau_{U}}); X_{\tau_{U}}
\in B(w_{Q,r}, a r) ]
\nonumber\\
&
\ge c_1f(w_{Q,r})\P_x(X_{\tau_{U}}\in 
B(w_{Q,r}, a r) )
\ge c_2f(w_{Q,r}) (\delta_D(x)/r)^p.
\end{align}

Fix a $w_0 =w_{Q,\eps_1r/2}$ so that 
$B(w_0, a \eps_1r/2) \subset (D(\eps_1r/2)\setminus  D(\eps_1r/8)) \cap B_D(Q,5\eps_1r/8)$.
Then, 
$$
J(z, y) \asymp J(w_0,y),  \quad \text{ for all } z \in B(w_0, a \eps_1r/2), \,\,y \in D\setminus V.
$$
Thus, 
by Lemma \ref{l:exit-time-lower-killed}, 
\begin{align}\label{e:POTAe7.27}
f(w_0)&\ge \E_{w_0}
\left[f(X_{\tau_{U}}); X_{\tau_{U}}\notin V
\right]
\ge \E_{w_0}\int^{\tau_{B(w_0, a \eps_1r/2)}}_0\int_{D\setminus V}
J(X_t, y)f(y)dydt\nonumber\\
&\ge c_{3}\E_{w_0}\tau_{B(w_0, a \eps_1r/2)}\int_{D\setminus V}J({w_0}, y)f(y)dy \ge c_{4} r^\alpha\int_{D\setminus V}J({w_0}, y)f(y)dy.
\end{align}
 
Note that $\eps_1r/8 \le \delta_D(w_0) \le \eps_1r/2$ and 
$$|z-y|\asymp |{w_0}-y| \ge (1-\eps_1)r/2 \ge \eps_1r \ge \delta_D(z) \text{ for any } z\in U,\,\,y\in D\setminus 
V.$$ Thus, by \eqref{e:J-Phi}, \eqref{e:Phi-scaling} and \eqref{e:without-A0}, for $z\in U$ and $y\in D\setminus 
V$,
\begin{align}\label{e:new-estimate-for-J}
J(z,y)&\le c_5 \frac{ \Phi(   \frac{(|y-w_0|\wedge A_0)^2}{(\delta_D(y)\wedge A_0) (\delta_D(z)\wedge A_0) }) }{|y-w_0|^{d+\alpha}}  
\le c_6  \left(\frac{r\wedge A_0}{\delta_D(z)\wedge A_0}\right)^{\beta}
\frac{ \Phi(   \frac{(|y-w_0|\wedge A_0)^2}{(\delta_D(y)\wedge A_0) (r\wedge A_0) }) }{|y-w_0|^{d+\alpha}}  \\
&   \le c_7
(r/\delta_D(z))^{\beta}J({w_0},y).
\end{align}

Combining \eqref{e:new-estimate-for-J} with  \eqref{e:POTAe7.27}
and Proposition \ref{p:bound-for-integral-new},  
 we now have that for $x\in B_D(Q, \eps_1 r / 2)$,
\begin{align}\label{e:POTAe7.29}
&\E_x\left[f(X_{\tau_{U}}); X_{\tau_{U}}\notin 
V\right]
=\E_x\int^{\tau_{U}}_0\int_{D\setminus V}
J(X_t, y)f(y)dydt\nonumber\\
&\le c_{7} r^\beta \E_x\int_0^{\tau_{U}} \delta_D(X_t)^{-\beta}dt \int_{D\setminus V}J(w_0, y)f(y)dy\nn\\
&\le c_{8} r^{\beta-\alpha} f(w_0)\E_x\int_0^{\tau_{U}} \delta_D(X_t)^{-\beta}dt \le c_9  f(w_0) (\delta_D(x)/r)^p.
\end{align}

On the other hand, by the Harnack inequality (\!\!\cite[Theorem 1.1]{KSV23})  and Theorem \ref{t:Carleson} in the first inequality, 
and Theorem \ref{t:explicit-equiv}  in the second, 
we have that for  $x\in B_D(Q, \eps_1 r / 2)$,
\begin{align}\label{e:POTAe7.30}
\E_x\left[f(X_{\tau_{U}}); X_{\tau_{U}}\in V
\right]&\le c_{10}f(w_0)
\P_x\left(X_{\tau_{U}}\in V
\right) \le c_{11} f(w_0) (\delta_D(x)/r)^p.
\end{align}
Combining \eqref{e:POTAe7.29} and \eqref{e:POTAe7.30}, we get that for $x\in B_D(Q, \eps_1 r / 2)$,
\begin{align*}
f(x)&=
\E_x\left[f(X_{\tau_{U}}); X_{\tau_{U}}\in V
\right]+
\E_x\left[f(X_{\tau_{U}}); X_{\tau_{U}} \in D \setminus  V
\right]\nonumber\\
&\le c_{11}f(w_0) (\delta_D(x)/r)^p+c_{9}f(w_0) (\delta_D(x)/r)^p \le c_{12} f(w_0) (\delta_D(x)/r)^p.
\end{align*}
 This with  
\eqref{e:TAMSe6.37}  and the Harnack inequality (\!\cite[Theorem 1.1]{KSV23}) implies that 
$f(x) \asymp   f(w_0)  (\delta_D(x)/r)^p$ for $x\in B_D(Q, \eps_1 r / 2)$.
Therefore we conclude that 
$$
 \qquad \qquad \qquad \qquad \qquad \qquad
\frac{f(x)}{f(y)}\asymp \frac{ \delta_D(x)^p}{ \delta_D(y)^{p}}, \quad \text{for} \;\, x,y\in B_D(Q, \eps_1 r / 2).
$$
\end{proof}

%%%%%%%%%%%%%%%%%%%%%%%%%%%%%%%%%%%%%%%%%%%%%%%%%%%%%%%%%%%%%%%%%%%%%%%%%%%%%%%%%%%%%%%%%%%%%%%%%%%%%%%%%%%%%%%%%%%%%%%%%%%%%%%%%%%%%
%%%%%%%%%%%%%%%%%%%%%%%%%%%%%%%%%%                         Applications                            %%%%%%%%%%%%%%%%%%%%%%%%%%%%%%%%%% 
%%%%%%%%%%%%%%%%%%%%%%%%%%%%%%%%%%%%%%%%%%%%%%%%%%%%%%%%%%%%%%%%%%%%%%%%%%%%%%%%%%%%%%%%%%%%%%%%%%%%%%%%%%%%%%%%%%%%%%%%%%%%%%%%%%%%%
\section{Applications}\label{s:appl}

Let 
\begin{align}
\label{e:Ada}
A_{d,-\alpha}:=2^\alpha \pi^{-d/2} \Gamma((d+\alpha)/2)/|\Gamma(-\alpha/2)|
\end{align}
 and let
\begin{align*}
	J^\alpha(x,y):= A_{d,-\alpha} |x-y|^{-d-\alpha}, \quad x,y \in \R^d,
\end{align*}
be the jump kernel of the fractional Laplacian $\Delta^{\alpha/2}:=-(-\Delta)^{\alpha/2}$ on $\R^d$.

\begin{defn}
	\rm 	We say that an open set $D\subset \R^d$ is \textit{half-space-like} if there are constants $a<b$ such that $\{(\wt x,x_d)\in \R^d: x_d>b\} \subset D \subset \{(\wt x,x_d)\in \R^d: x_d>a\}$.
\end{defn}

\subsection{Elliptic problems with nonlocal Neumann boundary condition}\label{ss:elliptic-problems}
This example is a continuation of \cite[Subsection 6.1]{CKSV25}. 
Throughout this subsection, we assume that $\alpha \in (0,2)$ and that   $D\subset \R^d$ is a $C^{1,\eta}$ open set  that is either bounded,  or half-space-like, or the complement of a bounded set.
Note that this implies that ${\rm dim}_{\rm A}(\partial D)=d-1$.

  Consider
the elliptic nonlocal problem with a \textit{nonlocal Neumann boundary condition} 
\begin{align}\label{e:nonlocal-Neumann}
	\begin{cases} 
		{\Delta^{\alpha/2} u=f}  & \text{ in $D$,}\\[1pt]
		{\sN^{\alpha/2}_D u =0}  & \text{ in $\R^d \setminus 
			\overline D$,} 
	\end{cases}    
\end{align}
where the \textit{nonlocal normal derivative} $\sN^{\alpha/2}_D $ is defined by
\begin{align*}
	N^{\alpha/2}_D u (x) = \int_D (u(x)-u(y)) J^\alpha(x,y)dy.
\end{align*}
This problem was introduced in \cite{DRV17} as a nonlocal analogue of the classical Neumann problem and has been studied extensively in subsequent works.  One of the key approaches to studying  \eqref{e:nonlocal-Neumann} is to  reformulate it as a problem on $D$ involving 
the regional-type nonlocal operator
\begin{align}\label{e:regional-nonlocal}
	\Delta^{\alpha/2,N}_D u(x)=p.v. \int_D (u(y)-u(x)) K^\alpha_D(x,y)dy,
\end{align}
where
\begin{align}
	K^\alpha_D(x,y)&:= J^\alpha(x,y) + \int_{\overline D^c} \frac{J^\alpha(x,z)J^\alpha(z,y)}{ \wh \kappa_{D}(z) } \, dz, \quad x,y \in D, \nn\\
	 \wh \kappa_D(z)&:=\int_D J^\alpha(z,w)dw, \quad z\in \overline D^c. \label{e:def-wh-kappa}
\end{align}
See \cite{Ab20}. Let $\sB_D^\alpha(x,y)$ denote the coefficient function associated with $K_D^\alpha(x,y)$, defined by
\begin{align*}
	\sB^\alpha_D(x,y):=\begin{cases}
		K_D^\alpha(x,y)|x-y|^{d+\alpha} &\mbox{ if $x\neq y$},\\[5pt]
	A_{d,-\alpha}&\mbox{ if $x=y$}.
	\end{cases}
\end{align*}

It was shown in \cite[Proposition 6.3]{CKSV25} that  $\sB_D^\alpha(x,y)$ satisfies \hyperlink{A}{{\bf (A)}} with 
$\Phi(r)=1\vee \log r$ 
and  $A_0= \diam(D^c)$. Note that the upper Matuszewska index of 
 $\Phi(r)=1\vee \log r$ 
is 0.

We prove that $\sB_D^\alpha(x,y)$ also satisfies  \hyperlink{B1}{{\bf (B1)}} and  \hyperlink{B2}{{\bf (B2)}}.
\begin{prop}\label{p:nonlocal-Neumann}
Suppose that $D\subset \R^d$ is a $C^{1,\eta}$ open set  that is either bounded,  or half-space-like, or the complement of a bounded set. Then $\sB^\alpha_D(x,y)$ satisfies   \hyperlink{B1}{{\bf (B1)}} and \hyperlink{B2}{{\bf (B2)}}.
\end{prop}

The proof of Proposition \ref{p:nonlocal-Neumann} is deferred to  
Subsection \ref{ss:proof-1}.

\begin{example}\label{ex:1}
\rm {\bf (Nonlocal Neumann problem with boundary Hardy potential)} Let $\alpha \in (0,2)$ and let   $D\subset \R^d$ be a $C^{1,\eta}$ open set  that is either bounded,  or half-space-like, or the complement of a bounded set. Let $\kappa \in [0,\infty)$. When $\alpha\le 1$, we assume that $\kappa>0$. 
Consider 
\begin{align}\label{e:nonlocal-Neumann-Hardy}
	\begin{cases} 
		{\Delta^{\alpha/2} u(x)=\kappa \delta_D(x)^{-\alpha}u(x) + f(x)}  & \text{ in $D$,}\\[1pt]
		u = 0 &  \text{ on $\partial D$,}\\[1pt]
		{\sN^{\alpha/2}_D u =0}    &\text{ in $\R^d \setminus 
			\overline D$.} 
	\end{cases}    
\end{align}
Following  \cite{Ab20} and \cite[Appendix A]{AFR20}, the problem \eqref{e:nonlocal-Neumann-Hardy} can be 
equivalently reformulated as the following  nonlocal problem on $D$ with a local Dirichlet boundary condition:
\begin{align}\label{e:nonlocal-Neumann-Hardy+}
	\begin{cases} 
		{\Delta^{\alpha/2,N}_D u(x)=\kappa \delta_D(x)^{-\alpha}u(x) + f(x)}  & \text{ in $D$,}\\[1pt]
		u = 0 &  \text{ on $\partial D$.}
	\end{cases}    
\end{align}
The bilinear form corresponding to the   variational formulation of  \eqref{e:nonlocal-Neumann-Hardy+} is given by
\begin{align*}
	\sE^{\alpha, N,\kappa}(u,v)&:= \frac12 \int_{D\times D} (u(x)-u(y))(v(x)-v(y))K_D^\alpha(x,y)dxdy   + \kappa \int_{D} u(x) v(x) \delta_D(x)^{-\alpha}dx,\\
	\sF^{\alpha,N,\kappa}&:=\text{$(	
	\sE^{\alpha, N,\kappa})_1
	$-closure of $C_c^2(D)$ in $L^2(D)$}.
\end{align*}

By Proposition \ref{p:nonlocal-Neumann}, the form $(\sE^{\alpha, N,\kappa}, \sF^{\alpha,N,\kappa})$ falls within our framework. In particular, by Theorem \ref{t:BHPnew},  there exists a unique $p_\kappa \in [(\alpha-1)_+, \alpha) \cap  (0,\alpha)$  such that $(\sE^{\alpha, N,\kappa}, \sF^{\alpha,N,\kappa})$  satisfies the BHP with  explicit decay rate  $\delta_D(x)^{p_\kappa}$. 
When $\alpha>1$ and $\kappa=0$, we have   $p_\kappa=\alpha-1$.
\end{example} 

Recall that  $\wh \kappa_D(z)$ is defined in \eqref{e:def-wh-kappa}.  Define
\begin{align*}
	\kappa_D(x):=\int_{D^c} J^\alpha(x,y)dy, \quad x\in D.
\end{align*}

\begin{example}\label{ex:2}
	\rm {\bf (Nonlocal Robin problem)} Let $\theta \in (0,1)$.  Consider the elliptic problem with a nonlocal Robin boundary condition (see \cite{Ab20}):
	\begin{align}\label{e:nonlocal-Robin}
		\begin{cases} 
			{\Delta^{\alpha/2} u=f}  & \text{ in $D$,}\\[1pt]
			{\theta \wh \kappa_D(z) u(z)+  (1-\theta)\sN^{\alpha/2}_D u(z) =0}  & \text{ in $\R^d \setminus 
				\overline D$.} 
		\end{cases}    
	\end{align}
	In view of \cite[Theorem 1.1]{Ab20}, and by adapting the argument in   \cite[Appendix A]{AFR20}, the problem \eqref{e:nonlocal-Robin} can be equivalently reformulated as the following  nonlocal problem on $D$:
	\begin{align}\label{e:nonlocal-Neumann-Robin+}
			{\Delta^{\alpha/2,R,\theta}_D u(x)= f(x)}  \quad \text{ in $D$,}
	\end{align}
	where
	\begin{align*}
			\Delta^{\alpha/2,R,\theta}_D u(x)&:=p.v. \int_D (u(y)-u(x))K^{\alpha,\theta}_D(x,y)dy -  \theta \kappa_D(x) u(x) ,\\
			K^{\alpha,\theta}_D(x,y)&:=\theta J^\alpha(x,y) + (1-\theta)K^\alpha_D(x,y), \quad x,y \in D.
	\end{align*}
	The corresponding bilinear form is
	\begin{align*}
		\sE^{\alpha, R,\theta}(u,v)&:= \frac12 \int_{D\times D} (u(x)-u(y))(v(x)-v(y))K_D^{\alpha,\theta}(x,y)dxdy   + \theta \int_{D} u(x) v(x) \kappa_D(x)dx,\\
		\sF^{\alpha,R,\theta}&:=\text{$(	\sE^{\alpha, R,\theta})_1$-closure of $C_c^2(D)$ in $L^2(D)$}.
	\end{align*}
	
	By  Lemma \ref{l:C1eta-killing} below, applied with  $D$ in place of $\overline D^c$,  there exists $C>0$ such that
	\begin{equation*}
		\left|\theta \kappa_D(x)   - 
	\theta	\kappa_{d,\alpha} \delta_{D}(x)^{-\alpha}  \right| \le   C\delta_{D}(x)^{-\alpha + (\eta \wedge (\alpha/2))} \quad \text{for all $x\in D$},
	\end{equation*}
	where  $\kappa_{d,\alpha}$ is defined in \eqref{e:def-kappa-alpha}. Further, for all $x\in D$ with $\delta_D(x)\ge 1$, we have $\theta \kappa_D(x) \le \theta \int_{B(x,1)^c} J^\alpha(x,y)dy = c$. Hence, the killing potential $\theta\kappa_D(x)$ satisfies  \hyperlink{K}{{\bf (K)}}. On the other hand, 	the function 
	$\sB_D^{\alpha,\theta}(x,y)$ associated with $K_D^{\alpha,\theta}(x,y)$ is  
	\begin{align*}
		\sB^{\alpha,\theta}_D(x,y)=\theta A_{d,-\alpha} + (1-\theta) \sB^\alpha_D(x,y), \quad x,y \in D.
	\end{align*} 
	Using Proposition \ref{p:nonlocal-Neumann}, 
	it is straightforward to verify that $\sB_D^{\alpha,\theta}(x,y)$ satisfies  \hyperlink{B1}{{\bf (B1)}} and \hyperlink{B2}{{\bf (B2)}}.  Consequently, the form $(\sE^{\alpha,R,\theta}, \sF^{\alpha,R,\theta})$ falls within our framework. Moreover, 
	by the monotonicity result in   
	Lemma \ref{l:constant}, the associated parameter $p_\theta$  satisfies $p_\theta \in 
	 [(\alpha-1)_+, \alpha/2) 
	$ with  $\lim_{\theta\to 0} p_\theta = (\alpha-1)_+$ and $\lim_{\theta \to 1} p_\theta = \alpha/2$. Indeed, as $\theta \to 0$ and $\theta \to 1$, respectively, the corresponding limit problems are the problem with $\kappa=0$ considered in Example \ref{ex:1} and the nonlocal Dirichlet problem, for which the boundary decay rate is known to be $\delta_D(x)^{\alpha/2}$.
\end{example}

\subsection{Trace of $\alpha$-stable processes}\label{ss:trace}

This example is a continuation of \cite[Subsection 6.2]{CKSV25}. 
Throughout this subsection, we assume that    $D\subset \R^d$ is a $C^{1,\eta}$ open set  that is either bounded,  or half-space-like, or the complement of a bounded set. 
Let $Y=(Y_t)_{t\ge 0}$ be an isotropic $\alpha$-stable process in $\R^d$ with generator  
$\Delta^{\alpha/2}$ where $\alpha \in (0,2)$.
Define a  continuous additive functional $\sA^{\overline D}=(\sA^{\overline D}_t)_{t\ge 0}$ of $Y$ by 
$$\sA^{\overline D}_t:=\int_0^t \1_{\{Y_s \in \overline D\}}ds,\quad t\ge 0.$$ 
Let $\tau^{\overline D}_t:=\inf\{s>0:\sA^{\overline D}_s>t\}$ be the right-continuous inverse of $\sA^{\overline D}$. Then the time-changed process $(X_t)_{t\ge 0}$ defined by $X_t:=Y_{\tau^{\overline D}_t}$ is a symmetric Markov process on $\overline D$,  known as the \textit{trace process} of $Y$ on $\overline D$. We refer to \cite[Chapter 5]{CF12} for basic results on trace processes. 
The  subprocess $X^{ D}$  of $X$ killed upon leaving  $D$, or equivalently, hitting $\partial D$, is referred to as the \textit{path-censored $\alpha$-stable process} on $D$. For $d=1$, path-censored $\alpha$-stable processes were studied in  \cite{KPW14}. 
When $\alpha\in (0, 1]$, $\partial D$ is polar for $X$, so $X^D$ coincides with $X$ in this case. In the remainder of this subsection, we assume $d\ge 2$ and $\alpha\in (1, 2)$.

In the case $D$  is either  half-space-like or the complement of a bounded set, the Dirichlet form of $X$ is purely jump and has no killing part.  
In the case $D$ is bounded, 
the Dirichlet form of $X$ is purely jump and  has a killing measure $V_D^\alpha(x)dx$, where
\[
V^\alpha_D(x):= \int_{ D^c}  J^\alpha(x,z)   \P_z(\sigma_{\overline D} = \infty) dz
\] 
and $\sigma_{\overline D}$ is the first hitting time of $\overline D$ for $Y$. 
Using the boundary Harnack principle for $Y$ (see, e.g., \cite[Corollary 1.9]{CS26}) and the transience of $Y$, 
one can obtain that  
$\P_z(\sigma_{\overline D} = \infty) \asymp (\delta_{D^c}(z) \wedge 1)^{\alpha/2}$ for $z\in \overline D^c$.  
Thus,  $\P_z(\sigma_{\overline D} = \infty) \le c_1 |x-z|^{\alpha/2}$ for $z\in \overline D^c$ and $x\in D$. 
Using this, we have that for $x\in D$, 
\begin{align}
	\label{e:nbabo}
	V_D^\alpha(x)\le c_1 \int_{B(x, \delta_D(x))^c} |x-z|^{-d-\alpha/2}
	dz \le c_2 \int_{\delta_{D}(x)}^\infty r^{-1-\alpha/2}
	dr
	\asymp
	\delta_{ D}^{-\alpha/2}(x).
\end{align}

 Let  $\wc J_D^\alpha(x,y)$  denote the jump kernel of $X$. Then for all $x,y \in D$,
 \begin{align*}
 	\wc J_D^\alpha(x,y) = J^\alpha(x,y) + \int_{\overline D^c\times \overline D^c} J^\alpha(x,z) g^\alpha_{\overline D^c} (z,w) J^\alpha(w,y) \, dzdw,
 \end{align*}
 where $g^\alpha_{\overline D^c} (z,w)$ denotes the Green function of the isotropic $\alpha$-stable process on $\overline D^c$; see the discussion at the beginning of \cite[Subsection 7.2]{CKSV25}. Let $\wc \sB_D^\alpha(x,y)$ denote the coefficient function associated with $\wc J_D^\alpha(x,y)$, defined by
 \begin{align*}
 \wc \sB_D^\alpha(x,y):=\begin{cases}
 	\wc J_D^\alpha(x,y)|x-y|^{d+\alpha} &\mbox{ if $x\neq y$},\\[5pt]
 		A_{d,-\alpha}&\mbox{ if $x=y$}.
 	\end{cases}
 \end{align*}
  By \cite[Proposition 6.4]{CKSV25}, $\wc \sB_D^\alpha(x,y)$ satisfies \hyperlink{A}{{\bf (A)}} with 
   $\Phi(r) = 1\vee r^{\alpha/2}$  
  and $A_0=\diam(D^c)$. We will prove that
\begin{prop}\label{p:nonlocal-trace}
	Suppose $d\ge 2$,  $\alpha\in (1, 2)$ and $D\subset \R^d$
	is a $C^{1,\eta}$ open set  that is either bounded,  or half-space-like, or the complement of a bounded set. Then $	\wc \sB_D^\alpha(x,y) $ satisfies   \hyperlink{B1}{{\bf (B1)}} and \hyperlink{B2}{{\bf (B2)}}.
\end{prop}

The proof of Proposition \ref{p:nonlocal-trace} is deferred to 
Subsection \ref{ss:proof-1}.

By \eqref{e:nbabo},  when $\alpha>1$ and $D$ is bounded, $V_D^\alpha$ satisfies the condition \hyperlink{K}{{\bf (K)}} with $\kappa_0=0$ and $\theta_3=\alpha/2$.
As an immediate consequence of Proposition \ref{p:nonlocal-trace} and Theorem \ref{t:BHPnew},  when 
$d\ge 2$ and  $\alpha\in (1, 2)$, 
the path-censored $\alpha$-stable process on $D$ satisfies the BHP with the explicit  decay rate $\delta_D(x)^{\alpha-1}$.

%%%%%%%%%%%%%%%%%%%%%%%%%%%%%%%%%%%%%%%%%%%%%%%%%%%%%%%%%%%%%%%%%%%%%%%%%%%%%%%%%%%%%%%%%%%%%%%%%%%%%%%%%%%%%%%%%%%%%%%%%%%%%%%%%%%%%
%%%%%%%%%%%%%%%%%%%%%%%%%%%%                           Appendix                         %%%%%%%%%%%%%%%%%%%%%%%%%%%%%%%%%%%%%%%%%%%%%
%%%%%%%%%%%%%%%%%%%%%%%%%%%%%%%%%%%%%%%%%%%%%%%%%%%%%%%%%%%%%%%%%%%%%%%%%%%%%%%%%%%%%%%%%%%%%%%%%%%%%%%%%%%%%%%%%%%%%%%%%%%%%%%%%%%%%
\section{Appendix}\label{s:appendix}

\subsection{Dynkin-type formula}\label{ss:dynkin}
We first check that our assumptions on $\sB(x,y)$ imply assumptions (H1)-(H4) in \cite{KSV23}. The assumption (H1) is simply the symmetry of $\sB$, and (H4) is essentially \hyperlink{B1}{{\bf (B1)}}. The assumption (H2) states that for any $a\in (0,1)$ there is a constant $C=C(a)>1$ such that $C^{-1}\le \sB(x,y)\le C$ for all $x,y\in D$ satisfying
$\delta_D(x)\wedge \delta_D(y)\ge a|x-y|$. This is straightforward from \hyperlink{A}{{\bf (A)}}. Finally, (H3) states that for any $a>0$ there is a constant $C=C(a)>0$ such that 
$\int_{D \setminus B(x, a\delta_D(x))} J(x,y)\, dy\le C\delta_D(x)^{-\alpha}. $
This is a direct consequence of  \eqref{e:jump-tail-estimate}  with $r=a\delta_D(x)$.

Recall that 
\begin{align*}
	L^0_\sB u(x)
	=\text{p.v.}\int_D (u(y)-u(x))J(x,y)\, dy 
	=\lim_{\epsilon \downarrow 0}\int_{D, |y-x|>\epsilon}(u(y)-u(x)) \frac{\sB(x,y)}{|y-x|^{d+\alpha}}dy
\end{align*}
and 
\begin{align*}
	L^{\kappa}_\sB u(x)=L^0_\sB
	u(x)-\kappa(x)u(x).
\end{align*}
Then 
%PK
%$L^0_\sB u$
%is well-defined for all $u \in C^2_c(D)$.  Moreover, if $(\sA^\kappa, \sD(\sA^\kappa))$ is the $L^2$-generator of $(\sE^\kappa,\sF^\kappa)$, then $C^{2}_c(D) \subset \sD(\sA^\kappa)$, and $\sA^\kappa u = L^\kappa_\sB u$  a.e.~in $D$ for $u \in C^2_c(D)$; see \cite[Proposition 4.2 and Corollary 4.5]{KSV23} for details. Using these facts and the strong Feller property, we obtain the following Dynkin-type formula; see  \cite[Theorem 4.8]{KSV23} for a proof.
$L^\kappa_\sB u(x)$
is well-defined and finite for all 
$u \in C^2_c(D)$ and $x \in D$; see \cite[Proposition 4.2]{KSV23}. In fact, 
%ZV One 
one
can  see easily from the proof of \cite[Theorem 4.8]{KSV23} (using \cite[Proposition 4.2(c)]{KSV23}) that,  
for any bounded open set $U$ with $\overline U \subset D$, 
$L^\kappa_\sB u(x)$
is well-defined and bounded on $\overline U$ for any bounded Borel
	function $u$ on $D$ such that $u|_{\overline U} \in C^2(\overline U)$.
Thus,  the following Dynkin-type formula holds true. see \cite[Theorem 4.8]{KSV23} for a proof.
Note that there are typos in 
	the statement of \cite[Corollary 4.4]{KSV23}: 
	 $C_c^2(D; \R^d)$ should be 
%ZV	 $C_c^2(D)$. 
$C_c^2(D)$,
	 cf. \cite[Corollary 4.5]{KSV23}.

\begin{prop}\label{p:Dynkin-formula}
	Let $U$ be a  bounded open set with $\overline U \subset D$. For any bounded Borel
	function $u$ on $D$ such that $u|_{\overline U} \in C^2(\overline U)$, we have
	\begin{align}\label{e:Dynkin-formula}
		\E_x \big[ u(X_{\tau_U})\big] = u(x) + \E_x \bigg[\int_0^{\tau_U} 
		L^\kappa_\sB
		u(X_s)ds \bigg] 
		\quad \text{for all} \;\, x \in U.
	\end{align}
\end{prop}
Here $\tau_U=\inf\{t>0:\, X_t\notin U\}$ is the first exit time of $X$ from $U$.

\begin{corollary}\label{c:Dynkin-local}
	 Let $U \subset D$ be a bounded open set and $u$ 
	 be a bounded Borel  function on $D$ such that 
	$u|_{\overline U \cap D} \in C^2(\overline U \cap D)$ and $u$  vanishes continuously at $\partial U \cap \partial D$. 
	  Suppose that either
	  $
	L^\kappa_\sB
	u(y) \ge 0$ in $U$ or $
	L^\kappa_\sB
	u(y) \le  0$ in $U$. Then  \eqref{e:Dynkin-formula} holds.
\end{corollary}
\begin{proof} 
	Let $x \in U$. For $j \ge 1$, define 
	$A_j:=\{y \in U: \delta_D(y) >2^{-j} \}$. 
	Clearly, $A_j \uparrow U$.  
	Hence, there exists $j_0 \ge 1$ such that $x \in A_j$ for all $j \ge j_0$, 
	and $\tau_{A_j}\uparrow \tau_U$. 
	Since  $u|_{\overline{A_j}} \in 
	C^2(\overline{A_j})$ for all $j\ge 1$, by Proposition \ref{p:Dynkin-formula}, we get that for all $j \ge j_0$,	
	\begin{align}\label{e:Dynkin-local-2}
		\E_x \big[ u(X_{\tau_{A_j}})\big] 
		= u(x) + \E_x \bigg[\int_0^{\tau_{A_j}} 
		L^\kappa_\sB
		u(X_s)ds \bigg].	
	\end{align}
	
	We claim that
	\begin{align}\label{e:Dynkin-formula-claim}
		\lim_{j\to \infty} u(X_{\tau_{A_j}}) = u(X_{\tau_U}) \quad \P_x \text{-almost surely}.
	\end{align}
	Indeed, if $\tau_{A_k} = \tau_U$ for some $k$, then  $\tau_{A_j}=\tau_U$ for all $j \ge k$, and  \eqref{e:Dynkin-formula-claim} follows immediately.
On the complementary event $\cap_{k=1}^\infty \{\tau_{A_k} < \tau_U\} $, we have $X_{\tau_{A_k}} \in U \setminus A_k$ and $\delta_D(X_{\tau_{A_k}}) \le 2^{-k}$ for all $k\ge 1$.  
Since $\tau_{A_j} \uparrow \tau_U$ and $X_t=\overline X_t$ for $t<\zeta$, we see that  $\lim_{j\to \infty} X_{\tau_{A_j}} = \overline X_{\tau_U - }\in \overline U  \cap \partial D =  \partial U \cap \partial D$ on $\cap_{k=1}^\infty \{\tau_{A_k} < \tau_U\} $. Further, by  quasi-left continuity, we deduce that $X_{\tau_U} \notin D$, or equivalently,  $\tau_U=\zeta$, $\P_x$-almost surely on  $\cap_{k=1}^\infty \{\tau_{A_k} < \tau_U\} $. By the  assumption that $u$ vanishes continuously at  $\partial U\cap \partial D$ and the convention $u(X_{\zeta})=0$, it follows that $$
\lim_{j \to \infty}  u(X_{\tau_{A_j}})=0 = u(X_{\tau_U})\quad \text{$\P_x$-almost surely on $\cap_{k=1}^\infty\{\tau_{A_k} < \tau_U\}$}.
$$
This completes the proof of \eqref{e:Dynkin-formula-claim}.
	
Since $u$ is bounded, by the dominated convergence theorem, we obtain from \eqref{e:Dynkin-formula-claim} that 
\begin{align*}
	\lim_{j\to \infty}\E_x\big[u(X_{\tau_{A_j}})\big]=
	\E_x\big[u(X_{\tau_U})\big].
\end{align*} 
 	
On the other hand, since $L^\kappa_\sB u$ 
is either positive in $U$  or negative in $U$, using the monotone convergence theorem, we see that  the right-hand side of \eqref{e:Dynkin-local-2} converges to that of \eqref{e:Dynkin-formula} as $j \to \infty$. Now we arrive at the result by letting $j \to \infty$ in \eqref{e:Dynkin-local-2}.
\end{proof}

\subsection{Proofs of Propositions \ref{p:nonlocal-Neumann} and \ref{p:nonlocal-trace}}\label{ss:proof-1}
Throughout this subsection, we assume that $\alpha \in (0,2)$ and that  $D\subset \R^d$ is a $C^{1,\eta}$ open set  with characteristics $(\wh R, \Lambda)$ satisfying $\wh R^\eta \le 1 \wedge (1/(4\Lambda))$, and that $D$ is either bounded,  or half-space-like, or the complement of a bounded set.  
We will use the following observation several times:
\begin{equation}\label{e:dist-D-Dc}
	|x-z| \ge  \delta_D(x) \vee \delta_{D^c}(z) \quad \text{for all} \;\, x \in D \text{ and } z \in D^c.
\end{equation}

For an open set $U\subset \R^d$ and $x,y \in U$, define
\begin{align*} 
	\delta_U^\wedge(x,y):=\delta_U(x) \wedge \delta_U(y) \quad \text{ and } \quad \delta_U^\vee(x,y):=\delta_U(x) \vee \delta_U(y).
\end{align*}

The following lemma is a particular case of  \cite[Lemma 7.3]{CKSV25}.
\begin{lemma}\label{l:technical-1}
	For any $a_1,a_2\in [-d,\infty)$ and 
	$a_3\ge 0$ with 
	  $a_1\vee a_2<a_3<d+a_1+a_2$, 
	there exists $C=C(a_1,a_2,a_3)>0$ such that for all $x,y \in D$,
	\begin{align*}
		&	\int_{D^c} \frac{\delta_{D^c}(z)^{a_3}}{|x-z|^{d+a_1}|z-y|^{d+a_2}}dz \le   \frac{ C}{(|x-y| + \delta_D^\vee(x,y))^{d+a_1+a_2-a_3} }.	\end{align*}
\end{lemma}

\begin{lemma}\label{l:technical-2}
	For any  $a_1\in (0,d]$, $a_2>0$ and $a_3\in [0,a_1\wedge 1)$,  there exists $C=C(a_1,a_2,a_3)>0$ such that for all $y \in D$ and $z \in \overline D^c$,
	\begin{align*}
		\int_{\overline D^c} \frac{dw}{|z-w|^{d-a_1}|w-y|^{d+a_2} \delta_{\overline D^c}(w)^{a_3}} \le   \frac{C}{|y-z|^{d-a_1} \delta_D(y)^{a_2+a_3}}.	\end{align*}
\end{lemma}
\begin{proof}
	Let $y \in D$ and $z \in \overline D^c$. 
	Define
	\begin{align*}
		S_1&:= B(z,\delta_{\overline D^c}(z)/2),\nn\\
		S_2&:=\{w \in \overline D^c \setminus S_1: |w-y| \ge 2|z-w|\} \quad \text{and} \quad 
		S_3:=\overline D^c \setminus (S_1 \cup S_2).
	\end{align*}  We decompose the integral as  
	\begin{align*}
		I_1+I_2+I_3
		&:= \bigg( \int_{S_1} + \int_{S_2}+\int_{S_3}\bigg)\,  \frac{dw}{|z-w|^{d-a_1}|w-y|^{d+a_2}\delta_{\overline D^c}(w)^{a_3}}.
	\end{align*}
	
	For all $w\in  S_1$, by \eqref{e:dist-D-Dc}, we have $|w-y| \ge |y-z|-|z-w| \ge |y-z| - \delta_{\overline D^c}(z)/2 \ge |y-z|/2$ and $\delta_{\overline D^c}(w) \ge \delta_{\overline D^c}(z) /2$. Hence, we get
	\begin{align*}
		I_1 &	\le \frac{2^{d+a_2+a_3}}{|y-z|^{d+a_2}\delta_{\overline D^c}(z)^{a_3}} \int_{B(z, \delta_{\overline D^c}(z)/2)} \frac{dw}{|z-w|^{d-a_1}} = \frac{c_1 \delta_{\overline D^c}(z)^{a_1-a_3}}{|y-z|^{d+a_2}} \\
		&\le  \frac{c_1}{|y-z|^{d-a_1+a_2+a_3}} \le  \frac{c_1}{|y-z|^{d-a_1}\delta_D(y)^{a_2+a_3}},
	\end{align*}
	where we used \eqref{e:dist-D-Dc} in the second and third inequalities.
	
	For all $w \in S_2$, we have $|y-z|\le |w-y| + |z-w| \le (3/2)|w-y|$, and $2|z-w|\le |w-y|\le |y-z|+|z-w|$ so that $|z-w|\le |y-z|$. Thus, using Lemma \ref{l:Aikawa-integrals}(i) (with $D$ replaced by $\overline D^c$) and \eqref{e:dist-D-Dc}, we obtain
	\begin{align*}
		I_2	&\le  \frac{(3/2)^{d+a_2}}{|y-z|^{d+a_2}}\int_{B_{\overline D^c}(z,|y-z|)}  \frac{dw}{|z-w|^{d-a_1}\delta_{\overline D^c}(w)^{a_3}}\\
		&\le  \frac{c_2}{|y-z|^{2d-a_1+a_2}}\sum_{n=0}^\infty 2^{n(d-a_1)} \int_{B_{\overline D^c}(z,2^{-n}|y-z|) \setminus B(z, 2^{-n-1}|y-z|)}  \frac{dw}{\delta_{\overline D^c}(w)^{a_3}} \\
		&\le  \frac{c_3}{|y-z|^{d-a_1+a_2+a_3}}\sum_{n=0}^\infty 2^{-n(a_1-a_3)}  =  \frac{c_4}{|y-z|^{d-a_1+a_2+a_3}} \le \frac{c_4}{|y-z|^{d-a_1}\delta_D(y)^{a_2+a_3}}.
	\end{align*}
	
	To estimate $I_3$,   recall that  $Q_y \in \partial D$ is a point satisfying $\delta_D(y)=|y-Q_y|$. 
	Observe that for all $w\in S_3$, we have  $	|y-z|\le |w-y|+|z-w|\le 3|z-w|$ and $|w-y| \ge |w-Q_y| - \delta_D(y)$.  Hence, using \eqref{e:dist-D-Dc} and Lemma \ref{l:Aikawa-integrals}(i), we get that
	\begin{align*}
		I_3&\le \frac{3^{d-a_1}}{|y-z|^{d-a_1}}	\int_{\overline D^c}   \frac{dw}{|w-y|^{d+a_2}\delta_{\overline D^c}(w)^{a_3}}\\
		&\le \frac{3^{d-a_1}}{|y-z|^{d-a_1}\delta_D(y)^{d+a_2}} \int_{B_{\overline D^c}(Q_y, 2\delta_D(y))}   \frac{dw}{\delta_{\overline D^c}(w)^{a_3}}  \\
		&\quad + \frac{3^{d-a_1}}{|y-z|^{d-a_1}\delta_D(y)^{d+a_2}} \sum_{n=2}^\infty \frac{1}{(2^{n-1} -1)^{d+a_2}}	\int_{B_{\overline D^c}(Q_y, 2^n \delta_D(y)) \setminus B(Q_y, 2^{n-1}\delta_D(y))}   \frac{dw}{\delta_{\overline D^c}(w)^{a_3}} \\
		&\le \frac{c_5}{|y-z|^{d-a_1}\delta_D(y)^{a_2+a_3} } + \frac{c_6}{|y-z|^{d-a_1}\delta_D(y)^{a_2+a_3}} \sum_{n=2}^\infty 2^{-n(a_2+a_3)} = \frac{c_7}{|y-z|^{d-a_1}\delta_D(y)^{a_2+a_3} }.
	\end{align*}The proof is complete.
\end{proof}

 Recall that $E_\nu^Q(r)$ is defined in \eqref{e:cone-like}. For $Q\in \partial D$, $0<\nu<\eta$ and $r \in (0,  \wh R /4)$,  define
\begin{align*}
	\wt E^Q_{\nu}(r)&:=\left\{y=(\wt{y},y_d) \text{ in CS$_Q$}: |\wt{y}|<r/4,\, 4r^{-\nu}|\wt{y}|^{1+\nu}<-y_d<r/2\right\}.
\end{align*}
Note that $\wt E^Q_{\nu}(r)$ is defined as $E^Q_{\nu}(r)$ with $D$ replaced by $D^c$.

\begin{lemma}\label{l:distance-E'} For every $0<\nu<\nu'<\eta$, 
	there exists $C>0$ such that for all $r \in (0,  \wh R /4)$, $Q\in \partial D$,  $x\in E_\nu^Q(r)$ and $z \in \overline D^c \setminus  \wt E_{\nu'}^Q(r)$,
	\begin{align*}
		\delta_{\overline D^c}(z) \le C  r^{-(\nu'-\nu)/(1+\nu)}  |x-z|^{(1+\nu')/(1+\nu)}.
	\end{align*}
\end{lemma}
\begin{proof}
	Let $Q \in \partial D$ and  let $\Psi=\Psi^Q$ be the function in Definition \ref{df-c11}.  We work in the coordinate system CS$_Q$. Let $x=(\wt x, x_d)\in E_\nu^Q(r)$ and	$z=(\wt z, z_d)\in  \overline D^c \setminus  \wt E_{\nu'}^Q(r)$. We consider the following three cases separately.

	\textit{Case 1:} Suppose $|z|\ge r$. In this case, 
	 $|x-z| \ge |z|-|x| \ge |z|-(3/4)r \ge |z|/4$. Since $\delta_{\overline D^c}(z) \le |z|$, it follows that 	\begin{align*}	 r^{-(\nu'-\nu)/(1+\nu)}  |x-z|^{(1+\nu')/(1+\nu)} \ge  |z|^{-(\nu'-\nu)/(1+\nu)} (|z|/4)^{(1+\nu')/(1+\nu)}\ge 4^{-(1+\nu')/(1+\nu)}\delta_{\overline D^c}(z),	\end{align*} 	which yields the desired inequality.

\textit{Case 2:} 
	Suppose $|z|< r$, and either  $-z_d \ge r/2$ or $|\wt z| \ge r/4$. We claim that
	\begin{align}\label{e:distance-E'-claim}
		|x-z| \ge  (2^{-3} \wedge 2^{-2-3\nu})  r.
	\end{align}  Indeed, if $-z_d \ge r/2$, then $|x-z| \ge x_d-z_d \ge -z_d \ge r/2$, which implies \eqref{e:distance-E'-claim}. If $|\wt x-\wt z|\ge r/8$, then \eqref{e:distance-E'-claim} is evident.  Assume $-z_d<r/2$, $|\wt z| \ge r/4$ and $|\wt x-\wt z| < r/8$. Then $|\wt x| \ge |\wt z|/2 \ge r/8$. Since $z\in B(0,r) \setminus \overline D$, $\Lambda r^\eta < \Lambda (\wh R/4)^\eta\le 2^{-3\eta}$ and $\eta>\nu$, by \eqref{e:Psi-maximum}, we have 
	$$z_d< \Psi(\wt z)\le \Lambda |\wt z|^{1+\eta} \le 2^{-3\eta}|\wt z| \le  2^{1-3\eta} |\wt x|<2^{1-3\nu} |\wt x|. $$ On the other hand, since $x\in E_\nu^Q(r)$ and $\nu<\eta$, we get that $
	x_d \ge 4r^{-\nu} |\wt x|^{1+\nu} \ge 4r^{-\nu} (r/8)^\nu |\wt x| = 2^{2-3\nu} |\wt x|.$ It follows that $	|x-z| \ge x_d -z_d > 2^{1-3\nu} |\wt x| \ge 2^{-2-3\nu}  r$, which completes the proof of \eqref{e:distance-E'-claim}. 
Now, since $\delta_{\overline D^c}(z)\le |z|<r$ in this case, we get from  \eqref{e:distance-E'-claim} that
	\begin{align*}
	 r^{-(\nu'-\nu)/(1+\nu)}  |x-z|^{(1+\nu')/(1+\nu)} \ge c_1 r >c_1  \delta_{\overline D^c}(z).
	\end{align*}
	
	\textit{Case 3:} Suppose that  $|\wt z| < r/4$ and $-z_d < r/2$. Since $(\wt z, \Psi(\wt z)) \in \partial D$, repeating the argument for \eqref{e:distance-in-Enu-2} and using the definition of $\wt E_{\nu'}^Q(r)$, we obtain 
	\begin{align*}
		\delta_{\overline D^c}(z) \le \Psi(\wt z)  - z_d \le  \Lambda |\wt z|^{1+\eta} -z_d \le \Lambda |\wt z|^{1+\eta}  + 4r^{-\nu'}|\wt z|^{1+\nu'} \le (\Lambda+4)r^{-\nu'}|\wt z|^{1+\nu'}.
	\end{align*}
	We used the facts that $\nu'\in (0,\eta)$ and $|\wt z| < r<1$ in the last inequality. Thus,  using Lemma \ref{l:distance-in-Enu}(i) and \eqref{e:dist-D-Dc}, we obtain that 
	\begin{align*}
		\delta_{\overline D^c}(z) &\le  (\Lambda+4)r^{-\nu'} (|\wt x - \wt z| + |\wt x|)^{1+\nu'}       \le c_2 r^{-\nu'} (|\wt x - \wt z|^{1+\nu'} + |\wt x|^{1+\nu'}) \\
		&\le  c_2 r^{-\nu'}|x-z|^{1+\nu'} + c_3 r^{-\nu'} (r^\nu x_d)^{(1+\nu')/(1+\nu)} \\
		&\le c_2 r^{-\nu'}|x-z|^{1+\nu'} + c_4 r^{-(\nu'-\nu)/(1+\nu)}  \delta_D(x)^{(1+\nu')/(1+\nu)}\\
		&\le c_2 r^{-\nu'}|x-z|^{1+\nu'} + c_4 r^{-(\nu'-\nu)/(1+\nu)}  |x-z|^{(1+\nu')/(1+\nu)}.
	\end{align*}
	Since $|x-z| \le |x| + |z| <r$, this yields the desired result. 
\end{proof}

In the remainder of this section,  we  let $R:=\wh R/8$ and fix  $0<\nu<\nu'<\eta$.
For simplicity,  for $Q\in \partial D$ we  write
$
E^Q_\nu:=E^Q_\nu(R),
$
and 
$	\wt E^Q_{\nu'}:=\wt E^Q_{\nu'}(R).
$

Recall that  $n_Q$ denotes the inward unit normal to $\partial D$ at $Q$ 
and that 
\begin{align*}
		\bH^Q =
	\left\{ y \in \R^d:  n_Q\cdot (y-Q) > 0 \right\}.
\end{align*}
Note that, for all $Q\in \partial D$, since the function $\Psi$ in Definition \ref{df-c11} satisfies $\nabla \Psi(\wt 0) = 0$, we have $n_Q=(\wt 0,1)$ in CS$_Q$. Consequently,  $\delta_{(\overline\bH^Q)^c}(y)=(-y_d)_+$ for all $y =(\wt y, y_d) \in B(Q,\wh R)$ in CS$_Q$.

Recall that $A_{d,-\alpha}$ is defined in  \eqref{e:Ada}. Define 
 $\kappa_{1,\alpha}:= A_{1,-\alpha}/\alpha$, and for $d\ge 2$, 
\begin{align}\label{e:def-kappa-alpha}
\kappa_{d,\alpha}&:=  \frac{ A_{d,-\alpha} \omega_{d-1}}{2} B \bigg( \frac{\alpha+1}{2}, \frac{d-1}{2} \bigg) \int_0^1 \frac{(t^{\alpha/2}-1)(1- t^{\alpha/2-1})}{(1-t)^{1+\alpha}}dt\nn\\
&= \frac{A_{d,-\alpha}\omega_{d-1}}{2\alpha} B \bigg( \frac{\alpha+1}{2}, \frac{d-1}{2} \bigg),
\end{align}
where $\omega_{d-1}$ denotes the $(d-2)$-dimensional Lebesgue measure of the unit sphere in $\R^{d-1}$ and $B(\cdot,\cdot)$ is the beta function.  
Recall that $\wh \kappa_D(z)$ is defined in \eqref{e:def-wh-kappa}. By \cite[(5.6)]{BBC03}, we have 
\begin{align}\label{e:halfspace-killing}
\wh \kappa_\bH(z)  = \kappa_{d,\alpha} \delta_{\overline \bH^c}(z)^{-\alpha}, \quad z\in \overline \bH^c.
\end{align}
Moreover, applying \cite[Lemma 5.6]{CS25} to $\overline D^c$,  since $\overline D^c$ is a $C^{1,\eta \wedge(\alpha/2)}$ open set, we obtain the following lemma.

\begin{lemma}\label{l:C1eta-killing}
	There exists $C>0$ such that for all  $z\in\overline D^c$,
	\begin{equation*}
		\left|\wh \kappa_D(z)   - 
	\kappa_{d,\alpha} \delta_{\overline D^c}(z)^{-\alpha}  \right| \le   C\delta_{\overline D^c}(z)^{-\alpha + (\eta \wedge (\alpha/2))}.
	\end{equation*}
\end{lemma}

\subsubsection{Proof of Proposition \ref{p:nonlocal-Neumann}}

Since $D$ is a $C^{1,\eta}$ open set that is either bounded, or half-space-like, or the complement of a bounded set, it is easy to see that for all $z\in \overline D^c$,
\begin{align}\label{e:Neumann-integral}
\wh \kappa_D(z) = A_{d,-\alpha}	\int_D \frac{dw}{|z-w|^{d+\alpha}} \asymp \delta_{\overline D^c}(z)^{-\alpha} \left(1+ \frac{\delta_{\overline D^c}(z)}{\diam(D)} \right)^{-d} .
\end{align}
\begin{lemma}\label{l:Neumann-geometry}
	There exists $C>0$ such that for all $y \in D$ and $z\in \overline D^c$,
	\begin{align*}
		 \frac{ J^\alpha(z,y)}{\wh \kappa_D(z)} \le  \frac{C}{(\delta_D(y)\wedge 1)^{d}}\left( \frac{\delta_{\overline D^c}(z)}{\delta_D(y)} \wedge 1 \right)^{\alpha}.
	\end{align*}
\end{lemma}
\begin{proof}
	Let $y \in D$ and $z\in \overline D^c$. If $\delta_{\overline D^c}(z) <\diam(D)$, then using    \eqref{e:Neumann-integral} and \eqref{e:dist-D-Dc}, we obtain
	\begin{align*}
		\frac{ J^\alpha(z,y)}{\wh \kappa_D(z)} \le \frac{c_1\delta_{\overline D^c}(z)^{\alpha}}{|z-y|^{d+\alpha}}   \le  \frac{c_1}{\delta_D(y)^d} \left(\frac{\delta_{\overline D^c}(z)}{\delta_D(y)}  \wedge  1 \right) ^{\alpha}.
	\end{align*}
 If $\delta_{\overline D^c}(z) \ge \diam(D)$, then using    \eqref{e:Neumann-integral} and \eqref{e:dist-D-Dc}, we also get that
	\begin{align*}
		\frac{ J^\alpha(z,y)}{\wh \kappa_D(z)} \le \frac{c_2\delta_{\overline D^c}(z)^{d+\alpha}}{|z-y|^{d+\alpha}}   \le  c_2 \left( \frac{\delta_{\overline D^c}(z)}{\delta_D(y)} \wedge 1 \right)^{\alpha}.
	\end{align*}
\end{proof}

Next, we establish that $\sB_D^\alpha(x,y)$  satisfies  \hyperlink{B2}{{\bf (B2)}}.  Observe that $\sB_\bH^\alpha(x,y)$ is symmetric and  invariant under dilations and horizontal translations.  
Thus, by the discussion at the beginning of Section \ref{s:F}, there exists $ \F_\bH^\alpha:\bH_{-1}\to [0,\infty)$   such that
\begin{align}\label{e:def-Neumann-case}
	\sB_\bH^\alpha(x,y) = A_{d,-\alpha}\F_\bH^\alpha((y-x)/x_d) \quad \text{for all $x,y \in \bH$}.
\end{align}
We will show that $\sB_D^\alpha(x,y)$ satisfies  \hyperlink{B2}{{\bf (B2)}} with $\F_0=\F_\bH^\alpha$.  To estimate $| \sB_D^\alpha(x,y) - \sB_D^\alpha(x,x)\F_\bH^\alpha((y-x)/x_d) | = |\sB_D^\alpha(x,y) - \sB_\bH^\alpha(x,y)|$, we first establish the following  two lemmas.

\begin{lemma}\label{l:Neumann-1}
There exist constants $\eta_1\in (0,\eta]$ and  $C>0$ such that for all $Q\in \partial D$ and $z \in \wt E^Q_{\nu'}$,
	\begin{equation*}
\big| \wh \kappa_D(z)^{-1} -\wh \kappa_{\bH^Q}(z)^{-1}\big| \le   C\delta_{\overline D^c}(z)^{\alpha +\eta_1}.
	\end{equation*}
\end{lemma}
\begin{proof} Let $Q\in \partial D$ and $z \in \wt E^Q_{\nu'}$. We work in the coordinate system CS$_{Q}$, so that $\bH^Q=\bH$. By Lemma \ref{l:distance-in-Enu}(i) (with $D$ replaced by $\overline D^c$), 
	\begin{align}\label{e:Neumann-1-1}
		\frac45 \delta_{\overline D^c}(z) \le  \delta_{\overline {\bH}^c}(z) \le 2(1+2\Lambda)  \delta_{\overline D^c}(z).
	\end{align}
	Moreover, by the mean value theorem, \eqref{e:Neumann-1-1} and Lemma \ref{l:distance-in-Enu}(ii),  we have
	\begin{align}\label{e:Neumann-1-2}
		\left| \delta_{\overline D^c}(z)^{-\alpha} - \delta_{\overline \bH^c}(z)^{-\alpha} \right| \le c_1 \delta_{\overline D^c}(z)^{-\alpha-1} \left| \delta_{\overline D^c}(z)- \delta_{\overline \bH^c}(z)\right| \le c_2 \delta_{\overline D^c}(z)^{-\alpha+ (\eta-\nu')/(1+\nu')} .
	\end{align}
	
	Set $\eta_1:=\eta \wedge (\alpha/2) \wedge ( (\eta-\nu')/(1+\nu'))$.
	Using Lemma \ref{l:C1eta-killing}, \eqref{e:halfspace-killing}
	and \eqref{e:Neumann-1-2},
	since $\delta_{\overline D^c}(z)<1$, we obtain that
\begin{align}\label{e:Neumann-1-3}
		&\big	| \wh \kappa_D(z) -\wh \kappa_{\bH}(z) \big| \nn\\
		&\le 	\big| 	 \wh \kappa_D(z)-  \kappa_{d,\alpha} \delta_{\overline D^c}(z)^{-\alpha}\big|   + \big| 	\wh \kappa_{\bH}(z) -  \kappa_{d,\alpha} \delta_{\overline \bH^c}(z)^{-\alpha}\big| + \kappa_{d,\alpha} \big| \delta_{\overline D^c}(z)^{-\alpha} - \delta_{\overline \bH^c}(z)^{-\alpha} \big| \nn\\
		&\le c_3\delta_{\overline D^c}(z)^{-\alpha + (\eta \wedge (\alpha/2))}
		+c_4 \delta_{\overline D^c}(z)^{-\alpha+ (\eta-\nu')/(1+\nu')} \le c_5\delta_{\overline D^c}(z)^{-\alpha+\eta_1}.
\end{align}	
	On the other hand,  
	 by \eqref{e:Neumann-integral},  \eqref{e:Neumann-1-1} and \eqref{e:halfspace-killing},  since $\delta_{\overline D^c}(z)<1$,  
	\begin{align}\label{e:Neumann-1-4}	\wh \kappa_D(z)  \asymp \kappa_{d,\alpha}  \delta_{\overline D^c}(z)^{-\alpha} \asymp 	 \kappa_{d,\alpha}  \delta_{\overline {\bH}^c}(z)^{-\alpha} = \wh \kappa_\bH(z) .
	\end{align}
	Combining \eqref{e:Neumann-1-3} with \eqref{e:Neumann-1-4}, we arrive at
\begin{align*}
\big| \wh \kappa_D(z)^{-1} -\wh \kappa_{\bH}(z)^{-1}\big| = 	\frac{| \wh \kappa_D(z) - \wh \kappa_\bH(z)|}{ \wh \kappa_D(z) \wh \kappa_\bH(z) } \le  c_6 \delta_{\overline D^c}(z)^{\alpha+\eta_1}.
\end{align*}
\end{proof}

\begin{lemma}\label{l:Neumann-2}
	There exist constants $\eta_2\in (0,\alpha/2)$ and  $C>0$  such that for all $Q\in \partial D$ and $x,y \in E_\nu^Q$,
	\begin{align*}
	\int_{\overline D^c\setminus \wt E_{\nu'}^Q} \frac{J^\alpha(x,z)J^\alpha(z,y)}{  \wh \kappa_D(z) } \, dz  \le \frac{CR^{-\eta_2}}{ \delta_D^\wedge(x,y)^{\alpha/2-\eta_2} \delta_D^\vee(x,y)^{d+\alpha/2}} \, .
	\end{align*}
\end{lemma}
\begin{proof}
	Let $Q\in \partial D$ and $x,y \in E_\nu^Q$. By symmetry, without loss of generality, we assume that $\delta_D(x)\le \delta_D(y)$.  Set $\eta_2:=(\alpha/2)(\nu'-\nu)/(1+\nu)$.    Using Lemmas  \ref{l:Neumann-geometry} and \ref{l:distance-E'},  we obtain
	\begin{align*}
&	\int_{\overline D^c\setminus \wt E_{\nu'}^Q} \frac{J^\alpha(x,z)J^\alpha(z,y)}{  \wh \kappa_D(z) } \, dz \le \frac{c_1}{\delta_D(y)^d} 	\int_{\overline D^c\setminus \wt E_{\nu'}^Q} \left( \frac{\delta_{\overline D^c}(z)}{\delta_D(y)} \wedge 1 \right)^{\alpha/2}\frac{1 }{|x-z|^{d+\alpha}} dz\nn\\
	&\le  \frac{c_2R^{-\eta_2}}{\delta_D(y)^{d+\alpha/2}} 	\int_{\overline D^c \setminus \wt E_{\nu'}^Q}\frac{ dz }{|x-z|^{d+\alpha - (\alpha/2 + \eta_2)}} \\
		&\le  \frac{c_2R^{-
		 \eta_2 
		} }{\delta_D(y)^{d+\alpha/2}} 	\int_{B(x,\delta_D(x))^c}\frac{dz }{|x-z|^{d+\alpha/2-
		 \eta_2 
		}} =  \frac{c_3R^{- 
		\eta_2 
		} }{\delta_D(x)^{\alpha/2-
		\eta_2 
		}\delta_D(y)^{d+\alpha/2}}.
	\end{align*}	
\end{proof}

\begin{proof}[Proof of Proposition \ref{p:nonlocal-Neumann}]
By \cite[Lemma 7.2]{KSV23},  
$\sB_D^\alpha(x,y)$ satisfies \eqref{e:ass-B1} with $\theta_0=d+\alpha$ so  \hyperlink{B1}{{\bf (B1)}} holds true. 
We now prove that $\sB_D^\alpha(x,y)$ satisfies \hyperlink{B2}{{\bf (B2)}} with $\theta_1=d+\alpha$ and $\theta_2= \eta_1\wedge \eta_2$, where $\eta_1 \in (0,\eta]$ and $\eta_2\in (0,\alpha/2)$ are the constants in Lemmas  \ref{l:Neumann-1} and \ref{l:Neumann-2}, respectively. 	

 Define $\F_\bH^\alpha:\bH_{-1}\to (0,\infty)$ by \eqref{e:def-Neumann-case}. Let $Q\in \partial D$ and $x,y \in E_\nu^Q$. We work in the coordinate system CS$_Q$, so that $\bH^Q=\bH$.   By \eqref{e:def-Neumann-case}, 
	\begin{align*}
		&\left| \sB_D^\alpha(x,y) - A_{d,-\alpha}
		\F_\bH^\alpha
		 ((y-x)/x_d) \right|  =	\left| \sB_D^\alpha(x,y) - \sB_\bH^\alpha(x,y) \right|\\
		&= |x-y|^{d+\alpha} \,\bigg| \int_{\overline D^c} \frac{J^\alpha(x,z)J^\alpha(z,y)}{ \wh \kappa_D(z) } \, dz -\int_{\overline \bH^c} \frac{J^\alpha(x,z)J^\alpha(z,y)}{ \wh \kappa_\bH(z) } \, dz  \bigg|\\
		&\le |x-y|^{d+\alpha}  \int_{\wt E_{\nu'}^Q} J^\alpha(x,z)J^\alpha(z,y)  \left| \frac{1}{\wh \kappa_D(z) }- \frac{1}{\wh \kappa_\bH(z) } \right| dz    \\
		&\quad +  |x-y|^{d+\alpha}\int_{\overline D^c\setminus \wt E_{\nu'}^Q} \frac{J^\alpha(x,z)J^\alpha(z,y)}{\wh \kappa_D(z)  } \, dz +  |x-y|^{d+\alpha}\int_{\overline \bH^c\setminus \wt E_{\nu'}^Q} \frac{J^\alpha(x,z)J^\alpha(z,y)}{ \wh \kappa_\bH(z)  } \, dz\\
		&=:I_1+I_2+I_3.
	\end{align*}
For $I_1$,	using Lemmas \ref{l:Neumann-1} and \ref{l:technical-1}, since $\eta_2<\alpha/2$, we get that
	\begin{align*}
		I_1& \le c_1|x-y|^{d+\alpha}   \int_{\wt E_{\nu'}^Q} \frac{\delta_{\overline D^c}(z)^{\alpha+(\eta_1 \wedge (\alpha/2))}}{|x-z|^{d+\alpha}|z-y|^{d+\alpha}} dz \le c_1|x-y|^{d+\alpha}   \int_{\overline D^c} \frac{\delta_{\overline D^c}(z)^{\alpha+(\eta_1 \wedge (\alpha/2))}}{|x-z|^{d+\alpha}|z-y|^{d+\alpha}} dz\\
		&\le  \frac{ c_2|x-y|^{d+\alpha}}{(|x-y| + \delta_D^\vee(x,y))^{d+\alpha - (\eta_1 \wedge (\alpha/2))} } \le c_2\left( \frac{|x-y| \vee \delta_D^\vee (x,y)}{|x-y| \wedge \delta_D^\wedge (x,y)} \right)^{d+\alpha} \left(|x-y| \vee \delta_D^\vee (x,y) \right)^{\eta_1 \wedge \eta_2}.
	\end{align*}
For $I_2+I_3$, applying Lemma \ref{l:Neumann-2} twice,  to $D$ and $\bH$, respectively, and using Lemma \ref{l:distance-in-Enu}(i), we obtain that
\begin{align*}
	I_2+I_3&\le \frac{c_3|x-y|^{d+\alpha}}{ \delta_D^\wedge(x,y)^{\alpha/2} \delta_D^\vee(x,y)^{d+\alpha/2}} \left( \frac{\delta_D^\wedge (x,y)}{R}\right)^{\eta_2} +  \frac{c_3|x-y|^{d+\alpha}}{ (x_d\wedge y_d)^{\alpha/2} (x_d\vee y_d)^{d+\alpha/2}} \left( \frac{	x_d\wedge y_d
	}{R}\right)^{\eta_2}\\
	&\le \frac{c_4|x-y|^{d+\alpha}}{ \delta_D^\wedge(x,y)^{\alpha/2} \delta_D^\vee(x,y)^{d+\alpha/2}} \left( \frac{\delta_D^\wedge (x,y)}{R}\right)^{\eta_2}\le \frac{c_4|x-y|^{d+\alpha}}{ \delta_D^\wedge(x,y)^{\alpha/2} \delta_D^\vee(x,y)^{d+\alpha/2}} \left( \frac{\delta_D^\wedge (x,y)}{R}\right)^{\eta_1\wedge \eta_2} \\
	&\le c_4\left( \frac{|x-y| \vee \delta_D^\vee (x,y)}{|x-y| \wedge \delta_D^\wedge (x,y)} \right)^{d+\alpha} \left( \frac{|x-y| \vee \delta_D^\vee (x,y)}{R}\right)^{\eta_1\wedge \eta_2}.
\end{align*}
We used $\delta_D^\wedge (x,y)<R$ in the third inequality.  The proof is complete.
\end{proof}

\subsubsection{Proof of Proposition \ref{p:nonlocal-trace}}
Throughout this subsection, we assume that 
$d\ge 2$ and $\alpha\in (1, 2)$. 
For an open set $U\subset \R^d$, 
let $g^\alpha_{U} (x,y)$ denote the Green function of the isotropic $\alpha$-stable process on $U$. It is well known that $g^\alpha_{\R^d}(x,y) =A_{d,\alpha} |x-y|^{\alpha-d}$, $x,y \in \R^d$.

\begin{lemma}\label{l:trace-pre}
There exists $C>0$ such that for all $y \in D$ and $z\in \overline D^c$,
	\begin{align*}
		\int_{\overline D^c}  g^\alpha_{\overline D^c} (z,w) J^\alpha(w,y) \,dw \le \frac{C}{|y-z|^{d-\alpha} \delta_D(y)^\alpha}.
	\end{align*}
\end{lemma}
\begin{proof}
	For all
	$y \in D$ and $z\in \overline D^c$, using  the domain monotoncity and Lemma \ref{l:technical-2}, we obtain 
	\begin{align*}
&			\int_{\overline D^c}  g^\alpha_{\overline D^c} (z,w) J^\alpha(w,y) \,dw \le 	\int_{\overline D^c}  g^\alpha_{\R^d} (z,w) J^\alpha(w,y) \,dw \\
&= c_1 \int_{\overline D^c}  \frac{dw}{|z-w|^{d-\alpha}|w-y|^{d+\alpha}} \le \frac{c_3}{|y-z|^{d-\alpha} \delta_D(y)^\alpha}.
	\end{align*}
\end{proof}

Since $\wc\sB_\bH^\alpha(x,y)$ is symmetric and  invariant under dilations and horizontal translations, by the discussion at the beginning of Section \ref{s:F}, there exists $ \wc\F_\bH^\alpha:\bH_{-1}\to [0,\infty)$   such that
\begin{align}\label{e:def-trace-case}
\wc	\sB_\bH^\alpha(x,y) = A_{d,-\alpha}\wc\F_\bH^\alpha((y-x)/x_d) \quad \text{for all $x,y \in \bH$}.
\end{align}

We next show that $\wc\sB_D^\alpha(x,y)$ satisfies  \hyperlink{B2}{{\bf (B2)}} with $\F_0=\wc\F_\bH^\alpha$. To this end, recall that $Y=(Y_t)_{t\ge 0}$ denotes the isotropic $\alpha$-stable process in $\R^d$ with  jump kernel $J^\alpha(x,y)$. The 
function $\sB$ 
associated with $J^\alpha(x,y)$ is the constant $A_{d,-\alpha}$, which  clearly satisfies  \hyperlink{B1}{{\bf (B1)}} and  \hyperlink{B2}{{\bf (B2)}} with $\F_0=1$. Moreover, by Lemma \ref{l:C1eta-killing}, $Y^{\overline D^c}$ satisfies \hyperlink{K}{{\bf (K)}}, with
the associated parameter $p$ in \eqref{e:C(alpha,p,F)} given by $p=\alpha/2$. 
  Consequently, the estimate in Proposition \ref{p:prelub} applies to $g^\alpha_{B_{\overline D^c} (x_0,R)}$ with $p=\alpha/2$ for every $x_0 \in D^c$ and $R>0$.

For an open set $U\subset \R^d$, define $\tau^\alpha_U:=\inf\{t>0: Y_t\notin U\}$.
\begin{lemma}\label{l:trace-green}
	There exist constants $\eta_3\in (0,d-\alpha]$ and  $C>0$ such that for all $Q\in \partial D$ and  $z,w \in \wt E^Q_{\nu'}$,
	\begin{align*}
		\left| g^\alpha_{\overline D^c} (z,w)  -  g^\alpha_{(\overline \bH^Q)^c} (z,w)  \right| \le \frac{CR^{-\eta_3}}{ \delta_{\overline D^c}^\vee (z,w)^{d-\alpha-\eta_3}}.
	\end{align*}
\end{lemma} 
\begin{proof}
	Let $Q\in \partial D$ and $z,w\in \wt E^Q_{\nu'}$. By symmetry, without loss of generality, we assume that $\delta_{\overline D^c}(z)\le \delta_{\overline D^c}(w)$. 
	We work in the coordinate system CS$_Q$, so that $\bH^Q=\bH$.  Set 
	$$\nu'':=\frac{\eta+\nu'}{2} \quad \text{and} \quad \eta_3:=(d-\alpha) \wedge \left( \frac{\alpha(\nu''-\nu')}{2(1+\nu')}\right). $$
	Let $U$ be either $\overline D^c$ or $\overline \bH^c$, and let $V$ denote the other set. 
	 By  the domain  monotonicity  and the Dynkin-Hunt formula,  
	\begin{align*}
&		g^\alpha_{U} (z,w)- 	g^\alpha_{V} (z,w) \le 	g^\alpha_{U} (z,w) -  g^\alpha_{U \cap V} (z,w)  = \E_z \Big[ g^\alpha_{U} (Y_{\tau^\alpha_{U \cap V}},w); \, Y_{\tau^\alpha_{U \cap V}} \in  U \setminus V \Big] \nn\\
		&= \E_z \Big[ g^\alpha_{U} (Y_{\tau^\alpha_{U \cap V}},w); \, Y_{\tau^\alpha_{U \cap V }} \in B_U(0,R) \setminus V \Big]+ \E_z \Big[ g^\alpha_{U} (Y_{\tau^\alpha_{U \cap V}},w); \, Y_{\tau^\alpha_{U \cap V }} \in U \setminus B(0,R) \Big]\nn\\
		&=: I_1+I_2.
	\end{align*}
	For all $y \in U\setminus B(0,R)$, since $|y-w| \ge |y|-|w| \ge R/4$, using the domain monotonicity,  we get that $g^\alpha_{U}  (y,w) \le g^\alpha_{\R^d}(y,w) \le c_2 R^{\alpha-d}$. Since $\delta_{\overline D^c}(w) <R$ and $\eta_3\le d-\alpha$, it follows that
	\begin{align} \label{e:censored-green-1}
	I_2 \le c_2R^{\alpha-d} \le  c_2R^{-\eta_3} \delta_{\overline D^c}(w)^{\alpha+\eta_3-d}.
	\end{align} 
	Next, we estimate $I_1$. Let $x \in B_U(0,R)\setminus V$. By the Dynkin-Hunt formula, 
	\begin{align*}
		g^\alpha_U(x,w) &=  g^\alpha_{B_U(0,2R)} (x,w)  + \E_x \Big[ g^\alpha_{U} (Y_{\tau^\alpha_{B_U(0,2R)}},w); \, Y_{\tau^\alpha_{B_U(0,2R)}} \in  U \setminus B(0,2R) \Big] \\
		&=:I_{1,1}+I_{1,2}.
	\end{align*}
	Repeating the argument leading to \eqref{e:censored-green-1}, since $|y-w| \ge |y|-|w| \ge 7R/4$ for all $y \in U \setminus B(0,2R)$, 
	we obtain $I_{1,2} \le c_3 R^{\alpha-d}\le c_3R^{-\eta_3} \delta_{\overline D^c}(w)^{\alpha+\eta_3-d}$. For $I_{1,1}$,  note that $x\in U \setminus V \subset U \setminus  E_{\nu''}^Q(R)$.  
	Using Proposition \ref{p:prelub} and 
	Lemma \ref{l:distance-E'} 
	(with $D$ replaced by $\overline D^c$ or $\overline \bH^c$),  we obtain 
	\begin{align*}
		I_{1,1} &\le c_4 \left( \frac{\delta_U(x)}{|x-w|}\right)^{\alpha/2} |x-w|^{\alpha-d} \le c_5 R^{-(\alpha/2)(\nu''-\nu')/(1+\nu')} |x-w|^{\alpha-d+(\alpha/2)(\nu''-\nu')/(1+\nu')}\\
	&	\le c_5 R^{-\eta_3} |x-w|^{\alpha+\eta_3-d} \le c_5 R^{-\eta_3} \delta_V(w)^{\alpha+\eta_3-d}.
	\end{align*}
	We used \eqref{e:dist-D-Dc} and the fact that $\alpha+\eta_3-d\le 0$ in the last inequality. 
	Since $\delta_V(w) \asymp \delta_{\overline D^c}(w)$ if $V=\overline \bH^c$ by Lemma \ref{l:distance-in-Enu}(i), we conclude that
	\begin{align*}
		I_1 \le  c_6R^{-\eta_3} \delta_{\overline D^c}(w)^{\alpha+\eta_3-d}  \P_z \big(  Y_{\tau^\alpha_{U \cap V }} \in B_U(0,R) \setminus V \big) \le c_6R^{-\eta_3} \delta_{\overline D^c}(w)^{\alpha+\eta_3-d}.
	\end{align*}
	The proof is complete.
\end{proof}

\begin{lemma}\label{l:trace-B2-1}
	There exist constants   $\eta_4>0$, 
	$\eta_5\in (0,d-\alpha]$
	 and $C>0$ such that for all $Q\in \partial D$, $y\in E_\nu^Q$ and $z \in \wt E_{\nu'}^Q$,
	\begin{equation*}
	\int_{\wt E_{\nu'}^Q} 	\Big|  g^\alpha_{\overline D^c} (z,w)  -g^\alpha_{(\overline \bH^Q)^c}  (z,w)  \Big|  J^\alpha(w,y) dw \le \frac{CR^{-\eta_4}}{|y-z|^{d-\alpha-\eta_5} \delta_D(y)^{\alpha-\eta_4+\eta_5}}.
	\end{equation*}
\end{lemma}
\begin{proof}
Let	$Q\in \partial D$, $y\in E_\nu^Q$ and $z \in \wt E_{\nu'}^Q$.   
By the domain monotonicity and  Lemma \ref{l:trace-green},   we see that  for all $w \in \wt E^Q_{\nu'}$,
\begin{align*}
	&	\Big|  g^\alpha_{\overline D^c} (z,w)  -g^\alpha_{(\overline \bH^Q)^c}  (z,w)  \Big|   \le 	\Big(  g^\alpha_{\R^d} (z,w)   + g^\alpha_{\R^d}  (z,w)  \Big)^{1-1/d}	\Big|  g^\alpha_{\overline D^c} (z,w)  -g^\alpha_{(\overline \bH^Q)^c}  (z,w)  \Big| ^{1/d}\\
	&\le  c_1|z-w|^{-d +\alpha +1 -\alpha/d } R^{-\eta_3/d} \delta_{\overline D^c}(w)^{-(1-\alpha/d-\eta_3/d)}.
\end{align*}
Note that $\alpha +1 -\alpha/d  = (d-1)\alpha/d + 1<d$. 
 Hence, using Lemma \ref{l:technical-2}, we obtain that
\begin{align*}
&	\int_{\wt E_{\nu'}^Q} 	\Big|  g^\alpha_{\overline D^c} (z,w)  -g^\alpha_{(\overline \bH^Q)^c}  (z,w)  \Big|  J^\alpha(w,y) dw\\
&\le  c_2R^{-\eta_3/d} \int_{\overline D^c} \frac{ dw }{|z-w|^{d - \alpha- 1 + \alpha/d }\,|w-y|^{d+\alpha}\delta_{\overline D^c}(w)^{1 -(\alpha+\eta_3)/d}}\\
&\le  \frac{c_3R^{-\eta_3/d}}{ |y-z|^{d-\alpha -1 + \alpha/d}\delta_D(y)^{\alpha+1-(\alpha+\eta_3)/d}}.
\end{align*}
This proves the lemma with $\eta_4= \eta_3/d$ and $\eta_5 = 1-\alpha/d$. 
\end{proof}

\begin{lemma}\label{l:trace-B2-2}
There exist constants $\eta_6\in (0,\alpha/2)$ and  $C>0$  such that for all $Q\in \partial D$ and   $x,y \in E_{\nu}^Q$,
	\begin{align*}     
	 \int_{(\overline D^c\times \overline D^c) \setminus (\wt E_{\nu'}^Q \times \wt E_{\nu'}^Q)} 	\int_{ \overline D^c}   J^\alpha(x,z)  g^\alpha_{\overline D^c}(z,w ) J^\alpha(w,y) \,  dwdz  \le 
	 \frac{CR^{-(\alpha/2-\eta_6)}}{ (\delta_D(x)\wedge \delta_D(y))^{d+\alpha/2+\eta_6}}. 	\end{align*}
\end{lemma}
\begin{proof}
	Let $Q\in \partial D$ and   $x,y \in E_{\nu}^Q$.
	Set $\eta_6:=(\alpha/2)(1+\nu)/(1+\nu')$. Note that $\eta_6<\alpha/2<1$.  By Lemma \ref{l:distance-E'}, we have  $|x-z|^{\alpha/2} \ge c_1 R^{\alpha/2-\eta_6} \delta_{\overline D^c}(z)^{\eta_6}$  for all $z\in \overline D^c \setminus \wt E_{\nu'}^Q$.
Hence,	using Lemma \ref{l:trace-pre}, \eqref{e:dist-D-Dc} and Lemma \ref{l:technical-2} (with $a_1=d$, $a_2=\alpha/2$ and $a_3=\eta_6$), we obtain
	\begin{align*}
			&
				\int_{\overline D^c\setminus \wt E_{\nu'}^Q}     	\int_{ \overline D^c}   J^\alpha(x,z)  g^\alpha_{\overline D^c}(z,w ) J^\alpha(w,y) \, dwdz\le  \frac{c_2}{\delta_D(y)^\alpha} \int_{\overline D^c \setminus \wt E_{\nu'}^Q} \frac{dz}{ |y-z|^{d-\alpha}|x-z|^{d+\alpha}} \\
				&\le  \frac{c_3R^{-(\alpha/2-\eta_6)}}{\delta_D(y)^d} \int_{\overline D^c \setminus \wt E_{\nu'}^Q} \frac{dz}{|x-z|^{d+\alpha/2} \delta_{\overline D^c}(z)^{\eta_6}} \le \frac{c_4R^{-(\alpha/2-\eta_6)}}{\delta_D(x)^{\alpha/2+\eta_6}\delta_D(y)^d} \le  \frac{c_4R^{-(\alpha/2-\eta_6)}}{(\delta_D(x) \wedge \delta_D(y))^{d+\alpha/2+\eta_6}}.
	\end{align*}
	Similarly, we obtain $		\int_{ \overline D^c} \int_{\overline D^c\setminus \wt E_{\nu'}^Q}       J^\alpha(x,z)  g^\alpha_{\overline D^c}(z,w ) J^\alpha(w,y) \, dwdz \le \frac{c_4R^{-(\alpha/2-\eta_6)}}{(\delta_D(x) \wedge \delta_D(y))^{d+\alpha/2+\eta_6}}$. Combining these yields the desired result.
\end{proof}

\begin{proof}[Proof of Proposition \ref{p:nonlocal-trace}]
By \cite[Lemma 7.2]{KSV23},  
$\wc\sB_D^\alpha(x,y)$ satisfies \eqref{e:ass-B1} with $\theta_0=d+\alpha$ so  \hyperlink{B1}{{\bf (B1)}} holds true. 
Thus, 
it suffices to prove  that $\wc\sB_D^\alpha(x,y)$ satisfies \hyperlink{B2}{{\bf (B2)}} with $\theta_1=d+\alpha$ and $\theta_2:= \eta_4 \wedge (\alpha/2-\eta_6)$, where $\eta_4>0$ and $\eta_6\in (0,\alpha/2)$ are the constants in Lemmas  \ref{l:trace-B2-1} and \ref{l:trace-B2-2}, respectively. 
		Define $\wc\F_\bH^\alpha:\bH_{-1}\to (0,\infty)$ by \eqref{e:def-trace-case}. Let $Q\in \partial D$ and $x,y \in E_\nu^Q$. We work in the coordinate system CS$_Q$, so that $\bH^Q=\bH$.   
		By \eqref{e:def-trace-case} and the symmetry of $g^\alpha_{\overline D^c}$ and $g^\alpha_{\overline \bH^c}$, 
	\begin{align*}
		&\left|\wc \sB_D^\alpha(x,y) - A_{d,-\alpha}\wc\F^\alpha_\bH ((y-x)/x_d) \right|  =	\left| \wc\sB_D^\alpha(x,y) - \wc\sB_\bH^\alpha(x,y) \right|\\
		&= |x-y|^{d+\alpha} \,\bigg| \int_{\overline D^c \times \overline D^c} J^\alpha(x,z)g^\alpha_{\overline D^c}(z,w)J^\alpha(w,y) \, dzdw -\int_{\overline \bH^c\times \overline \bH^c} J^\alpha(x,z)g^\alpha_{\overline \bH^c}(z,w) J^\alpha(w,y) \, dzdw  \bigg|\\
		&\le |x-y|^{d+\alpha}  \int_{\wt E_{\nu'}^Q\times \wt E_{\nu'}^Q} J^\alpha(x,z) \left|  g^\alpha_{\overline D^c}(z,w) - g^\alpha_{\overline \bH^c}(z,w) \right|J^\alpha(w,y) \, dzdw     \\
			&\quad +  |x-y|^{d+\alpha}	\int_{(\overline D^c\times \overline D^c) \setminus (\wt E_{\nu'}^Q \times \wt E_{\nu'}^Q)}   J^\alpha(x,z)g^\alpha_{\overline D^c}(z,w)J^\alpha(w,y) \, dzdw\\
		&\quad  +  |x-y|^{d+\alpha}	\int_{(\overline \bH^c\times \overline \bH^c) \setminus (\wt E_{\nu'}^Q \times \wt E_{\nu'}^Q)}    J^\alpha(x,z)g^\alpha_{\overline \bH^c}(z,w)J^\alpha(w,y) \, dzdw\\
		&=:I_1+I_2+I_3.
	\end{align*}
	For $I_1$,	using Lemma \ref{l:trace-B2-1} and \eqref{e:dist-D-Dc}, we get that
	\begin{align*}
		I_1& \le \frac{c_1R^{-\eta_4}|x-y|^{d+\alpha}  }{\delta_D(y)^{\alpha-\eta_4+\eta_5}} \int_{\wt E_{\nu'}^Q} \frac{dz}{|x-z|^{d+\alpha}|y-z|^{d-\alpha-\eta_5}} \le \frac{c_1R^{-\eta_4}|x-y|^{d+\alpha}  }{\delta_D(y)^{d-\eta_4}} \int_{B(x,\delta_D(x))^c} \frac{dz}{|x-z|^{d+\alpha}} \\
		&= \frac{c_2R^{-\eta_4}|x-y|^{d+\alpha}  }{\delta_D(x)^\alpha \delta_D(y)^{d-\eta_4}} \le c_2\left( \frac{|x-y| \vee \delta_D^\vee (x,y)}{|x-y| \wedge \delta_D^\wedge (x,y)} \right)^{d+\alpha} \left( \frac{|x-y| \vee \delta_D^\vee (x,y)}{R}\right)^{\eta_4}.
	\end{align*}
	For $I_2+I_3$, applying Lemma \ref{l:trace-B2-2} twice,  to $D$ and $\bH$, respectively, and using Lemma \ref{l:distance-in-Enu}(i), we obtain that
	\begin{align*}
		I_2+I_3&\le  \frac{c_3|x-y|^{d+\alpha}}{ \delta_D^\wedge(x,y)^{d
				+\alpha}} \bigg( \frac{\delta_D^\wedge (x,y)}{R}\bigg)^{\alpha/2-\eta_6} +  \frac{c_3|x-y|^{d+\alpha}}{ (x_d\wedge y_d)^{d+\alpha} } \bigg( \frac{	x_d\wedge y_d}{R}\bigg)^{\alpha/2-\eta_6}\\
		&\le  \frac{c_4|x-y|^{d+\alpha}}{ \delta_D^\wedge(x,y)^{d+\alpha}} \left( \frac{\delta_D^\wedge (x,y)}{R}\right)^{\alpha/2-\eta_6} \le \frac{c_4|x-y|^{d+\alpha}}{ \delta_D^\wedge(x,y)^{d+\alpha}} \left( \frac{\delta_D^\wedge (x,y)}{R}\right)^{\theta_2} \\
		&\le c_4\left( \frac{|x-y| \vee \delta_D^\vee (x,y)}{|x-y| \wedge \delta_D^\wedge (x,y)} \right)^{d+\alpha} \left( \frac{|x-y| \vee \delta_D^\vee (x,y)}{R}\right)^{\theta_2}.
	\end{align*}
	We used $\delta_D^\wedge (x,y)<R$ in the third inequality.  The proof is complete.
\end{proof}

\subsection{A general $\sB$ satisfying {\bf (B2)}}

Suppose that $\Phi$ is a blow-up weight function satisfying \eqref{e:beta-range}. Suppose further that 
$\Phi$ satisfies the additional assumption that $\Phi$ is differentiable with 
\begin{align}\label{e:additional-Phi}
\sup_{t>0}\,t\Phi'(t)/\Phi(t)=:c_1<\infty \, . 
\end{align}
In this subsection we will show that, 
 for any $A_0\ge \wh R/4$,   the function 
$$
\sB^D(x,y):=\Phi\left(\frac{(|x-y| \wedge A_0 )^2}{(\delta_D(x) \wedge A_0)(\delta_D(y) \wedge A_0)}\right), \quad x,y \in D,
$$
satisfies \hyperlink{B2}{{\bf (B2)}} for all $
\nu_0 
\in (0,\eta)$ with 
$$\F_0(z):=\Phi\left(\frac{|z|^2 }{z_d+1}\right), \quad z \in \bH_{-1}.$$
The choice of this $\F_0$ can be viewed as $\F_0(x):=\sB^{\bH}(\e_d,x+ \e_d)$
where $\sB^\bH$ is defined by 
$ \sB^{\bH}(x,y):=\Phi({|x-y| ^2}/{(x_dy_d) })$. 
The function $\sB^{\bH}(x,y)$ 
is invariant under dilations and horizontal translations in $\bH$; 
see the discussion at the beginning of Section \ref{s:F}.
This is analogous to \cite[Subsection 11.2]{CKSV23}.

Let
 $\nu_0 \in (0,\eta)$  and  $Q\in \partial D$. We work in the coordinate system CS$_{Q}$.  
Set $R:=\wh R/8$ and  let $x,y\in E^Q_{\nu_0}(R)$. 
Without loss of generality, we assume that $\delta_D(x)\le \delta_D(y)$.  
 Observe that $\sB^D(x,x) =\Phi(0)=1$ and $\delta_D(x) \vee \delta_D(y) \vee |x-y| \le 2R\le A_0$.  Let 
$I:=[(\delta_D(x)\delta_D(y))\wedge (x_dy_d), (\delta_D(x)\delta_D(y))\vee (x_dy_d)]$.
Since $\F_0((y-x)/x_d)=\Phi\left(\frac{|x-y|^2}{x_dy_d}\right)$, 
 by the mean value theorem and  \eqref{e:additional-Phi},  we obtain that
\begin{align}\label{e:difference-1}
&\big|\sB^D(x,y)-\sB^D(x,x)\F_0((y-x)/x_d)\big|+\big|\sB^D(x,y)-\sB^D(y,y)\F_0((y-x)/x_d)\big|\nn \\
&\quad=
 2\left|\Phi\left(\frac{|x-y|^2}{\delta_D(x)\delta_D(y)}\right) 
-\Phi\left(\frac{|x-y|^2}{x_dy_d}\right) \right|    \nn \\
&\quad\le \frac{2|x-y|^2 |x_dy_d - \delta_D(x)\delta_D(y)|}{\delta_D(x)\delta_D(y)x_dy_d}\sup_{s\in I}\Phi'\bigg(\frac{|x-y|^2}{s}\bigg) \nn \\
&\quad\le  \frac{c_1 |x_dy_d - \delta_D(x)\delta_D(y)|}{\delta_D(x)\delta_D(y)x_dy_d} \sup_{s\in I}s\Phi\bigg(\frac{|x-y|^2}{s}\bigg) .
\end{align}
By Lemma \ref{l:distance-in-Enu}(i), we have $c_2^{-1} x_dy_d\le \delta_D(x) \delta_D(y) \le c_2 x_dy_d$ for some constant independent of $Q,x,y$. Thus, by \eqref{e:Phi-scaling-2} (with $\beta=1$),  we have
\begin{align}\label{e:difference-2}
&	\sup_{s\in I}s\Phi\bigg(\frac{|x-y|^2}{s}\bigg)  \le c_3 \delta_D(x) \delta_D(y)\Phi\bigg(\frac{|x-y|^2}{\delta_D(x)\delta_D(y)}\bigg) \nn\\
&\le c_4 (|x-y|^2 \vee (\delta_D(x)\delta_D(y))) \le c_4 (|x-y| \vee \delta_D(y))^2.
\end{align}
On the other hand, using Lemma \ref{l:distance-in-Enu}(i) and (ii), we get that
\begin{align}\label{e:difference-3}
&	 \frac{|x_dy_d - \delta_D(x)\delta_D(y)|}{\delta_D(x)\delta_D(y)x_dy_d}\le    \frac{|x_d - \delta_D(x)|}{\delta_D(x)\delta_D(y)x_d} +  \frac{|y_d-\delta_D(y)|}{\delta_D(y)x_dy_d}\nn\\
&\le  \frac{c_5 \delta_D(x)^{\frac{1+\eta}{1+\nu_0}}}{\delta_D(x)^2\delta_D(y)} +  \frac{c_6\delta_D(y)^{\frac{1+\eta}{1+\nu_0}}}{\delta_D(x)\delta_D(y)^2} \le \frac{(c_5+c_6) \delta_D(x)^{\frac{\eta-\nu_0}{1+\nu_0}} }{\delta_D(x)^{2}}.
\end{align}
Combining \eqref{e:difference-1} with \eqref{e:difference-2} and \eqref{e:difference-3}, we arrive at
\begin{align*}
&\big|\sB^D(x,y)-\sB^D(x,x)\F_0((y-x)/x_d)\big|+\big|\sB^D(x,y)-\sB^D(y,y)\F_0((y-x)/x_d)\big|\\
&\le  c_7 \left(\frac{|x-y|\vee  \delta_D(y)}{\delta_D(x)}\right)^{2} 
\delta_D(x)^{\frac{\eta-\nu_0}{1+\nu_0}} \le c_7 \left(\frac{|x-y|\vee  \delta_D(y)}{|x-y|\wedge \delta_D(x)}\right)^{2} 
\big(|x-y|\vee \delta_D(y)\big)^{\frac{\eta-\nu_0}{1+\nu_0}}. 
\end{align*}
Hence, \hyperlink{B2}{{\bf (B2)}}  
holds with $\theta_1=2$ and $\theta_2=( \eta-\nu_0 )/(1+\nu_0)>0$.

%%%%%%%%%%%%%%%%%%%%%%%%%%%%%%%%%%%%%%%%%%%%%%%%%%%
\vspace{.1in}
\begin{singlespace}

\end{singlespace}

\end{document}